\documentclass[12pt]{article}
\usepackage{latexsym}
\usepackage{graphics}
\DeclareGraphicsExtensions{eps,ps,pdf}
\usepackage{amssymb}
\usepackage[mathscr]{eucal}
\input{amssym.def}
\input{amssym.tex}

\newtheorem{theorem}{Theorem}
\newtheorem{lemma}{Lemma}
\newtheorem{proposition}{Proposition}
\newtheorem{remark}{Remark}

\newtheorem{corollary}{Corollary}
\newtheorem{definition}{Definition}

\newcommand{\Aut}{{\rm Aut}}
\newcommand{\aut}{{\rm Aut}}

\def\bigDiamond{\mathop{\hbox{{\Large $\mathop{\Diamond}$}}}}
\def\blim{\mathop{{\bf lim}}}

\newcommand{\Dec}{\mathscr{D}\mathrm{ec}}
\newcommand{\dec}{{\rm Dec}}

\newcommand{\DecA}{\Dec_{\A}}

\newcommand{\decA}{\dec_{\A}}
\newcommand{\decAS}{\dec_{\AS}}

\newcommand{\deci}{\dec_{\infty}}

\newcommand{\decv}{\dec_{v}}

\newcommand{\End}{{\rm End}}

\newcommand{\GL}{{\rm GL}}

\newcommand{\GLn}{{\rm GL}_n}

\newcommand{\GLnAS}{{\rm GL}_n(\AS)}

\newcommand{\Hom}{{\rm Hom}}

\newcommand{\ord}{{\rm ord}}

\newcommand{\proof}{{\it Proof. }}
\newcommand{\pr}{{\rm pr}}

\newcommand{\rep}{{\rm Rep}}

\newcommand{\rk}{{\rm rk}}

\newcommand{\SL}{{\rm SL}}

\newcommand{\St}{{\rm St}}

\newcommand{\Sp}{{\rm Sp}}
\newcommand{\sign}{{\rm sign}}

\newcommand{\tr}{{\rm Tr}}

\newcommand{\bi}{{\bf i}}
\newcommand{\bj}{{\bf j}}

\newcommand{\cC}{{\mathcal C}}

\newcommand{\cCo}{{\mathcal C}^{o}}

\newcommand{\cF}{{\mathcal F}}

\newcommand{\cFo}{{\mathcal F}^{o}}

\newcommand{\cI}{{\mathcal I}}

\newcommand{\cM}{{\mathcal M}}

\def\path{{\wp}}
\def\scA{{\mathscr A}}
\def\scB{{\mathscr B}}
\def\scC{{\mathscr C}}
\def\scD{{\mathscr D}}
\def\scE{{\mathscr E}}
\def\scF{{\mathscr F}}
\def\scG{{\mathscr G}}
\def\scH{{\mathscr H}}

\def\scI{{\mathscr I}}
\def\scM{{\mathscr M}}
\def\scN{{\mathscr N}}
\def\scP{{\mathscr P}}
\def\scPt{\tilde{\mathscr P}}
\def\scQ{{\mathscr Q}}

\def\scS{{\mathscr S}}
\def\scT{{\mathscr T}}
\def\scU{{\mathscr U}}
\def\scV{{\mathscr V}}
\def\scW{{\mathscr W}}

\def\sq{\hfill$\Box$}

\newcommand{\A}{{\mathbb A}}
\newcommand{\Af}{{\mathbb A}_{f}}
\newcommand{\AS}{{\mathbb A}(S)}
\newcommand{\C}{{\mathbb C}}

\newcommand{\N}{{\mathbb N}}

\newcommand{\Q}{{\mathbb Q}}

\newcommand{\R}{{\mathbb R}}
\newcommand{\Z}{{\mathbb Z}}

\def\ga{{\mathfrak a}}

\def\gf{{\mathfrak f}}
\def\gF{{\mathfrak F}}
\def\gg{{\mathfrak g}}

\def\go{{\mathfrak o}}

\def\gp{{\mathfrak p}}

\def\gs{{\mathfrak s}}
\def\gS{{\mathfrak S}}
\def\gt{{\mathfrak t}}

\def\gX{{\mathfrak X}}

\def\xx{\underline{x}}
\def\yy{\underline{y}}
\def\para{\scP{\rm ar}}
\def\parat{\widetilde{\scP}{\rm ar}}
\def\Mu{\overline\mu}

\def\deca{{\rm Dec}_{\rm arith}}
\def\decm{{\rm Dec}_{m}}
\def\Ga{G_{\rm arith}}
\def\gal{{\rm Gal}}
\def\GLnA{{\rm GL}_{n}(\A)}
\def\GLt{\widetilde{\rm GL}}
\def\SLt{\widetilde{\rm SL}}

\def\Va{V_{\rm arith}}
\def\Vert{{\rm Vert}}

\def\Span{{\rm span}}

\begin{document}

\title{Geometric construction of metaplectic covers of $\GL_{n}$
       in characteristic zero}
\author{Richard Hill}
\date{March 2002}
\maketitle

\begin{abstract}
\noindent
This paper presents a new construction
 of the $m$-fold metaplectic cover of $\GL_{n}$ over
 an algebraic number field $k$, where $k$ contains a primitive $m$-th
 root of unity.
A 2-cocycle on $\GL_{n}(\A)$ representing
 this extension is given and the splitting of the 
 cocycle on $\GL_{n}(k)$ is found explicitly.
The cocycle is smooth at almost all places of $k$.
As a consequence, a formula for the
 Kubota symbol on $\SL_{n}$ is obtained.
The construction of the paper
 requires neither class field theory
 nor algebraic K-theory, but relies instead on
 naive techniques from the geometry of numbers introduced
 by W. Habicht and T. Kubota.
The power reciprocity law for a number field
 is obtained as a corollary.
\end{abstract}

\tableofcontents

\section{Introduction.}

\subsection{Metaplectic Groups}

Let $k$ be a global field with ad\`ele ring $\A$
 and let $G$ be a linear algebraic group over $k$.
We shall regard $G(\A)$ as a locally compact topological group with
 the topology induced by that of $\A$.
For a finite abelian group $A$, a \emph{metaplectic extension}
 of $G$ by $A$ is a topological central extension:
$$
 1 \to A \to \tilde{G}(\A) \to G(\A) \to 1,
$$
 which splits over the discrete subgroup $G(k)$ of $G(\A)$.
Such extensions are of use in the theory of
 automorphic forms since certain automorphic forms
 (for example the classical theta-functions,
 see \cite{weil})
 may be regarded as functions on $\tilde{G}(\A)$,
 which are invariant under translation by the lift of
 $G(k)$ to $\tilde{G}(\A)$.
The groups $\tilde{G}(\A)$ are topological groups
 but they are not in general groups of ad\`ele valued points of
 an algebraic group.

If $k$ contains a primitive $m$-th root of unity,
 then the group $\SL_{n}$ has a canonical metaplectic
 extension with kernel the group $\mu_{m}$ of all $m$-th
 roots of unity in $k$.
This extension is always non-trivial;
 in fact if $m$ is the total number of roots of unity in $k$
 and if $n\ge 3$,
 then the canonical extension
 is universal amongst metaplectic extensions.

As the groups $\GL_{n}(\A)$ and $\GL_{n}(k)$
 are not perfect, $\GL_{n}$ has no universal
 metaplectic extension.
However the canonical extension of $\SL_{n}$
 may be continued in various ways to give a metaplectic
 extension of $\GL_{n}$ by $\mu_{m}$.
This has been done by embedding $\GL_{n}$ in
 $\SL_{r}$ for $r>n$ (see \cite{kazpat});
 we shall call the metaplectic extensions on $\GL_{n}$
 obtained in this way the \emph{standard twists}.
In this paper we shall give an elementary construction
 of a metaplectic extension $\GLt_{n}(\A)$
 of $\GL_{n}$, which in fact is not one of the standard twists,
 but which nevertheless restricts to the canonical metaplectic
 extension of $\SL_{n}$.
The method used gives an independent construction
 of the canonical metaplectic extension of $\SL_{n}$.

There are other ways of constructing the group
 $\GLt_{n}(\A)$.
The advantages of the method of construction employed in this
 paper are as follows:
\begin{itemize}
    \item
    Other methods of construction require class field theory and
     algebraic K-theory.
    In contrast the method here is very elementary.
    In fact one can deduce certain theorems of class field theory
     as corollaries of the results here.

    \item
    The method described here is very explicit in the sense that
     a 2-cocycle $\decA$ is given which represents the group extension.
    This means $\GLt_{n}(\A)$ may be realized as a set of pairs
     $(\alpha,\chi)\in \GL_{n}(\A)\times\mu_{m}$ with multiplication
     given by
    $$
     (\alpha,\chi)(\beta,\xi)
     =
     (\alpha\beta,\chi\xi \decA(\alpha,\beta)).
    $$
    Thus the cohomology class of $\decA$ is an element of
     $H^{2}(\GL_{n}(\A),\mu_{m})$ which splits on the subgroup
     $\GL_{n}(k)$.
    An expression for $\decA$ as a coboundary on $\GL_{n}(k)$
     is also obtained.
    In contrast the usual method of construction gives only
     a cocycle on the standard Borel subgroup.
    An expression for the whole cocycle has been obtained
     (after various incorrect formulae obtained by
     other authors) in \cite{banksetc}, but the cocycle
     obtained there is more complicated than ours.
    In partiular the formula of \cite{banksetc}
     involves first decomposing
     $\alpha$, $\beta$ and $\alpha\beta$ in the Bruhat
     decomposition, and then decomposing each of the three
     Weyl group elements as a minimal product of simple
     reflections.

    \item
    The cocycle $\decA$ is smooth on the non-archimedean part of
     $\GL_{n}(\A)$; in fact if $k$ has no real places
     then we obtain a cocycle which is smooth everywhere.
    The cocycle may therefore be used to study the smooth
     representations of the metaplectic group.
    More precisely suppose $\pi$ is an
     irreducible representation of $\GLt_{n}(\A)$
     on a space $V$ of smooth functions
     on $\GLt_{n}(\A)$, on which $\GLt_{n}(\A)$ acts by
     right translation:
    \begin{equation}
	\label{actionA}
	(\pi(g)\phi)(h):= \phi(hg),\quad
	g,h\in\GLt_{n}(\A).
    \end{equation}
    Let $\epsilon:\mu_{m}\to\C^{\times}$ be the restriction of the
     central character to the subgroup $\mu_{m}$.
    Then $V$ is isomorphic to a space $V'$ of smooth functions on
     $\GL_{n}(\A)$ with the twisted action:
    \begin{equation}
	\label{actionB}
	(\pi(\alpha,\chi) \phi')(\beta)
	=
	\epsilon(\chi \decA(\beta,\alpha)) \phi'(\beta\alpha).
    \end{equation}
    The isomorphism $V\to V'$ is given by $\phi\mapsto\phi'$,
     where $\phi'$ is defined by
    $$
     \phi'(\alpha)
     =
     \phi(\alpha,1).
    $$
    Although the action (\ref{actionB}) looks more complicated
     than (\ref{actionA}), it is perhaps easier to
     use in calculations as one is dealing with elements of $\GL_{n}$.
    One cannot do this with the cocycle of \cite{banksetc}
     as it is not smooth (in fact on $\GL_{n}(\A)$ with $n\ge 2$
     it is nowhere continuous).
\end{itemize}
There are two disadvantages to the method of construction
 described in this paper.
First, the construction is long and quite difficult.
Second, the part of the
 cocycle on the subgroup $\GL_{n}(k_{v})$ for $v|m$ is
 not very explicit.
In a sense one has the same problem with the cocycle
 of \cite{banksetc}, since it is expressed in
 terms of ramified Hilbert symbols.

\paragraph{The method of construction.}
The case that $k$ is a function field is described in
 \cite{hill2}, \cite{hill3};
 in this paper we shall deal with the more difficult case that $k$
 is an algebraic number field.
Let $S$ be the set of places $v$ of $k$ for which $|m|_{v}\ne 1$,
 and let $\AS$ denote the restricted topological product of the
 fields $k_{v}$ for $v\notin S$.
Let $k_{\infty}$ be the sum of the archimedean completions of $k$
 and let $k_{m}$ be the sum of the fields $k_{v}$
 for non-archimedean places $v\in S$.
We then have
$$
 \A
 =
 \AS \oplus k_{\infty} \oplus k_{m},
$$
 and hence:
$$
 \GL_{n}(\A)
 =
 \GL_{n}(\AS) \oplus \GL_{n}(k_{\infty}) \oplus \GL_{n}(k_{m}).
$$
We shall write down explicit 2-cocycles $\decAS$ on $\GL_{n}(\AS)$
 and $\deci$ on $\GL_{n}(k_{\infty})$.
Then, for a certain compact, open subgroup $U_{m}$ of $\GL_{n}(k_{m})$,
 we find a function $\tau:\GL_{n}(k)\cap U_{m}\to \mu_{m}$ such that
\begin{equation}
    \label{ratsplit}
    \decAS(\alpha,\beta)
    \deci(\alpha,\beta)
    =
    \frac{\tau(\alpha)\tau(\beta)}{\tau(\alpha\beta)},
    \quad
    \alpha,\beta\in\GL_{n}(k)\cap U_{m}.
\end{equation}
It follows fairly easily from this that we can extend $\tau$ to $\SL_{n}(k)$
 in such a way that there is a unique continuous cocycle $\dec_{m}$ on
 $\SL_{n}(k_{m})$ defined by the formula
$$
 \decAS(\alpha,\beta)
 \deci(\alpha,\beta)
 \dec_{m}(\alpha,\beta)
 =
 \frac{\tau(\alpha)\tau(\beta)}{\tau(\alpha\beta)},
 \quad
 \alpha,\beta\in\SL_{n}(k).
$$
Finally we show that there is a cocycle $\decA$ on $\GL_{n}(\A)$
 which is metaplectic and which extends all our cocycles.

Note that this definition of $\dec_{m}$ is global;
 it is defined on the dense subgroup $\SL_{n}(k)$ and then
 extended by continuity to $\SL_{n}(k_{m})$.
It would be of some interest to find a local construction
 of the cocycle $\dec_{m}$, as the ramified Hilbert
 symbols may be expressed in terms of this cocycle
 via the isomorphism (\ref{K2}) below.
However I do not know how to make such a construction.

If $m$ is even then the cocycle $\decA$
 has the surprising property that it is not, even up to a coboundary,
 a product of cocycles $\decv$ on the groups $\GL_{n}(k_{v})$
 (in contrast the standard twists are products of local cocycles).
In fact if we write $\GLt_{n}(k_{v})$ for the preimage of
 $\GL_{n}(k_{v})$ in $\GLt_{n}(\A)$, then the various
 subgroups $\GLt_{n}(k_{v})$ do not even commute with each other.
This means that irreducible representations $\pi$
 of $\GLt_{n}(\A)$ cannot be expressed as
 restricted tensor products of
 irreducible representations of the groups $\GLt_{n}(k_{v})$.
Thus the usual local-to-global approach
 to studying automorphic representations must be
 modified to deal with $\GLt_{n}(\A)$.

\paragraph{Matsumoto's Construction.}
We now review the usual construction of the canonical metaplectic
 extension of $\SL_{n}$.
Let $F$ be any field.
Recall (or see \cite{milnor} or \cite{hahnomeara})
 that for $n\ge 3$ there is a universal central extension
$$
 1 \to K_{2}(F) \to \St_{n}(F) \to \SL_{n}(F) \to 1,
$$
 where $\St_{n}$ denotes the Steinberg group.
Hence, for any abelian group $A$ we have
\begin{equation}
    \label{K2}
    H^{2}(\SL_{n}(F),A) \cong \Hom(K_{2}(F),A).
\end{equation}
Matsumoto (see \cite{matsumoto} or \cite{milnor})
 proved that for any field $F$
 the group $K_{2}(F)$ is the quotient of
 $F^{\times}\otimes_{\Z} F^{\times}$ by the subgroup generated
 by $\{\alpha\otimes (1-\alpha) :\alpha\in F\setminus\{0,1\}\}$.
The isomorphism (\ref{K2}) may be described as follows.
If $\sigma$ is a 2-cocycle representing a cohomology class
 in $H^{2}(\SL_{n}(F),A)$ then for diagonal matrices
 $\alpha,\beta\in\SL_{n}(F)$ we have
$$
 \frac{\sigma(\alpha,\beta)}{\sigma(\beta,\alpha)}
 =
 \prod_{i=1}^{n}\phi(\alpha_{i},\beta_{i}),
 \quad
 \alpha
 =
 \left(
 \matrix{\alpha_{1}\cr &\ddots\cr &&\alpha_{n}}
 \right),
 \quad
 \beta
 =
 \left(
 \matrix{\beta_{1}\cr &\ddots\cr &&\beta_{n}}
 \right),
$$
 where $\phi:K_{2}(F)\to A$ is the image of $\sigma$.
Here we are writing the group law in $A$ multiplicatively.

Suppose that $k$ is a global field
 containing a primitive $m$-th root of unity.
For any place $v$ of $k$ the $m$-th power
 Hilbert symbol gives a map $K_{2}(k_{v})\to \mu_{m}$.
Corresponding to this map there is a cocycle
 $\sigma_{v}\in H^{2}(\SL_{n}(k_{v}),\mu_{m})$.
For any place $v$ of $k$ we shall write $\go_{v}$ for the ring of
 integers in $k_{v}$.
For almost all places $v$ the cocycle $\sigma_{v}$ splits on
 $\SL_{n}(\go_{v})$.
One may therefore define a cocycle $\sigma_{\A}$ on $\SL_{n}(\A)$
 by $\sigma_{\A}=\prod_{v}\sigma_{v}$.
This corresponds to a topological central extension:
\begin{equation}
    \label{extension}
    1 \to \mu_{m} \to \SLt_{n}(\A) \to \SL_{n}(\A) \to 1.
\end{equation}
Now recall the following.

\begin{theorem}
    [Power Reciprocity Law]
    For any place $v$ of $k$ let $(-,-)_{v,m}$ denote the
     $m$-th power Hilbert symbol on $k_{v}$.
    For $\alpha,\beta\in k^{\times}$ we have
    $$
     \prod_{v} (\alpha,\beta)_{v,m}
     =
     1,
    $$
    where the product is taken over all places of $k$.
\end{theorem}

(For a proof, see Chapter 12, Verse 4, Theorem 12 of \cite{ArtinTate}.)
    
If one restricts $\sigma_{\A}$ to $\SL_{n}(k)$,
 this restriction corresponds under (\ref{K2})
 to the map $K_{2}(k)\to \mu_{m}$
 induced by the $m$-th power Hilbert symbol on $k_{v}$.
Hence the restriction of $\sigma_{\A}$ to $\SL_{n}(k)$
 corresponds to the map $K_{2}(k)\to \mu_{m}$ induced by the
 product of all the $m$-th power Hilbert symbols.
By the reciprocity law this map is trivial.
Therefore $\sigma_{\A}$ splits on $\SL_{n}(k)$,
 so the extension (\ref{extension}) is metaplectic.

\subsection{The Kubota symbol}

One of the results of this paper is a formula
 (see \S6.2) for the Kubota symbol on $\SL_{n}$.
We recall here the definition of the Kubota symbol.

Let $k$ be an algebraic number field and let
 $\go$ denote the ring of algebraic integers in $k$.
Given an ideal $\ga$ of $\go$,
 we define $\SL_{n}(\go,\ga)$ to be the subgroup of
 matrices in $\SL_{n}(\go)$ which are congruent
 to the identity matrix modulo $\ga$.
A subgroup of $\SL_{n}(k)$ is said to be an \emph{arithmetic} subgroup
 if it is commensurable with $\SL_{n}(\go)$.
An arithmetic subgroup is said to be a \emph{congruence subgroup}
 if it contains $\SL_{n}(\go,\ga)$ for some non-zero ideal $\ga$.
The \emph{congruence subgroup problem} is the question:
 ``is every arithmetic subgroup a congruence subgroup?''
This question has been answered by Bass--Milnor--Serre
 \cite{bassmilnorserre} for $n\ge 3$ and by Serre \cite{serre2} in the
 more difficult case $n=2$
 (see for example \cite{PrasadRaghunathan} for generalizations).
To rule out some pathological cases,
 assume either that $k$ has at least two archimedean places
 or that $n\ge 3$.
If $k$ has a real place then the answer to the congruence
 subgroup problem is ``yes'' whereas if $k$ is totally complex
 the answer is ``no''.
We shall describe what happens when $k$ is totally complex.
To make statements clearer suppose $\mu_{m}$ is the group of all roots of
 unity in $k$.
In this case there is an ideal $\gf$ and a surjective homomorphism
 $\kappa_{m}:\Gamma(\gf)\to \mu_{m}$ with the following properties:
\begin{itemize}
    \item
    $\ker(\kappa_{m})$ is a non-congruence subgroup.
    \item
    For any arithmetic subgroup $\Gamma$ of $\SL_{n}(k)$ there
     is an ideal $\ga$
     such that $\Gamma\supset\SL_{n}(\go,\ga)\cap\ker(\kappa_{m})$.
\end{itemize}
The map $\kappa_{m}$ is called the Kubota symbol.
In the case of $\SL_{2}$ it is given by the following formula
 (see \cite{kubota1}):
$$
 \kappa_{m}\left(\matrix{a & b \cr c & d}\right)
 =
 \left\{
 \begin{array}{ll}
     \displaystyle{
     \left(\frac{c}{d}\right)_{m}} &
     \hbox{if $c\ne 0$},\medskip \\
     1 & \hbox{if $c=0$.}
 \end{array}
 \right.
$$
Here $\displaystyle{\left(\frac{c}{d}\right)_{m}}$ denotes the $m$-th power
 residue symbol, which we recall is defined by
$$
 \left(\frac{c}{d}\right)_{m}
 =
 \prod_{v|d} (c,d)_{v,m}
 =
 \left[\frac{k(\sqrt[m]{c})/k}{d}\right].
$$
The symbol on the right is the global Artin symbol;
 the Galois group $\gal(k(\sqrt[m]{c})/k)$ may be
 identified with $\mu_{m}$.

\paragraph{Connection between the congruence subgroup problem
 and metaplectic groups.}
As before suppose the number field $k$ is totally complex.
We shall introduce two topologies on $\SL_{n}(k)$.
For the first topology we take the congruence subgroups
 as a basis of neighbourhoods of the identity.
The completion of $\SL_{n}(k)$ with respect to this topology
 is the group $\SL_{n}(\Af)$, where $\Af$ denotes the ring of finite
 ad\`eles of $k$.
For our second topology we take the arithmetic subgroups
 of $\SL_{n}(k)$ as a basis of neighbourhoods of the identity.
This is a finer topology, and the completion of $\SL_{n}(k)$ with
 respect to this topology is a group extension $\SLt_{n}(\Af)$ of
 $\SL_{n}(\Af)$.
The homomorphism $\kappa_{m}$ identifies the kernel of the extension
 with $\mu_{m}$.
We therefore have a short exact sequence:
\begin{equation}
    \label{SLtAf}
    1 \to \mu_{m} \to \SLt_{n}(\Af) \to \SL_{n}(\Af) \to 1.
\end{equation}
This must be a central extension as $\SL_{n}(\Af)$ is perfect.
On the other hand $\SL_{n}(k)$ is contained in both completions,
 so the extension splits on $\SL_{n}(k)$.
Adding the group $\SL_{n}(k_{\infty})$
 to both $\SL_{n}(\Af)$ and $\SLt_{n}(\Af)$
 we obtain the canonical metaplectic extension of $\SL_{n}$.

Conversely we may reconstruct the Kubota symbol from the
 metaplectic group as follows.
Let $\SLt_{n}(\Af)$ be the
 preimage of $\SL_{n}(\Af)$ in the
 canonical metaplectic extension of $\SL_{n}$.
By restriction we have an exact sequence (\ref{SLtAf}).
Since $k$ is totally complex,
 the group $\SL_{n}(k_{\infty})$ is
 both connected and simply connected.
This implies that the restriction map gives an
 isomorphism
$$
 H^{2}(\SL_{n}(\A),\mu_{m})
 \to
 H^{2}(\SL_{n}(\Af),\mu_{m}).
$$
Hence, if we regard $\SL_{n}(k)$
 as a subgroup of $\SL_{n}(\Af)$,
 then this subgroup lifts to a subgroup of $\SLt_{n}(\Af)$.

There are two subgroups of
 $\SL_{n}(\Af)$ on which
 the extension (\ref{SLtAf})
 splits.
First, since (\ref{SLtAf}) is a topological central extension,
 there are neighbourhoods $\tilde U$ of the
 identity in $\SLt_{n}(\Af)$
 and $U$ of the identity in $\SL_{n}(\Af)$
 such that the projection map restricts to
 a homeomorphism $\tilde U\to U$.
As $\SL_{n}(\Af)$ is totally disconnected
 we may take $\tilde U$, and hence $U$,
 to be a compact open subgroup.
The inverse map $U\to \tilde U$ is a
 continuous splitting of the extension.
Secondly we have a splitting of the extension on $\SL_{n}(k)$.
By the Strong Approximation Theorem
 (Theorem 3.3.1 of \cite{bump}),
 $\SL_{n}(k)$ is dense in $\SL_{n}(\Af)$.
Hence the splitting on $\SL_{n}(k)$ cannot be continuous,
 since otherwise it would extend by continuity to a splitting
 of the whole extension.

Let $\Gamma=U\cap \SL_{n}(k)$.
The group $\Gamma$ is a congruence subgroup of $\SL_{n}(k)$,
 and we have two different splittings of the extension on $\Gamma$.
Dividing one splitting by the other we obtain a homomorphism
 $\kappa_{m}:\Gamma\to\mu_{m}$.
This homomorphism is not continuous with respect to the
 induced topology from $\SL_{n}(\Af)$,
 so its kernel is a non-congruence subgroup.
The map $\kappa_{m}$ is the Kubota symbol.

In the construction of this paper the
 splittings of the extension are described explicitely.
As a consequence we obtain a formula for
 the Kubota symbol on $\SL_{n}$.

\subsection{Organization of the paper}

The paper is organized into the following sections:
\begin{itemize}
    \item[\S2]
    We fix some standard notation from group cohomology and
     singular homology.
    To avoid confusions of signs later, the normalizations
     of various maps are fixed.
    Some known results are stated for later reference.
    \item[\S3]
    The cocycles $\decAS$ and $\deci$ on the
     groups $\GL_{n}(\AS)$ and $\GL_{n}(k_{\infty})$
     are defined.
    The cocycle $\decAS$ has been studied in \cite{hill2};
     some of the results from there are recalled.
    Analogous results are obtained for the cocycle $\deci$.
    On their own the results of this section
     concerning $\deci$ are of little interest since they
     describe group extensions which are already
     well understood.
    It is the relation (\ref{ratsplit}) between $\deci$ and $\decAS$
     that is interesting.
    This relation  is stated in \S4 and proved in \S6.
    \item[\S4]
    The function $\tau$ is defined, and the formalism
     used in the proof of the relation (\ref{ratsplit})
     between $\decAS$, $\deci$ and $\tau$ is introduced.
    \item[\S5]
    This is a technical section on the existence of
     certain limits.
    \item[\S6]
    The relation (\ref{ratsplit}) between $\decAS$,
     $\deci$ and $\tau$ is proved.
    \item[\S7]
    The cocycles are extended to $\GL_{n}(\A)$.
\end{itemize}
Sections 3 and 4 are easy, although
 section 4 has a lot of notation.
Section 7 is quite formal.
In contrast, sections 5 and 6 are difficult.

\section{Notation}

\subsection{Acyclic $(\Z/m)[\mu_{m}]$-modules.}
Throughout the paper, $\mu_{m}$ will denote a cyclic group of order $m$.
As $\mu_{m}$ will often be taken to be a group of roots of unity,
 we shall write the group law in $\mu_{m}$ multiplicatively.
We shall deal with various modules over the group-ring $(\Z/m)[\mu_{m}]$.
We shall write $[\mu_{m}]$ for the sum of the elements of $\mu_{m}$.
We also fix once and for all a generator $\zeta$ of $\mu_{m}$.
There is an exact sequence (see \S7 of \cite{AtiyahWall}):
$$
 (\Z/m)[\mu_{m}]
 \stackrel{[\mu_{m}]}{\leftarrow}
 (\Z/m)[\mu_{m}]
 \stackrel{1-[\zeta]}{\leftarrow}
 (\Z/m)[\mu_{m}]
 \stackrel{[\mu_{m}]}{\leftarrow}
 (\Z/m)[\mu_{m}].
$$
Applying the function $\Hom_{(\Z/m)[\mu_{m}]}(-,M)$
 we obtain the chain complex:
\begin{equation}
    \label{star}
    M
    \stackrel{[\mu_{m}]}{\to}
    M
    \stackrel{1-[\zeta]}{\to}
    M
    \stackrel{[\mu_{m}]}{\to}
    M.
\end{equation}
The cohomology of this complex is the Tate cohomology
 $\hat H^{\bullet}(\mu_{m},M)$ (see \cite{weibel}).
We shall call $M$ an \emph{acyclic} module if (\ref{star})
 is exact, i.e. if its Tate cohomology is trivial.

Free modules are clearly acyclic.
Injective modules are acyclic,
 since for these the functor
 $\Hom(-,M)$ is exact.
More generally, by Shapiro's Lemma
 (see \cite{weibel}),
 any $(\Z/m)[\mu_{m}]$-module
 which is induced from a $\Z/m$-module
 is acyclic.

\begin{lemma}
    \label{easy}
    Let $M$ be an acyclic $(\Z/m)[\mu_{m}]$-module
     and suppose we have a $(\Z/m)[\mu_{m}]$-module homomorphism
    $$
     \Phi:M \to \Z/m.
    $$
    Then there is a $(\Z/m)[\mu_{m}]$-homomorphism
     $\hat\Phi:(1-[\zeta])M\to \mu_{m}$
     defined by
    $$
     \hat\Phi ( (1-[\zeta]) a ) = \zeta^{\Phi(a)}.
    $$
    The lifted map $\hat\Phi$ is independent of the choice of 
     generator $\zeta$.
    Here $\Z/m$ and $\mu_{m}$ are regarded as
     $(\Z/m)[\mu_{m}]$-modules with the trivial action of $\mu_{m}$.
\end{lemma}

The definition of $\hat\Phi$ looks more natural if we identify
 $\mu_{m}$ with $\ga/\ga^{2}$, where $\ga$ denotes the
 augmentation ideal of $(\Z/m)[\mu_{m}]$.
\medskip

\proof
We need only show that $\hat\Phi$ is well-defined.
Suppose $(1-[\zeta]) a=(1-[\zeta])b$.
By the exact sequence (\ref{star}), there is a $c\in M$
 such that $b-a=[\mu_{m}] c$.
Therefore $\Phi(b)-\Phi(a)=[\mu_{m}] \Phi(c)$.
Since the action of $[\mu_{m}]$ on $\Z/m$ is zero,
 we have $\Phi(b)=\Phi(a)$.
\sq
\medskip

\subsection{Central extensions.}
Let $G$ be an abstract group.
We shall regard $\mu_{m}$ as a $G$-module with the trivial action.
Given an inhomogeneous 2-cocycle $\sigma$
 on $G$ with values in $\mu_{m}$, one defines a central
 extension of $G$ by $\mu_{m}$ normalized as follows:
$$
 \tilde G
 =
 G\times \mu_{m},\quad
 (g,\xi)(h,\psi)
 :=
 (gh,\xi\psi\sigma(g,h)).
$$
Conversely, given a central extension
$$
 1 \to \mu_{m} \to \tilde G \to G \to 1,
$$
 we may recover a 2-cocycle $\sigma$ by choosing,
 for every $g\in G$,
 a preimage $\hat g\in\tilde G$
 and defining
$$
 \sigma(g,h)
 =
 \hat{g}\hat{h}\widehat{gh}^{-1}.
$$

If $G$ is a locally compact topological group
 then by a \emph{topological central extension}
 of $G$ by $\mu_{m}$ we shall mean a short exact sequence of
 topological groups and continuous homomorphisms:
$$
 1 \to \mu_{m} \to \tilde{G} \to G \to 1,
$$
 such that the map $\tilde{G}\to G$ is locally
 a homeomorphism.
The isomorphism classes of
 topological central extensions of $G$ by $\mu_{m}$
 correspond to elements of $H^{2}_{\rm meas}(G,\mu_{m})$,
 where $H^{\bullet}_{\rm meas}$ denotes
 the group cohomology theory based on Borel-measurable cochains
 (see \cite{moore}).
As all our cochains will be Borel-measurable
 we shall write $H^{\bullet}$ instead of
 $H^{\bullet}_{\rm meas}$.

\subsection{Commutators and symmetric cocycles}

Proofs of the following facts may be found in \cite{klose},
 where they were used to determine the metaplectic extensions of
 $D^{\times}$, where $D$ is a quaternion algebra over $k$.

Let $G$ be a group and suppose
 $G_{1}$ and $G_{2}$ are two subgroups of $G$,
 such that every element of $G_{1}$ commutes with every element
 of $G_{2}$.
Given a 2-cocycle $\sigma\in Z^{2}(G,\mu_{m})$,
 we define for $a\in G_{1}$ and $b\in G_{2}$
 the \emph{commutator}:
$$
 [a,b]_{\sigma}
 :=
 \frac{\sigma(a,b)}{\sigma(b,a)}.
$$
The commutator map is bimultiplicative and skew symmetric,
 and depends only on the cohomology class of $\sigma$.
If $G$ is a locally compact topological group
 and $\sigma$ is Borel-measurable,
 then the commutator map is a continuous function on $G_{1}\times G_{2}$.

If $G$ is an abelian group,
 then we shall call a 2-cocycle $\sigma$
 on $G$ \emph{symmetric} if $[\cdot,\cdot]_{\sigma}$ is trivial on
 $G\times G$.
This amounts to saying that the corresponding central extension
 $\tilde G$ is abelian.
We shall write $H^{2}_{sym}(G,\mu_{m})$ for the subgroup
 of symmetric classes.
If $G$ and $H$ are two abelian groups, then
 the restriction maps give an isomorphism:
$$
 H^{2}_{sym}(G\oplus H,\mu_{m})
 \cong
 H^{2}_{sym}(G,\mu_{m})
 \oplus
 H^{2}_{sym}(H,\mu_{m}).
$$
The restriction map gives an isomorphism:
$$
 H^{2}_{sym}(G,\mu_{m})
 \cong
 H^{2}_{sym}(G[m],\mu_{m}),
$$
 where $G[m]$ denotes the subgroup of $m$-torsion elements of $G$.
Furthermore there is a canonical isomorphism (independent of $\zeta$):
$$
 \begin{array}{rcl}
     H^{2}_{sym}(\mu_{m},\mu_{m})
     &\cong&
     \Z/m\\
     \sigma
     &\mapsto&
     b,\quad \hbox{ where }\quad
     \displaystyle{
     \prod_{i=1}^{m}\sigma(\zeta,\zeta^{i})
     =
     \zeta^{b}.
     }
 \end{array}
$$

\subsection{Singular homology groups.}

We will need some notation from singular homology,
 which we now introduce.
In order to avoid sign errors we must be
 clear about the precise definition of our chain complex,
 which is rather non-standard.

For $r\ge 0$ we define the \emph{$r$-simplex} $\Delta^r$ by
$$
 \Delta^r
 =
 \left\{
 \xx\in\R^{r+1} :
 x_0,\ldots,x_r\ge 0 \hbox{ and } \sum x_i =1
 \right\}.
$$
By an \emph{$r$-prism} we shall mean a product of finitely many simplices,
 i.e. an expression of the form $\prod_{i=1}^s \Delta^{a(i)}$,
 where $\sum a(i)=r$.
Let $Y$ be a topological space.
By a \emph{singular $r$-cell} in $Y$ we shall mean a continuous map from an
 $r$-prism to $Y$.
The cell $\scT:\prod_{i=1}^s \Delta^{a(i)}\to Y$
 will be said to be \emph{degenerate} if $\scT(\xx_1,\ldots,\xx_s)$
 is independent of one of the variables $\xx_i\in\Delta^{a(i)}$
 ($a(i)>0$).
We shall write $C_r(Y)$ for the $\Z/m\Z$-module
 generated by the set of all singular $r$-cells in $Y$,
 with the following relations:
\begin{itemize}
    \item
    Suppose $A$ is an $a$-prism and $B$ is a $b$-prism with $a+b=r$.
    If $\scT:A\times B \to Y$ is a singular $r$-cell
     then we define another singular $r$-cell
     $\scT^\prime:B\times A\to Y$ by $\scT^\prime(b,a)=\scT(a,b)$.
    For every such $\scT,\scT^\prime$ we impose a relation:
    $$
     \scT  = (-1)^{ab} \scT^\prime;
    $$
    \item
    $\scT=0$ for every degenerate $r$-prism $\scT$.
\end{itemize}
For any simplex $\Delta^r$ we define the $j$-th face map ($j=0,\ldots,r$)
 to be the map
$$
 \gF_j:\Delta^{r-1}\to \Delta^r,\quad
 (x_0,\ldots,x_{r-1})
 \mapsto
 (x_0,\ldots,x_{j-1},0,x_{j+1},\ldots,x_r).
$$
The boundary of a singular cell
 $\scT:\prod_{i=1}^s \Delta^{a(i)} \to Y$
 is defined to be the sum
\begin{eqnarray*}
 \partial \scT
 &=&
 \sum_{i=1}^s
 \sum_{j=0}^{a(i)}
 (-1)^{a(1)+\cdots+a(i-1)+j}
 \scT
 \circ
 \gF_{i,j},\\
 \hbox{where }
 \gF_{i,j}(\xx_{1},\ldots,\xx_{s})
 &=&
 (\xx_1,\ldots,\xx_{i-1},\gF_{j}(\xx_{i}),\xx_{i+1},\ldots,\xx_s).
\end{eqnarray*}
This boundary map extends by linearity
 to a map $\partial : C_r(Y) \to C_{r-1}(Y)$.

For any subspace $Z\subseteq Y$
 there is an inclusion
 $C_r(Z) \subseteq C_r(Y)$,
 and we define $C_{\bullet}(Y,Z) = C_{\bullet}(Y) / C_{\bullet}(Z)$.
The homology groups of the complexes
 $C_{\bullet}(Y)$ and $C_{\bullet}(Y,Z)$
 are the usual singular homology groups with coefficients in
 $\Z/m\Z$ (see for example \cite{massey}).
We have taken a non-standard definition of the chain complex because
 we will write down singular cells explicitly and these will
 be as described above.

The base set $|\scT|$ of a singular $r$-cell $\scT$ is defined
 to be the image of $\scT$ if $\scT$ is non-degenerate, and the empty set
 if $\scT$ is degenerate.
The base set of an element of $C_r(Y)$ is defined to be the union
 of all base-sets of singular $r$-cells in its support.

If $Y$ is a vector space over $\R$ then
 for vectors $v_{0},\ldots,v_{r}\in Y$
 we shall denote by $[v_0,\ldots,v_r]$
 the $r$-cell $\Delta^r\to Y$ given by
$$
 [v_0,\ldots,v_r](\xx)
 :=
 \sum_{i=0}^{r} x_i v_i.
$$
The image of this map is the convex hull
 of the set $\{v_{0},\ldots,v_{r}\}$.
To simplify our formulae we shall sometimes substitute
 the closed interval $I=[0,1]$ for $\Delta^1$,
 by identifying $0$ with $(1,0)$ and $1$ with $(0,1)$.

\paragraph{Orientations.}
Let $Y$ be a $d$-dimensional manifold.
If $y\in Y$ then $H_d(Y,Y\setminus \{y\})$ is non-canonically
 isomorphic to $\Z/m\Z$.
The manifold $Y$ is said to be $\Z/m\Z$-orientable
 if one can associate to each point $y\in Y$ an isomorphism
$$
 \ord_y:H_d(Y,Y\setminus \{y\}) \longrightarrow \Z/m\Z,
$$
 with the property that for every $y\in Y$ there is a neighbourhood $U$
 of $y$, such that for every $z\in U$ the following diagram commutes.

\unitlength0.8cm
\begin{picture}(15,5)
 \thicklines
 \put(5.8,4){$H_d(Y,Y\setminus U)$}
 \put(2.8,2){$H_d(Y,Y\setminus \{y\})$}
 \put(9.2,2){$H_d(Y,Y\setminus \{z\})$}
 \put(6.5,0){$\Z/m\Z$}
 \put(4.5,0.5){$\ord_y$}
 \put(9,0.5){$\ord_z$}
 \put(4.5,1.8){\vector(3,-2){2}}
 \put(10,1.8){\vector(-3,-2){2}}
 \put(8,3.8){\vector(3,-2){2}}
 \put(6.5,3.8){\vector(-3,-2){2}}
\end{picture}\\
\medskip

\noindent
Such a collection of isomorphisms will be called an \emph{orientation}.

Suppose $Y$ is $\Z/m\Z$-orientable
 and fix an orientation $\{\ord_y\}$ on $Y$.
Let $\scT \in C_d(Y)$.
Suppose that $y\in Y$ does not
 lie in the base set of $\partial \scT$.
Then $\scT$ represents a homology class in $H_{d}(Y,Y\setminus\{y\})$.
From our condition on $\ord$, the map
 $y\mapsto \ord_y(\scT)$ is a locally constant
 function $Y\setminus |\partial \scT| \to \Z/m\Z$.
For a discrete subset $M\subseteq Y$ we shall use the notation
$$
 \{ \scT | M\}
 =
 \sum_{y\in M} \ord_{y}\scT.
$$

\section{The arithmetic and geometric cocycles.}

\subsection{The arithmetic cocycle}

Let $\mu_{m}$ be a cyclic group of order $m$.
Let $\Va$ be a totally disconnected locally compact abelian group
 and suppose that multiplication by $m$ induces a measure-preserving
 automorphism of $\Va$.
Suppose $\mu_{m}$ acts on $\Va$ in such a way that
 every non-zero element of $\Va$ has trivial stabilizer in $\mu_{m}$.
Given this data we shall construct a 2-cocycle $\deca$
 on the group $\Ga=\aut_{\mu_{m}}(\Va)$ with values in $\mu_{m}$.
This cocycle $\deca$ was first studied in \cite{hill2}
 and we shall keep to the notation of that paper.
Those results, which are stated here without proof,
 are proved in \cite{hill2}.

\paragraph{The Cocycle.}
We choose a compact, open, $\mu_{m}$-invariant
 neighbourhood $L$ of $0$ in $\Va$ and we normalize
 the Haar measure $dv$ on $\Va$ so that $L$ has measure $1$.
By our condition on $\Va$, it follows that for any compact
 open subset $U$ of $\Va$, the measure of $U$ is a rational number
 and is integral at every prime dividing $m$.
Thus for any locally constant
 function $\varphi:\Va\to\Z/m$ of compact support,
 we may define its integral to be an element of $\Z/m$.
This ``modulo $m$ integration'' is independent of the neighbourhood
 $L$ used to normalize the measure.

Finally, we choose an open and closed fundamental domain
 $F$ for the action of $\mu_{m}$ on $\Va\setminus\{0\}$;
 we write $f$ for the characteristic function
 of $F$.
The cocycle is given by the formula:
\begin{equation}
    \label{def1}
    \deca^{f,L}(\alpha,\beta)
    =
    \prod_{\xi\in\mu_{m}}
    \xi^{\textstyle{\left\{\int_{L}-\int_{\beta L}\right\}
    f(\alpha v)f(\xi v)dv}},
    \qquad
    \alpha,\beta\in\Ga.
\end{equation}
Up to a coboundary, this is independent of the choices
 of $f$ and $L$.
If we regard $\Ga$ as a topological group with the
 compact-open topology then the cocycle $\deca$ is
 a locally constant function.
This may be deduced using the cocycle relation from the
 following fact.

\begin{lemma}
    If $\alpha,\beta\in \Ga$ and $\beta L=L$
     then we have $\deca^{(f,L)}(\alpha,\beta)=1$.
\end{lemma}

\proof
This follows immediately from the definition of $\deca$.
\sq
\medskip

\paragraph{A Pairing.}
Before proceeding,
 we shall reformulate the definition of $\deca$
 in a more useful notation.
We shall call a function $\varphi:\Va\to\Z/m$ \emph{symmetric}
 if $\varphi(\xi v)= \varphi(v)$ for all $\xi\in\mu_{m}$.
We define $\cC$ to be the space of locally constant
 symmetric functions of compact support on $\Va$.
On the other hand a function $g:\Va\setminus\{0\}\to\Z/m$
 will be called \emph{cosymmetric} if the sum
$$
 \sum_{\xi\in \mu_{m}} g(\xi v)
$$
 is a constant, independent of $v$.
The value of the constant will be called the \emph{degree} of $g$.
We define $\cF$ to be the space of locally constant, cosymmetric
 functions $\Va\setminus\{0\}\to\Z/m$.
A cosymmetric function $f\in\cF$ of degree $1$
 will be called a \emph{fundamental function}.
For example the function $f$ described above is a fundamental function.
Thus fundamental functions generalize
 the notion of a fundamental domain.
The group $\Ga$ acts on $\cC$ and $\cF$ on the right by composition.
Let $\cCo$ be the submodule of functions $M\in \cC$ satisfying
 $M(0)=0$ and let $\cFo$ be the submodule of functions
 in $\cF$ of degree $0$.
There is a pairing $\cFo\times \cCo\to\mu_{m}$
 given by
$$
 \langle
 g|M
 \rangle
 =
 \prod_{\xi\in\mu_{m}}
 \xi^{\textstyle{\int g(v)f(\xi v)M(v)dv}},
$$
 where $f$ is a fundamental function.
The pairing is independent of $f$,
 and is $\Ga$-invariant in the sense that for $\alpha\in \Ga$
 we always have:
\begin{equation}
    \label{covariant}
    \langle h\alpha|M\alpha \rangle
    =
    \langle h|M \rangle.
\end{equation}

With this notation we can express $\deca$ as follows:
\begin{equation}
    \label{def2}
    \deca^{(f,L)}(\alpha,\beta)
    =
    \langle f-f\alpha | \beta L-L \rangle.
\end{equation}
Here we are writing $L$ and $\beta L$ for the characteristic
 functions of these sets.
It is now clear that $\deca$ is a 2-cocycle,
 since it is the cup-product of the 1-cocycles $f-f\alpha$
 and $\beta L-L$.

\paragraph{Another expression for the pairing.}
To aid calculation we shall describe the pairing
 in a different way.
Let $\cM^{c}$ be the space of locally constant functions
 of compact support $\varphi:\Va\setminus\{0\}\to \Z/m$.
Let $\cM$ be the space of locally constant functions
 $\varphi:\Va\setminus\{0\}\to \Z/m$.
There is a right action of $\Ga$ by composition
 on the spaces $\cM^{c}$ and $\cM$:
$$
 (\phi\alpha)(v)
 :=
 \phi(\alpha v),
 \quad
 \alpha\in\Ga.
$$
We shall also regard $\cM$ and $\cM^{c}$
 as left $(\Z/m)[\mu_{m}]$-modules with the action given by
$$
 (\xi\phi)(v) = \phi(\xi^{-1}v),
 \quad
 \xi\in\mu_{m}.
$$
The two actions commute.

As $(\Z/m)[\mu_{m}]$-modules,
 both $\cM$ and $\cM^{c}$ are acyclic,
 since they are induced from spaces of
 functions on the fundamental domain $F$.
Hence, by Lemma \ref{easy},
 the integration map $\int:\cM^{c}\to \Z/m$
 lifts to a map
$$
 \widehat{\int}: (1-[\zeta])\cM^{c} \to \mu_{m},
$$
 defined by
$$
 \widehat{\int} (\varphi-\varphi\zeta^{-1})(v) dv
 =
 \widehat{\int} (1-[\zeta])\varphi dv
 :=
 \zeta^{-\textstyle{\int \varphi(v)dv}}.
$$

The modules $\cFo$ and $\cCo$ introduced above may be expressed as
 follows:
$$
 \cCo
 =
 \ker\Big((1-[\zeta]):\cM^{c}\to\cM^{c}\Big)
 =
 [\mu_{m}]\cM^{c},
$$
$$
 \cFo
 =
 \ker\Big([\mu_{m}]:\cM\to\cM\Big)
 =
 (1-[\zeta])\cM.
$$

\begin{proposition}
    Given $g\in\cFo$ and $M\in\cCo$, the product $g\cdot M$
     is in $(1-[\zeta])\cM^{c}$.
    The pairing $\cFo\times\cCo\to\mu_{m}$
     is given by
     \begin{equation}
	 \langle g|M \rangle
	 =
	 \widehat{\int} (g\cdot M)(v)dv.
     \end{equation}
\end{proposition}

\proof
Since $[\mu_{m}]g=0$, we have $g=(1-[\zeta]) h$ for some $h\in \cM$.
As $M$ is symmetric we have $g\cdot M = (1-[\zeta]) (h\cdot M)$.
Since $M$ has compact support, so does $h\cdot M$.
Therefore $g\cdot M\in (1-[\zeta])\cM^{c}$
 and we have
$$
 \widehat\int g\cdot M
 =
 \zeta^{\textstyle \int h(v) M(v) dv}.
$$
We now calculate the pairing:
$$
 \langle g | M \rangle
 =
 \prod_{\xi\in\mu_{m}}
 \xi^{\textstyle\int f(\xi v)(h(v)-h(\zeta^{-1}v))M(v)dv}.
$$
Replacing $\xi$ by $\zeta^{-1} \xi$ in the second term
 we obtain:
$$
 \langle g | M \rangle
 =
 \prod_{\xi\in\mu_{m}}
 \xi^{\textstyle
 \int f(\xi v) h(v) M(v) dv}
 (\zeta^{-1}\xi)^{\textstyle
 -\int f(\zeta^{-1}\xi v) h(\zeta^{-1} v) M(v) dv}.
$$
Replacing $v$ by $\zeta v$ in the second term
 and using the symmetry of $M$ we obtain:
$$
 \langle g | M \rangle
 =
 \prod_{\xi\in\mu_{m}}
 \xi^{\textstyle
 \int f(\xi v) h(v) M(v) dv}
 (\zeta^{-1}\xi)^{\textstyle
 -\int f(\xi v) h( v) M(v) dv}.
$$
This cancels to give:
$$
 \langle g | M \rangle
 =
 \prod_{\xi\in\mu_{m}}
 \zeta^{\textstyle
 \int f(\xi v) h( v) M(v) dv}.
$$
Since $f$ is fundamental we have:
$$
 \langle g | M \rangle
 =
 \zeta^{\textstyle \int h(v)M(v)dv}.
$$
\sq
\medskip

\subsection{Reduction to a discrete space.}

The arithmetic cocycle $\deca$ depends on a choice of
 an open and closed fundamental domain $F$ for $\mu_{m}$ in
 $\Va\setminus\{0\}$.
In practice it is unnecessary to describe such an $F$ completely.
We shall show that a large part of
 the cocycle depends only on a fundamental
 domain in a discrete quotient of $\Va$.

We shall fix once and for all a $\mu_{m}$-invariant, compact,
 open subgroup $L\subset \Va$.
We shall write $X$ for the (discrete) quotient group $\Va/L$.
The idea is that it is easier to find a fundamental domain
 in $X\setminus\{0\}$ than in $\Va\setminus\{0\}$.

We introduce modules of functions on $X$
 analogous to $\cF^{o}$ and $\cC^{o}$ above.
We let $\cM_{X}$ denote the space of all functions
 $X\setminus\{0\}\to \Z/m$
 and we let $\cM_{X}^{c}$ denote the space of all
 such functions with finite support.
As the action of $\mu_{m}$ on $X\setminus\{0\}$
 has no fixed points, these are acyclic
 $(\Z/m)[\mu_{m}]$ modules.
We define
$$
 \cFo_{X}
 =
 \ker\Big([\mu_{m}]:\cM_{X}\to\cM_{X}\Big),\quad
 \cCo_{X}
 =
 \ker\Big(1-[\zeta]:\cM_{X}^{c}\to\cM_{X}\Big).
$$
There is a pairing $\cFo_{X}\times \cCo_{X}\to\mu_{m}$ given by
\begin{equation}
    \label{simplepairing}
    \langle (1-[\zeta])g | M \rangle_{X}
    =
    \zeta^{\textstyle{\sum_{x\in X\setminus\{0\}}g(x)M(x)}}.
\end{equation}
We have canonical inclusions
 $\iota:\cM_{X}\to \cM$ and $\iota:\cM^{c}_{X}\to \cM^{c}$
 and we have
\begin{equation}
    \label{Xism2}
    \langle g | M \rangle_{X}
    =
    \langle \iota(g) | \iota(M) \rangle.
\end{equation}
We shall write $\Ga^{+}$ for the semi-group of
 elements $\alpha\in\Ga$ such that $\alpha L \supseteq L$.
Let $F_{X}$ be a fundamental domain for the
 action of $\mu_{m}$ on $X\setminus\{0\}$.
We shall suppose our fundamental domain $F$ is chosen
 so that for $v\notin L$ we have
 $v\in F$ if and only if $v+L\in F_{X}$.
As before we shall write $f$ for the characteristic function of $F$.

\begin{lemma}
    \label{Xism4}
    For $\alpha,\beta\in \Ga^{+}$ we have:
    \begin{itemize}
	\item
	The restriction of $f\alpha^{-1}$ to $\Va\setminus \alpha L$
	 is $L$-periodic, and therefore induces
	 a function on $X\setminus \alpha L$.
	\item
	The set $\alpha\beta L-\alpha L$ is $L$-periodic.
	Its characteristic function therefore
	 induces a function on $X$ which is zero on $\alpha L$.
	\item
	We have
	$$
	 \deca^{(f,L)}(\alpha,\beta)
	 =
	 \langle f\alpha^{-1}-f |\alpha\beta L - \alpha L \rangle_{X}.
	$$
    \end{itemize}
\end{lemma}

\proof
The first two statements are easy to check;
 the third follows from (\ref{covariant}), (\ref{def2})
 and (\ref{Xism2}).
\medskip
\sq

\subsection{Examples of arithmetic cocycles}

Let $k$ be a global field containing a primitive $m$-th
 root of unity, and let $\mu_{m}$ be the group of
 $m$-th roots of unity in $k$.
Consider the vector space $V=k^{n}$.

\paragraph{Local examples.}
Let $v$ be a place of $k$ such that $|m|_{v}=1$,
 where $|\cdot|_{v}$ is the absolute value
 on $k_{v}$,
 normalized so that $d(mx)=|m|_{v}dx$
 for any Haar measure $dx$ on $k_{v}$.
We shall write $V_{v}$ for the vector space $V\otimes_{k}k_{v}=k_{v}^{n}$.
The action of $\mu_{m}$ by scalar multiplication on $V_{v}$
 satisfies the conditions of $\Va$ of \S3.1.
We therefore obtain a 2-cocycle on
 $\GL_{n}(k_{v})$ with values in $\mu_{m}$.
We shall refer to this cocycle as $\decv$.

The commutator of $\decv$ on the maximal torus
 has been calculated in Theorems 3 and 5 of \cite{hill2}.
It is given by
\begin{equation}
    \label{localcomm}
    [\alpha,\beta]_{\decv}
    =
    (-1)^{\textstyle{\frac{(|\det\alpha|_{v}-1)(|\det\beta|_{v}-1)}{m^{2}}}}
    \prod_{i=1}^{n}
    (-1)^{\textstyle{\frac{(|\alpha_{i}|_{v}-1)(|\beta_{i}|_{v}-1)}{m^{2}}}}
    (\alpha_{i},\beta_{i})_{m,v},
\end{equation}
 where
$$
 \alpha
 =
 \left(
 \matrix{\alpha_{1} \cr &\ddots\cr&&\alpha_{n}}
 \right),\quad
 \beta
 =
 \left(
 \matrix{\beta_{1} \cr &\ddots\cr&&\beta_{n}}
 \right)
$$
 and $(\cdot,\cdot)_{m,v}$ denotes the $m$-th power Hilbert symbol
 on $k_{v}$.
Our restriction on $v$ amounts to requiring that $(\cdot,\cdot)_{m,v}$
 is the tame symbol.

\paragraph{Semi-global examples.}
Let $S$ be the set of all places $v$ of $k$,
 such that $|m|_{v}\ne 1$.
We shall write $\AS$ for the ring of $S$-ad\`eles of
 $k$, i.e. the restricted topological product
 of the fields $k_{v}$ for $v\notin S$.
Let $V_{\AS}=V\otimes_{k} \AS=\AS^{n}$.
The action of $\mu_{m}$ on $V_{\AS}$ by
 scalar multiplication satisfies the conditions
 of $\Va$ of \S3.1.
We therefore obtain a 2-cocycle on $\GL_{n}(\AS)$
 with values in $\mu_{m}$.
We shall refer to this cocycle as $\decAS$.
The cocycle $\decAS$ is not quite the product of
 the local cocycles $\decv$ for $v\notin S$.
In fact we have (Theorem 3 of \cite{hill2})
 up to a coboundary:
$$
 \decAS(\alpha,\beta)
 =
 \prod_{v\notin S}\decv(\alpha,\beta)
 \prod_{v<w}
 (-1)^{
 \textstyle{\frac{(|\det\alpha|_{v}-1)(|\det\beta|_{w}-1)}{m^{2}}}}.
$$
Here we have chosen an ordering on the set of places $v\notin S$;
 up to a coboundary the right hand side of the above formula is
 independent of the choice of ordering.
As a consequence of this and (\ref{localcomm}),
 we have on the maximal torus in $\GL_{n}(\AS)$:
$$
 [\alpha,\beta]_{\decAS}
 =
 (-1)^{\textstyle{\frac{(|\det\alpha|_{\AS}-1)(|\det\beta|_{\AS}-1)}{m^{2}}}}
 \prod_{v\notin S}
 \prod_{i=1}^{n}
 (-1)^{\textstyle{\frac{(|\alpha_{i}|_{v}-1)(|\beta_{i}|_{v}-1)}{m^{2}}}}
 (\alpha_{i},\beta_{i})_{v,m}.
$$

\paragraph{The positive characteristic case.}
Suppose for a moment that $k$ has positive characteristic.
In this case the set $S$ is empty,
 so in fact we have a cocycle $\decA$
 on the whole group $\GLnA$,
 where $\A$ is the full ad\`ele ring of $k$.
The proof that $\decA$ splits on $\GL_{n}(k)$
 is easy to understand.
 
We may also take $\Va=\A^{n}/k^{n}$.
This is acted on by $\GL_{n}(k)$, so we obtain
 a corresponding cocycle $\dec_{k}$
 on $\GL_{n}(k)$ with values in $\mu_{m}$.
One easily shows that, up to a coboundary,
 $\dec_{k}$ is the restriction of $\decA$.
In this case however, $\Va$ is compact, so we may take $L=\Va$
 in the definition (\ref{def1}) of the cocycle.
With this choice of $L$ one sees immediately that
 the cocycle is trivial.
This shows that the restriction of $\decA$ to $\GL_{n}(k)$ is
 a coboundary.


\subsection{The Gauss--Schering Lemma}

We next indicate the connection between the cocycle $\decAS$
 and the Gauss-Schering Lemma.
Let $k$, $\mu_{m}$ and $S$
 be as in \S3.3.
We shall write $\go^{S}$ for the ring of $S$-integers in $k$.
Recall that for non-zero, coprime $\alpha,\beta\in \go^{S}$, the
 $m$-th power Legendre symbol is defined by
$$
 \left(\frac{\alpha}{\beta}\right)_{S,m}
 =
 \prod_{v\notin S,\; v|\beta}(\alpha,\beta)_{m,v}.
$$
The Gauss-Schering Lemma is a formula for the Legendre symbol,
 commonly used to prove the quadratic reciprocity law in
 undergraduate courses.
Choose a set $R_{\beta}$ of representatives
 of the non-zero $\mu_{m}$-orbits in $\go^{S}/(\beta)$.
Such sets are called $m$-th sets modulo $\beta$.
For any $\xi\in\mu_{m}$ define
$$
 r(\xi)
 =
 |\{x\in R_{\beta}:\alpha x\in\xi R_{\beta}\}|.
$$
The Gauss-Schering Lemma is the statement
$$
 \left(\frac{\alpha}{\beta}\right)_{S,m}
 =
 \prod_{\xi\in\mu_{m}} \xi^{r(\xi)}.
$$

Consider the cocycle $\decAS$ on $\GL_{1}(\AS)$ of \S3.3.
Let $L=\prod_{v\notin S}\go_{v}$;
 we shall also write $L$ for the characteristic function of this set.
The subset $L\subset\AS$ satisfies the conditions of \S3.2.
The semi-group $\Ga^{+}$ of \S3.2 contains $\beta^{-1}$ for
 all non-zero $\beta\in\go^{S}$.
Choose a fundamental domain $F$ for the action of
 $\mu_{m}$ on $(\AS/L) \setminus\{0\}$
 and let $f$ be its characteristic function.
We shall write $f_{\AS}$ for an extension of $f$ to $\AS\setminus\{0\}$.
The Gauss-Schering Lemma may be reformulated as follows.

\begin{proposition}
    For non-zero, coprime $\alpha,\beta\in\go_{S}$ we have:
    $$
     \left(\frac{\alpha}{\beta}\right)_{S,m}
     \;=\;
     \decAS^{f_{\AS},L}(\alpha,\beta^{-1})^{-1}.
    $$
\end{proposition}

\proof
As in \S3.2 we let $X=\AS/L$
 and we let $F$ be a set of representatives for $\mu_{m}$-orbits
 in $X$.
Let $X[\beta]=\{x\in X:\beta x=0\}$.
We shall identify $X[\beta]$ with $\go^{S}/(\beta)$.
Let $F_{\beta}=F\cap X[\beta]$.
The set $F_{\beta}$ is an $m$-th set modulo $\beta$
 in the sense described above.
Therefore the Gauss-Schering Lemma
 states that
 $\left(\alpha/\beta\right)_{S,m} = \prod \xi^{r(\xi)}$,
 where $r(\xi)$ may be rewritten as:
$$
 r(\xi)
 =
 \sum_{x\in F_{\beta}}
 f(\xi^{-1}\alpha x).
$$
This is equivalent to
$$
 r(\xi)
 =
 \sum_{x\in X[\beta]\setminus\{0\}}
 f(x)f(\xi^{-1}\alpha x).
$$
As we have normalized the Haar measure to give $L$
 measure 1, the sum can be rewritten as an integral:
$$
 r(\xi)
 =
 \left\{\int_{\beta^{-1}L}-\int_{L}\right\}
 f(v)f(\xi^{-1}\alpha v)dv.
$$
The proposition follows from the definition
 (\ref{def1}) of $\deca$.
\sq
\medskip


It is common in proofs of reciprocity laws using the
 Gauss--Scher\-ing Lemma to change from a summation
 over $\beta$-division points to a summation over
 $\alpha\beta$-division points when calculating
 $(\alpha/\beta)_{m}$ (see \cite{thesis,habicht}).
The reason why this technique is useful,
 is because whereas the sum over $\beta$-division points is
 $\decAS(\alpha,\beta^{-1})$,
 the sum over $\alpha\beta$-division points is,
 as we can see from Lemma \ref{Xism4},
 related to $\decAS(\alpha,\beta)$.

\subsection{The geometric cocycle}

In this section we let $V_\infty$ be a vector space
 over $\R$
 with an action of the cyclic group $\mu_m$
 such that every non-zero vector in $V_{\infty}$
 has trivial stabilizer in $\mu_m$
 (such representations are called linear space forms).
We shall write $d$ for the dimension of $V_\infty$ as
 a vector space over $\R$.
Let $G_\infty$ be the group $\aut_{\mu_m}(V_\infty)$
 and let $\gg$ be its Lie algebra $\End_{\mu_m}(V_\infty)$.
We shall construct a 2-cocycle $\deci$ on $G_{\infty}$
 with values in $\mu_{m}$.
The construction is analogous to the construction of
 $\deca$ of \S3.1.

If $m=2$ then the group $\mu_{2}$ acts on $V_{\infty}$
 with the non-trivial element acting by multiplication by $-1$.
In this case $G_{\infty}=\GL_{d}(\R)$.
In contrast if $m\ge 3$ then $V_{\infty}$ is a direct sum
 of irreducible two-dimensional representations of $\mu_{m}$.
The group $G_{\infty}$ is then a direct sum of groups
 isomorphic to $\GL_{r(i)}(\C)$,
 with $d=2\sum_{i} r(i)$.

For any singular $r$-cell $\scT:A\to\gg$ in $\gg$ and any
 singular $s$-cell
 $\scU:B\to V_\infty$ we define a singular $r+s$-cell
 $\scT\cdot\scU:A\times B\to V_\infty$ by
$$
 (\scT\cdot\scU)(\xx,\yy)
 =
 \scT(\xx)\cdot \scU(\yy).
$$
This operation extends to a bilinear map
 $C_r(\gg) \times C_s(V_\infty) \to C_{r+s}(V_\infty)$.

We fix an orientation $\ord$ on $V_{\infty}$.
The following lemma will often be used in what follows.

\begin{lemma}
    \label{orient}
    (i) For any $\alpha\in G_{\infty}$, $v\in V_{\infty}$
     and any $\scT\in Z_{d}(V_{\infty},V_{\infty}\setminus\{v\})$
     we have in $\Z/m$:
    $$
     \ord_{\alpha v}(\alpha\scT)
     =
     \ord_{v}(\scT).
    $$
    (ii) For every $\scT\in C_{d}(V_{\infty},V_{\infty}\setminus\{0\})$
     we have in $\Z/m$:
    $$
     \ord_{0}((1-[\zeta])\scT)
     =
     \ord_{0}([\mu_{m}]\scT)
     =
     0.
    $$
\end{lemma}

\proof
(i)
Since $\Z/2\Z$ has no non-trivial automorphisms,
 there is nothing to prove in the case $m=2$.
Assume now $m>2$.
The group $G_{\infty}$ is a direct sum of groups
 isomorphic to $\GL_{r}(\C)$, and is therefore connected.
Hence the action of $G_{\infty}$ on the
 set of orientations of $V_{\infty}$ is trivial.

(ii)
By (i) we have $\ord_{0}((1-[\zeta])\scT)=0$.
On the other hand,
 the following holds in $(\Z/m)[\mu_{m}]$:
$$
 [\mu_{m}]
 =
 \Big([\zeta]+2[\zeta^{2}]+\ldots+(m-1)[\zeta^{m-1}]\Big)
 (1-[\zeta]).
$$
Hence $\ord_{0}([\mu_{m}]\scT)=0$.
\sq
\medskip

Let $d$ be the dimension of $V_{\infty}$ as a vector space over $\R$.
We may choose a finite $d-1$-dimensional cell complex
 $\gS$ in $V_{\infty}\setminus\{0\}$ with the following
 properties:
\begin{itemize}
    \item
    The inclusion $\gS\hookrightarrow V_{\infty}\setminus\{0\}$
     is a homotopy equivalence;
    \item
    The group $\mu_{m}$ permutes the cells in $\gS$
     and each cell has trivial stabilizer in
     $\mu_{m}$.
\end{itemize}
One could for example take $\gS$ to be
 a triangulation of the unit sphere in $V_{\infty}$.
We shall write $\gS_{\bullet}$ for the corresponding
 chain complex with coefficients in $\Z/m$.
It follows from the second condition
 above that each $\gS_{r}$ is a free $(\Z/m)[\mu_{m}]$-module.
A basis consists of a set of representatives of $\mu_{m}$-orbits
 of $r$-cells.
Thus each $\gS_{r}$ satisfies the exact sequence (\ref{star}).

Choose a cycle $\omega\in\gS_{d-1}$ whose
 homology class is mapped to 1 by the isomorphisms:
$$
 H_{d-1}(\gS)
 \to
 H_{d-1}(V_{\infty}\setminus\{0\})
 \stackrel{\partial^{-1}}{\to}
 H_{d}(V_{\infty},V_{\infty}\setminus\{0\})
 \stackrel{\ord_{0}}{\to}
 \Z/m.
$$
By Lemma \ref{orient} we know that $(1-[\zeta])[\omega]=0$
 in $H_{d-1}(\gS)$.
As $\gS_{d}=0$, this implies that $(1-[\zeta])\omega=0$
 in $\gS_{d-1}$.
Thus by the exact sequence (\ref{star}) there is a $\scD\in\gS_{d-1}$
 such that $\omega=[\mu_{m}]\scD$.
We may think of $\scD$ as a fundamental domain
 for $\mu_{m}$ in $\gS$.
As $\omega$ is a cycle we have:
$$
 [\mu_{m}]\partial\scD
 =
 \partial[\mu_{m}]\scD
 =
 \partial\omega
 =
 0.
$$
Thus, by the exact sequence (\ref{star}),
 there is an $\scE\in\gS_{d-2}$ such that
$$
 \partial \scD = (1-[\zeta])\scE.
$$

\begin{definition}
    We may define the geometric cocycle $\deci^{(\scD)}$
     for generic $\alpha,\beta\in G_{\infty}$ as follows:
    $$
     \deci^{(\scD)}(\alpha,\beta)
     =
     \zeta^{\textstyle{-\ord_{0}([1,\alpha,\alpha\beta]\cdot\scE)}}.
    $$
\end{definition}

\begin{remark}
    In fact we have
     $\deci^{(\scD)}(\alpha,\beta)
      =\widehat{\ord}_{0}(-[1,\alpha,\alpha\beta]\cdot\partial\scD)$,
     where $\widehat{\ord}_{0}$ is a lifted map
     in the sense of Lemma \ref{easy}.
    The $(\Z/m)[\mu_{m}]$-module 
     $Z_{d}(V_{\infty},V_{\infty}\setminus\{0\})$
     is free.
    However we shall not use this interpretation.
\end{remark}

If the choice of $\scD$ is clear or irrelevant,
 then we shall omit it from the notation.
In order to define $\deci (\alpha,\beta)$
 on all pairs $(\alpha,\beta)$ rather than just generically,
 we must define
 $\ord_{0}([1,\alpha,\alpha\beta]\cdot\scE)$ when $0$ is
 on the boundary.
To do this we slightly modify our definition.
We define $[1,\alpha,\alpha\beta]$ to be the map $\gg^3\to C_2(\gg)$
 defined by
$$
 (\epsilon_1,\epsilon_2,\epsilon_3)
 \mapsto
 [1+\epsilon_1,\alpha(1+\epsilon_2),\alpha\beta(1+\epsilon_3)].
$$
Then $[1,\alpha,\alpha\beta]\cdot\scE$ becomes a map $\gg^3\to C_d(V_\infty)$.
We fix a basis $\{b_1,\ldots,b_r\}$ for $\gg$ over $\R$.
For $\epsilon=\sum x_i b_i$ we shall use the notation:
$$
 ``{\blim_{\epsilon\to 0^+}}"
 =
 \lim_{x_1\to 0^+}
 \cdots
 \lim_{x_r\to 0^+}.
$$
With this notation we define
$$
 \ord_0([1,\alpha,\alpha\beta]\cdot\scE)
 =
 \blim_{\epsilon_1\to 0^+}
 \blim_{\epsilon_2\to 0^+}
 \blim_{\epsilon_3\to 0^+}
 \ord_0\Big(
 [1+\epsilon_1,\alpha(1+\epsilon_2),\alpha\beta(1+\epsilon_3)]\cdot\scE
 \Big).
$$
We will prove in \S5 that, for a suitable choice of $\gS$,
 the above limit exists for all $\alpha,\beta\in G_\infty$.
It is worth mentioning that the above limits do not
 necessarily commute.

\begin{theorem}
    $\deci$ is a 2-cocycle.
\end{theorem}

\proof
Let $\alpha,\beta,\gamma\in G_\infty$
 and consider the $3$-cell in $\gg$:
$$
 \scA
 =
 [1+\epsilon_1,\alpha(1+\epsilon_2),
 \alpha\beta(1+\epsilon_3),\alpha\beta\gamma(1+\epsilon_4)].
$$
We have straight from the definition:
$$
 \partial\deci(\alpha,\beta,\gamma)
 =
 \blim_{\epsilon_1\to 0^+}
 \blim_{\epsilon_2\to 0^+}
 \blim_{\epsilon_3\to 0^+}
 \blim_{\epsilon_4\to 0^+}
 \zeta^{\textstyle{-\ord_{0}((\partial\scA)\cdot\scE)}}.
$$
The boundary of $\scA\cdot\scE$ is given by
$$
 \partial(\scA\cdot\scE)
 =
 \partial(\scA)\cdot\scE
 +
 \scA\cdot\partial\scE.
$$
(Since coefficients are in $\Z/m$,
 the sign is only important for $m\ge3$, and in these cases
 $d$ is even).
Taking the order of this at 0 and then taking limits
 we obtain:
$$
 \partial\deci(\alpha,\beta,\gamma)
 =
 \blim_{\epsilon_1\to 0^+}
 \blim_{\epsilon_2\to 0^+}
 \blim_{\epsilon_3\to 0^+}
 \blim_{\epsilon_4\to 0^+}
 \zeta^{\textstyle{\ord_{0}(\scA\cdot\partial\scE)}}.
$$
We must show that the right hand side here is 1.
We know that $(1-[\zeta])\scE=\partial \scD$.
Therefore $(1-[\zeta])\partial\scE=0$,
 so by the exact sequence (\ref{star}),
 there is a $\scB\in\gS_{d-3}$ such that
$$
 \partial\scE = [\mu_{m}]\scB,
$$
This implies
$$
 \scA\cdot\partial\scE
 =
 [\mu_{m}](\scA\cdot\scB).
$$
The result now follows from Lemma \ref{orient}.
\sq
\medskip

We next verify that our definition is independent of the various
 choices made.

\begin{proposition}
    (i) The cocycle $\deci^{(\scD)}$ is independent of the choice of
     $\scE$, the orientation $\ord$ and the generator $\zeta$ of
     $\mu_{m}$.\\
    (ii) The cohomology class of $\deci^{(\scD)}$ is independent of the
     choice of $\scD$.
\end{proposition}

\proof
(i)
We first fix $\scD$ and choose another $\scE^{\prime}$ such that
$$
 (1-[\zeta])\scE
 =
 (1-[\zeta])\scE^{\prime}
 =
 \partial\scD.
$$
By the exact sequence (\ref{star}) there is an $\scA\in\gS_{d-2}$ such that
$$
 \scE^{\prime}-\scE
 =
 [\mu_{m}]\scA.
$$
This implies by Lemma \ref{orient}:
\begin{eqnarray*}
 \deci^{\scE^{\prime}}(\alpha,\beta)
 &=&
 \deci^{\scE}(\alpha,\beta)
 \zeta^{\textstyle{-\ord_{0}([\mu_{m}][1,\alpha,\alpha\beta]\cdot\scA)}}\\
 &=&
 \deci^{\scE}(\alpha,\beta).
\end{eqnarray*}

Now suppose we choose a different orientation $\ord'$.
We have $\ord'=u\cdot\ord$ for some $u\in(\Z/m)^{\times}$.
We may therefore choose $\omega'=u^{-1}\omega$,
 $\scD'=u^{-1}\scD$ and $\scE'=u^{-1}\scE$.
With these choices we have
$$
 \ord_{0}'([1,\alpha,\alpha\beta]\cdot\scE')
 =
 \ord_{0}([1,\alpha,\alpha\beta]\cdot\scE).
$$

Finally suppose $\zeta^{\prime u}=\zeta$ for some
 $u\in(\Z/m)^{\times}$.
We then have in $(\Z/m)[\mu_{m}]$:
$$
 (1-\zeta)
 =
 (1-\zeta')(1+\zeta'+\ldots+\zeta^{\prime u-1}).
$$
We may therefore take
$$
 \scE'
 =
 (1+\zeta'+\ldots+\zeta^{\prime u-1})\scE.
$$
This implies by Lemma \ref{orient}
$$
 \ord_{0}([1,\alpha,\alpha\beta]\cdot\scE')
 =
 u \cdot \ord_{0}([1,\alpha,\alpha\beta]\cdot\scE).
$$
Therefore
$$
 \zeta^{\prime\textstyle{-\ord_{0}([1,\alpha,\alpha\beta]\cdot\scE')}}
 =
 \zeta^{\textstyle{-\ord_{0}([1,\alpha,\alpha\beta]\cdot\scE)}}.
$$

(ii)
We now allow $\scD$ to vary.
We choose $\scD^{\prime}$ to satisfy
$$
 [\mu_{m}]\scD^{\prime}
 =
 [\mu_{m}]\scD.
$$
By the exact sequence (\ref{star}) there is a $\scB\in\gS_{d-1}$ such that
$$
 \scD^{\prime}
 =
 \scD + (1-[\zeta])\scB.
$$
Thus
$$
 \partial\scD^{\prime}
 =
 \partial\scD + (1-[\zeta])\partial \scB.
$$
We may therefore choose $\scE^{\prime}=\scE+\partial\scB$.
Note that we have
\begin{eqnarray*}
 \partial([1,\alpha,\alpha\beta]\cdot\scB)
 &=&
 \partial([1,\alpha,\alpha\beta])\cdot\scB
 +
 [1,\alpha,\alpha\beta]\cdot\partial\scB\\
 &=&
 [1,\alpha]\cdot\scB-[1,\alpha\beta]\cdot\scB+[\alpha,\alpha\beta]\cdot\scB
 +
 [1,\alpha,\alpha\beta]\cdot\partial\scB.
\end{eqnarray*}
This implies using Lemma \ref{orient}:
\begin{eqnarray*}
 \ord_{0}([1,\alpha,\alpha\beta]\cdot\partial\scB)
 &=&
 \ord_{0}(
 [1,\alpha\beta]\cdot\scB
 -[1,\alpha]\cdot\scB
 -[\alpha,\alpha\beta]\cdot\scB)\\
 &=&
 \ord_{0}([1,\alpha\beta]\cdot\scB)
 -
 \ord_{0}([1,\alpha]\cdot\scB)
 -
 \ord_{0}([1,\beta]\cdot\scB).
\end{eqnarray*}
Putting things back together we obtain:
$$
 \deci^{(\scD')}(\alpha,\beta)
 =
 \deci^{(\scD)}(\alpha,\beta)
 \frac{\tau(\alpha)\tau(\beta)}{\tau(\alpha\beta)},
$$
 where
$$
 \tau(\alpha)
 =
 \zeta^{\textstyle{\ord_{0}([1,\alpha]\cdot\scB)}}.
$$
\sq
\medskip

\subsection{Example: $\GL_{2}(\R)$}
If $V_{\infty}$ is $1$-dimensional then the fundamental domain
 $\scD$ is zero-dimensional, and so we have $\scE=0$
 and $\deci$ is always 1.
The smallest non-trivial example is the case $m=2$
 and $V_{\infty}=\R^{2}$ with the group $\mu_{2}$
 acting on $\R^{2}$ by scalar multiplication.
We shall consider this example now.
We have $G_{\infty}=\GL_{2}(\R)$.

As $m=2$ there is no need to worry about a choice of orientation.
We may take our fundamental domain $\scD$ to be the half-circle:
$$
 \scD
 =
 \left\{\left(\matrix{x \cr y}\right):
 x^{2}+y^{2}=1,\; y\ge 0
 \right\}.
$$
The boundary of this consists of the two points $\left(\matrix{1 \cr 0}\right)$
 and $\left(\matrix{-1 \cr 0}\right)$.
We may therefore take
 $\scE=[v]$, where $v=\left(\matrix{1 \cr 0}\right)$.
The cocycle is then given (generically) by:
$$
 \deci^{(\scD)}(\alpha,\beta)
 =
 \left\{
 \begin{array}{rl}
      1 & \hbox{if } 0\notin [v,\alpha v,\alpha\beta v],\\
     -1 & \hbox{if } 0\in [v,\alpha v,\alpha\beta v].
 \end{array}
 \right.
$$
By choosing a different $\scD$ we may replace $v$ by any
 other non-zero vector to obtain a cohomologous cocycle.

For later use we calculate the commutator of $\deci$
 on the standard torus in $\GL_{2}(\R)$.

\begin{proposition}
    \label{twocom}
    The commutator of $\deci$ on the
     subgroup of diagonal matrices in $\GL_{2}(\R)$
     is given by:
    $$
     \left[
     \left(\matrix{\alpha_{1}\cr & \alpha_{2}}\right),
     \left(\matrix{\beta_{1}\cr & \beta_{2}}\right)
     \right]_{\deci}
     =
     (\alpha_{1},\beta_{1})_{\R,2}
     (\alpha_{2},\beta_{2})_{\R,2}
     (\det\alpha,\det\beta)_{\R,2}.
    $$
    The right hand side here consists of quadratic Hilbert
     symbols on $\R$.
\end{proposition}

\proof
Commutators are is continuous, bimultiplicative
 and alternating.
The right hand side of the above formula is also continuous,
 bimultiplicative and alternating in $\alpha,\beta$.
It is therefore sufficient to check the formula
 in the case $\alpha=\left( \matrix{1 \cr &-1}\right)$,
 $\beta=\left(\matrix{-1 \cr &1}\right)$.
To calculate the commutator there we choose
 $v=\left(\matrix{1\cr 1}\right)$.
With this choice we have
\begin{eqnarray*}
    \left[(1+\epsilon_{1})v,\alpha(1+\epsilon_{2})
    v,\alpha\beta(1+\epsilon_{3}) v\right]
    &=&
    \left[
    (1+\epsilon_{1}) \left(\matrix{1\cr 1}\right),
    (1+\epsilon_{2})^{\alpha} \left(\matrix{1\cr -1}\right),
    (1+\epsilon_{3}) \left(\matrix{-1\cr -1}\right)
    \right],\\
    \left[(1+\epsilon_{1})v,\beta(1+\epsilon_{2})
    v,\beta\alpha(1+\epsilon_{3}) v\right]
    &=&
    \left[
    (1+\epsilon_{1})\left(\matrix{1\cr 1}\right),
    (1+\epsilon_{2})^{\beta}\left(\matrix{-1\cr 1}\right),
    (1+\epsilon_{3})\left(\matrix{-1\cr -1}\right)
    \right].
\end{eqnarray*}
We therefore have in $\Z/2$:
$$
 \blim_{\epsilon_{1}\to 0^{+},\epsilon_{2}\to 0^{+},\epsilon_{3}\to 0^{+}}
 \left(
 \begin{array}{ll}
     \ord_{0}(
     [(1+\epsilon_{1})v,
      \alpha(1+\epsilon_{2}) v,
      \alpha\beta(1+\epsilon_{3}) v])\\
     -\ord_{0}(
     [(1+\epsilon_{1})v,
      \alpha(1+\epsilon_{2})v,
      \alpha\beta(1+\epsilon_{3}) v])
 \end{array}
 \right)
 =
 1.
$$
This implies $\frac{\deci(\alpha,\beta)}{\deci(\beta,\alpha)}=-1$,
 which verifies the result.
\sq

\subsection{Stability of the cocycles}

Suppose $V_\infty$ is the direct sum
 of the representations $V_1$ and $V_2$ of $\mu_m$.
The group $G_\infty$ contains the direct sum
 of $G_1=\Aut_{\mu_m}(V_1)$ and $G_2=\Aut_{\mu_m}(V_2)$.
We have defined cocycles $\deci$ on $G_\infty$,
 $\deci^{(1)}$ on $G_1$ and $\deci^{(2)}$ on $G_2$.
The next result describes how these are related.

\begin{proposition}
    \label{stability}
    Let $\alpha=(\alpha_1,\alpha_2)$ and $\beta=(\beta_{1},\beta_{2})$
     denote elements of $G_1\oplus G_2$.
    We have up to a coboundary:
    $$
     \deci(\alpha,\beta)
     =
     \deci^{(1)}(\alpha_1,\beta_1)
     \deci^{(2)}(\beta_1,\beta_2)
     (\det(\alpha_1),\det(\beta_2))_{\R,2}.
    $$
    Here $(\cdot,\cdot)_{\R,2}$ denotes the quadratic
     Hilbert symbol on $\R$
     and $\det$ is the determinant over the base field $\R$.
\end{proposition}

Note that for $m\ge 3$, $\alpha_{i}$ and $\beta_{i}$ have
 positive determinant, so in this case the final term above
 vanishes.
\medskip

\proof
We shall first consider the case $m=2$.
Thus we have $G_1=\GL_a(\R)$, $G_2=\GL_b(\R)$
 and $G_\infty=\GL_{a+b}(\R)$.
There is an isomorphism:
$$
 H^2(G_1\oplus G_2,\mu_2)
 \cong
 H^2(G_1,\mu_2)
 \oplus
 H^1(G_1,H^1(G_2,\mu_2))
 \oplus
 H^2(G_2,\mu_2)
$$
The middle component of the isomorphism
 is given by the commutator
 $[\alpha_{1},\beta_{2}]$ ($\alpha_{1}\in G_1$, $\beta_{2}\in G_2$);
 the other two components are given by restriction.
We must show that the image of $\deci$
 is that described in the proposition.

We first examine the restriction of $\deci$ to $G_1$.
We may assume without loss of generality
 that $V_2=\R$.
We shall assume for the moment
 that $V_1=\R^n$ with $n\ge 2$;
The case $n=1$ will be dealt with separately.

For any $r$-cell $\scA:A\to\R^n$ we define
 we define two $r+1$-cells
 $\scA^+,\scA^-:A\times I\to \R^{n+1}$ by:
$$
 (\scA^+)(\xx,t)
 =
 (1-t)\scA(\xx)+ te_{n+1},
$$
$$
 (\scA^-)(\xx,t)
 =
 (1-t)\scA(\xx)- te_{n+1},
$$
 and we shall write $\scA^\pm=\scA^+-\scA^-$.
Here $e_{n+1}$ is the $n+1$-st standard basis element in $\R^{n+1}$.
The above construction has the following properties
 which are easily checked:
\begin{itemize}

\item[1.]
For $r\ge 1$, we have modulo degenerate cells:
 $\partial (\scA^\pm) = (\partial\scA)^\pm$.

\item[2.]
For $\alpha\in G_1=\GL_n(\R)$
 we have $\alpha(\scA^\pm)=(\alpha\scA)^\pm$.
Furthermore $([-1]\cdot \scA)^{\pm}=[-1]\cdot(\scA^\pm)$
 (here $[-1]$ is acting by scalar multiplication on $V_{\infty}$,
 rather than on the coefficients of chains).

\item[3.]
If $\scA$ is an $n$-cell in $\R^n$
 then $\ord_{0,\R^n}(\scA) = \ord_{0,\R^{n+1}}(\scA^\pm)$.
(The choice of orientation here is unnecessary as $m=2$).
One should understand this formula as meaning that if
 one side is defined then so is the other and they are equal.

\end{itemize}
Let $\omega_{1}$ be the generator of $H_{n-1}(V_{1}\setminus\{0\})$
 as in \S3.5.
By the first and third properties above,
 we may take $\omega=\omega_{1}^{\pm}$
 as our generator for $H_{n}(V_{\infty}\setminus\{0\})$.
By the second property,
 we may choose $\scD=\scD_1^\pm$.
Since $n\ge 2$, the first property
 implies that we may take $\scE=\scE_1^\pm$.

Let $\alpha_{1},\beta_{1}\in G_1$.
By the second property we have:
$$
 [1,\alpha_{1},\alpha_{1}\beta_{1}]\cdot\scE
 =
 ([1,\alpha_{1},\alpha_{1}\beta_{1}]\cdot\scE_{1})^\pm.
$$
Hence by the third property, it follows that:
$$
 \deci^{(\scD_{1})}(\alpha_{1},\beta_{1})
 =
 \deci^{(\scD)}(\alpha_{1},\beta_{1}),
$$
 so the restriction of $\deci$ to $G_1$ is $\deci^{(1)}$.

We may check by hand that this still holds in the case $n=1$,
 where $\deci^{(1)}$ is trivial
 (simply take $\scE=\left[\left(\matrix{0\cr 1}\right)\right]$
 and draw a picture).
By the same reasoning we also know that the restriction
 of $\deci$ to $G_2$ is cohomologous to $\deci^{(2)}$.

It remains only to prove the formula for
 the commutator $[\alpha_{1},\beta_{2}]_{\deci}$,
 where $\alpha_{1}\in G_1$, $\beta_{2}\in G_2$.
As commutators are bimultiplicative,
 this only depends on $\det(\alpha_{1})$ and $\det(\beta_{2})$
 (as $\SL_{n}(\R)$ is the commutator subgroup of $\GL_{n}(\R)$).
Furthermore, by what we have already proved,
 we may assume without loss of generality that
 $V_1=V_2=\R$.
We are reduced to calculating the commutator:
$$
 \left[
 \left(\matrix{-1 & 0\cr 0 & 1}\right),
 \left(\matrix{1 & 0\cr 0 & -1}\right)
 \right]_{\deci}.
$$
The result now follows from Proposition \ref{twocom}.

Finally suppose $m\ge 3$.
In this case $G_{1}$ and $G_{2}$ are connected,
 so the middle commutator term is trivial.
We need only verify that the
 restriction of $\deci$ to $G_{1}$ is $\deci^{(1)}$.
By induction it is sufficient to prove this in the case
 that $V_{2}$ is a simple $\mu_{m}$-module.
Since $V_{2}$ is a linear space form of $\mu_{m}$,
 this implies that $V_{2}$ is 2-dimensional.
The result follows as in the case $m=2$ but replacing
 the construction $\scA^{\pm}$ by a construction which
 increases the dimensions of cells by 2.
This is left to the reader.
\sq
\medskip

\begin{corollary}
Let $m=2$ and let $\deci$ be the cocycle
 on $\GL_n(\R)$ constructed from the action of
 $\mu_2$ on $\R^n$.
Then, on the standard torus in $\GL_n(\R)$,
 the commutator of $\deci$ is given by
$$
 \left[ \alpha, \beta \right]_{\deci}
 =
 (\det\alpha,\det\beta)_{\R,2}
 \prod_{i=1}^{n}
 (\alpha_i,\beta_i)_{\R,2},
$$
 where
$$
 \alpha
 =
 \left(
 \matrix{\alpha_1 &&\cr &\ddots&\cr &&\alpha_n}
 \right),
 \quad
 \beta
 =
 \left(
 \matrix{\beta_1 &&\cr &\ddots&\cr &&\beta_n}
 \right).
$$
\end{corollary}

\proof
This follows by induction from Proposition \ref{stability}.
\sq
\medskip

\subsection{Example: $\GL_{1}(\C)$}

We now calculate a specific example which we shall need
 in the next section.
We choose an embedding of $\iota:\mu_m\hookrightarrow\C^\times$ and we
 let $V_\infty=\C$ with the action of $\mu_m$ given by $\iota$.
Thus $\dec_\infty$ defines an element of $H^2(\C^\times,\mu_m)$
 and we shall now calculate this element.
As $\C^\times$ is abelian, 2-cocycles on this group may be studied by
 studying their commutators.
However since $\C^\times$ is connected, the commutator of
 every $2$-cocycle is trivial.
We therefore have
$$
 H^2(\C^\times,\mu_m)
 =
 H^2_{sym}(\C^\times,\mu_m).
$$
By Klose's isomorphism (see \S2.3), it follows that
\begin{eqnarray*}
 H^2(\C^\times,\mu_m)
 &\cong&
 \Z/m,\\
 \sigma
 &\mapsto&
 a,\quad\hbox{where }
 \zeta^a=\prod_{i=1}^{m-1}\sigma(\iota(\zeta)^i,\iota(\zeta)).
\end{eqnarray*}

\begin{proposition}
The image under the above isomorphism of $\deci$ is $1$.
\end{proposition}

\proof
We choose the generator $\zeta$ so that
 $\iota(\zeta)=e^{2\pi i/m}$.
The group $H_{2}(\C,\C^{\times})$ is generated by
 the following element:
$$
 \scA
 =
 [1,\iota(\zeta),\iota(\zeta)^2]
 +
 [1,\iota(\zeta)^2,\iota(\zeta)^3]
 +
 \ldots
 +
 [1,\iota(\zeta)^{m-1},\iota(\zeta)^m].
$$
We shall fix our orientation on $\C$ such that $\ord_{0}(\scA)=1$.
We may therefore take
$$
 \omega
 =
 \partial\scA
 =
 [1,\iota(\zeta)]
 +
 [\iota(\zeta),\iota(\zeta)^2]
 +
 \ldots
 +
 [\iota(\zeta)^{m-1},1].
$$
We may then choose $\scD$ to be the line segment $[1,e^{2\pi i/m}]$.
Thus
$$
 \partial\scD
 =
 [e^{2\pi i/m}]-[1]
 =
 (1-[\iota(\zeta)])(-[1]).
$$
We may therefore choose $\scE=-[1]$.
With this choice of $\scE$ we have as required:
$$
 \prod_{i=1}^{m-1}\deci(\iota(\zeta)^i,\iota(\zeta))
 =
 \zeta^{\textstyle{\ord_0(\scA)}}
 =
 \zeta.
$$
\sq
\medskip

The corresponding central extension (normalized
 as described in \S2.2) is as follows:
$$
 \begin{array}{ccccccccc}
 1 &\to& \mu_{m}& \stackrel{\iota}{\to}&
 \C^{\times} &\to &\C^{\times}& \to& 1,\\
 &&&& \alpha &\mapsto&\alpha^{m}.
 \end{array}
$$

\subsection{The totally complex case.}

%

In this section we suppose $k$
 is a totally complex number field
 containing a primitive $m$-th root of unity and let
 $\mu_{m}$ be the group of all $m$-th roots of unity in $k$.
We let $V_{\infty}=k_{\infty}^{n}$,
 where $k_{\infty}=k\otimes_{\Q}\R$.
The action of $\mu_{m}$ by scalar multiplication on $V_{\infty}$
 satisfies the conditions of \S3.5.
The group $G_{\infty}$ contains $\GL_{n}(k_{\infty})$.
We therefore obtain by restriction
 a cocycle $\deci$ on $\GL_{n}(k_{\infty})$.
In this section we shall study this cocycle.

\subsubsection{The cocycle up to a coboundary.}

The determinant map gives as isomorphism:
$$
 \GL_{n}(k_{\infty})/\SL_{n}(k_{\infty})
 \stackrel{\det}{\cong}
 k_{\infty}^{\times}.
$$
This gives us an inflation map:
$$
 H^{2}(k_{\infty}^{\times},\mu_{m})
 \to
 H^{2}(\GL_{n}(k_{\infty}),\mu_{m}).
$$
As $\SL_{n}(k_{\infty})$ is both connected and simply connected,
 the groups $H^{1}(\SL_{n}(k_{\infty}),\mu_{m})$
 and $H^{2}(\SL_{n}(k_{\infty}),\mu_{m})$
 are both trivial.
Hence by from the Hochschild--Serre spectral sequence,
 it follows that the above map is an isomorphism.

As $k_{\infty}^{\times}$
 is abelian we may speak of the
 commutators of cocycles; however as
 $k_{\infty}^{\times}$
 is connected, these commutators are all trivial.
Thus we have
$$
 H^{2}(\GL_{n}(k_{\infty}),\mu_{m})
 \cong
 H^{2}_{sym}(k_{\infty}^{\times},\mu_{m}),
$$
Let $S_{\infty}$ be the set of archimedean places of $k$.
We have a decomposition:
$$
 k_{\infty}^{\times}
 =
 \bigoplus_{v\in S_{\infty}} k_{v}^{\times}.
$$
By the results described in \S2.3, we have
$$
 H^{2}(\GL_{n}(k_{\infty}),\mu_{m})
 \cong
 \bigoplus_{v\in S_{\infty}}
 H^{2}_{sym}(k_{v}^{\times},\mu_{m}).
$$
Klose's isomorphism (\S2.3) now gives:
$$
 H^{2}(\GL_{n}(k_{\infty}),\mu_{m})
 \cong
 \bigoplus_{v\in S_{\infty}}
 \Z/m.
$$
By the results of the previous two sections 3.7 and 3.8,
 we know that the image of $\deci$ under the above
 isomorphism is $(1,\ldots,1)$.

\subsubsection{The group extension.}
We now describe the central extension of $\GL_{n}(k_{\infty})$,
 corresponding to $\deci$.
Fix $v\in S_{\infty}$, so $k_{v}$ is non-canonically isomorphic
 to $\C$.
Define a subgroup $\widetilde{\GL}_{n}(k_{v})$ of $\GL_{n+1}(k_{v})$
 as follows:
$$
 \widetilde{\GL}_{n}(k_{v})
 =
 \left\{\left(\matrix{\alpha&0\cr 0&\beta}\right) :
 \begin{array}{c}
     \alpha\in\GL_{n}(k_{v}),\beta\in k_{v}^{\times},\\
     \det\alpha=\beta^{m}
 \end{array}
 \right\}.
$$
The group extension of $\GL_{n}(k_{v})$
 defined by the restriction of $\deci$,
 is concretely realized as follows:
$$
 \begin{array}{ccccccccc}
     1&\to&\mu_{m}&\to&\widetilde{\GL}_{n}(k_{v})
     &\to&\GL_{n}(k_{v})&\to&1,\\
     &&\zeta&\mapsto&
     \left(\matrix{I_{n}&0\cr 0&\iota_{v}(\zeta)}\right),\\
     &&&&\left(\matrix{\alpha&0\cr 0&\beta}\right)&\mapsto&\alpha.
 \end{array}
$$
Here $\iota_{v}:k\hookrightarrow k_{v}$ denotes the embedding
 corresponding to the place $v$.
Let $\mu_{m}(k_{\infty})$ be the $m$-torsion subgroup of 
 $k_{\infty}^{\times}$.
Define a subgroup $K$ of $\mu_{m}(k_{\infty})$ to be the kernel of
 the homomorphism $h:\mu_{m}(k_{\infty})\to\mu_{m}$ defined by
$$
 h(\xi_{v})
 =
 \prod_{v\in S_{\infty}}\iota_{v}^{-1}(\xi_{v}).
$$
The full group extension $\widetilde{\GL}_{n}(k_{\infty})$ is the quotient
$$
 \left(
 \bigoplus_{v\in S_{\infty}}
 \widetilde{\GL}_{n}(k_{v})
 \right)/K.
$$

\subsubsection{How the cocycle splits.}

We will need to calculate precisely how the cocycle
 $\deci$ splits.
This is essential in order to find a formula for
 the Kubota symbol.
Consider the set
$$
 U
 =
 \{\alpha\in \GL_{n}(k_{\infty}):
 \hbox{$\alpha$ has no negative real eigenvalue}\}.
$$
The set $U$ is contractible, as it is a star body from $1$.
It is also a dense open subset of $\GL_{n}(k_{\infty})$,
 as may be seen from the Jordan canonical form.

\begin{lemma}
    If $\alpha\in U$ then the function $\deci$
     is locally constant at the point $(1,\alpha)$.
\end{lemma}

\proof
By definition we have:
$$
 \deci(1,\alpha)
 =
 \zeta^{\textstyle{\ord_{0}([1,1,\alpha]\cdot\scE)}}.
$$
If we suppose that $\deci$ is discontinuous at
 $(1,\alpha)$ then this implies that $0$ is on the
 base set of $[1,1,\alpha]\cdot\scE$.
Thus there is a $v$ in the base set of $\scE$ such that
$$
 (t+(1-t)\alpha) \cdot v
 =
 0,
 \quad t\in(0,1).
$$
This implies that $v$ is a (non-zero) eigenvector of $\alpha$
 with eigenvalue $\frac{-t}{1-t}<0$.

~
\sq
\medskip

Let $\Mu=k_{\infty}^{\times}/K$.
Using the map $h$,
 we may identify $\mu_{m}$ with the subgroup $\mu_{m}(k_{\infty})/K$
 of $\Mu$,
 and hence we may regard $\deci$ as taking values in $\Mu$;
 as such it is a coboundary.
We define a function $w:\GL_{n}(k_{\infty})\to \Mu$ which splits
 $\deci$.
First suppose $\alpha\in U$ and consider the path in $\gg$:
$$
 \path_{\alpha}=[1,\alpha].
$$
As $\alpha$ has no negative real eigenvalues,
 this path is contained in $\GL_{n}(k_{\infty})$.
There is a unique continuous path $q_{\alpha}:I\to k_{\infty}^{\times}$
 defined by:
$$
 q_{\alpha}(0)=1,\quad
 q_{\alpha}(t)^{m}=\det\path_{\alpha}(t).
$$
We define $w(\alpha)=q_{\alpha}(1)$.
More generally for $\alpha\in \GL_{n}(k_{\infty})$ we define
$$
 w(\alpha)
 =
 \blim_{\epsilon\to 0^{+}}
 w(\alpha(1+\epsilon)).
$$
Clearly $w(\alpha)^{m}=\det\alpha$.
Hence, for $\alpha\in \SL_{n}(k_{\infty})$
 we have $w(\alpha)\in \mu_{m}(k_{\infty})$.

\begin{theorem}
    \label{complexsplit}
    For $\alpha,\beta\in \GL_{n}(k_{\infty})$ we have
     $\frac{w(\alpha)w(\beta)}{w(\alpha\beta)}\in \mu_{m}(k_{\infty})$.
    Furthermore in $\mu_{m}$ we have:
    $$
     \deci(\alpha,\beta)
     =
     h\left(\frac{w(\alpha)w(\beta)}{w(\alpha\beta)}\right).
    $$
    For $\alpha,\beta\in \SL_{n}(k_{\infty})$ we have
    $$
     \deci(\alpha,\beta)
     =
     \frac{hw(\alpha)hw(\beta)}{hw(\alpha\beta)}.
    $$
\end{theorem}

The theorem gives an explicit splitting of the
 image of $\deci$ in $Z^{2}(\GL_{n}(k_{\infty}),\Mu)$
 and also a splitting of the restriction of
 $\deci$ to $\SL_{n}(k_{\infty})$.
\medskip

\proof
The first statement is trivial and the third statement follows
 immediately from the second.
We shall prove the second statement.
By the limiting definition of both sides
 of the formula, it is sufficient to
 prove this in the case $\alpha,\beta,\alpha\beta\in U$.
As $U$ is simply connected
 there is a unique continuous section $\tau:U\to \widetilde{\GL}_{n}(k_{\infty})$
 such that $\tau(1)=1$.
This section is given by
$$
 \alpha
 \mapsto
 \left(\matrix{\alpha&0\cr 0&w(\alpha)}\right).
$$
Now consider the realization of $\widetilde{\GL}_{n}(k_{\infty})$
 as $\GL_{n}(k_{\infty})\times\mu_{m}$ with multiplication
 given by
 $(\alpha,\xi)(\beta,\chi)=(\alpha\beta,\xi\chi\deci(\alpha,\beta))$.
To prove the theorem it is sufficient to show that the map
 $U\to \widetilde{\GL}_{n}(k_{\infty})$ given by
 $\alpha\mapsto (\alpha,1)$ is continuous, and hence coincides
 with $\tau$.

There is a neighbourhood $U_{0}$ of 1 in $\GL_{n}(k_{\infty})$ such that for
 $\alpha,\beta\in U_{0}$ we have $\dec(\alpha,\beta)=1$.
Therefore the map $\beta\to (\beta,1)$ is a local
 homomorphism on $U_{0}$, and is therefore continuous.
On the other hand by the previous lemma,
 for $\alpha\in U$ there is a neighbourhood
 $U_{\alpha}$ of 1 in $U_{0}$,
 such that for $\beta\in U_{\alpha}$ we have:
$$
 (\alpha\beta,1)=(\alpha,1)(\beta,1).
$$
Thus the left hand side is a continuous function
 of $\beta\in U_{\alpha}$, so the theorem is proved.
\sq
\medskip

\begin{remark}
    The above theorem shows that
     in the complex case the cocycle $\deci^{(\scD)}$
     does not depend on the fundamental domain
     $\scD$.
    The cocycle does depend on the basis of
     $\gg$ used to define the limits, but this
     dependence is only for $\alpha$, $\beta$
     or $\alpha\beta$ in $\GL_{n}(k_{\infty})\setminus U$.
\end{remark}

\subsection{The real case.}

If $k$ is a number field, which is not totally complex,
 then $k$ contains only two roots of unity.
We describe the result analogous to Theorem \ref{complexsplit}
 in this case.

Suppose that $\mu_{m}=\{1,-1\}$.
We shall identify $V_{\infty}$ with $\R^{d}$ for purposes
 of notation.
With this identification, $G_{\infty}$ is the group $\GL_{d}(\R)$.
We may choose $\gS$ to be a triangulation 
 of the unit sphere $S^{d-1}$ in $\R^{d}$.
Our fundamental domain $\scD$ may be taken to be the cell:
$$
 \scD
 =
 \left\{
 \left(\matrix{x_{1}\cr\vdots\cr x_{d}}\right)\in S^{d-1} :
 x_{1}\ge 0
 \right\}.
$$
Thus $\scE$ can be taken to be the cell
$$
 \scE
 =
 \left\{
 \left(\matrix{x_{1}\cr\vdots\cr x_{d}}\right)\in S^{d-1} :
 x_{1}= 0, \; x_{2}\ge 0
 \right\}.
$$
The cell $\scE$ is contained in the following subspace
 $W = \left\{
 \left(\matrix{0\cr x_{2}\cr\vdots\cr x_{d}}\right)
 \right\}$.
Define
$$
 U^{\prime}
 =
 \left\{\alpha\in\GL_{n}(\R):
 \begin{array}{c}
     \hbox{ $\alpha$ has no eigenvector in $W$}\\
     \hbox{with a negative real eigenvalue}
 \end{array}
 \right\}.
$$
The set $U^{\prime}$ is a dense, open, contractible subset of $\GL_{n}(\R)$.
One may prove analogously to the totally complex case:

\begin{theorem}
    If $\alpha$, $\beta$ and $\alpha\beta$ are all in $U^{\prime}$
     then $\deci^{(\scD)}$ is locally constant at $(\alpha,\beta)$.
\end{theorem}

This shows that $\deci^{(\scD)}$ is the cocycle
 obtained from a section
 $G_{\infty}\to \tilde G_{\infty}$,
 which takes the identity to the identity
 and is continuous on $U'$.

\section{Construction of fundamental functions.}

Let $k$ be an algebraic number field
 containing a primitive $m$-th root of unity
 and consider the vector space $V=k^{n}$.
As before, we let $S$ be the set of places $v$ of $k$ such that
 $|m|_{v}\ne 1$.
We define $V_{\AS}=\AS^{n}$ and $V_{\infty}=k_{\infty}^{n}$,
 where $k_{\infty}=k\otimes_{\Q}\R$.
From the previous section
 we have an arithmetic cocycle $\decAS^{(f,L)}$ on $\GLnAS$
 and a geometric cocycle $\deci^{(\scD)}$ on $\GLn(k_\infty)$.
We shall relate the two.
However the arithmetic cocycle is dependent on a choice
 of fundamental function $f$ on $\AS^n\setminus\{0\}$
 and the geometric cocycle is dependent (in the real case at least)
 on a fundamental domain $\scD$.
In order to describe the relation between $\decAS$ and $\deci$,
 we must first fix our choices of $f$ and $\scD$.
In this section we choose a fundamental function $f$
 (at least generically) and a related fundamental domain $\scD$
 and describe the mechanism by which the two cocycles will be related.
In section 5 we deal with the problem of defining $f$ everywhere,
In section 6 we prove the relation between
 the cocycles,
 based on the methods of this section.

From now on we shall assume that $m$ is a power
 of a prime $p$.
There is no loss of generality here,
 but we will save on notation by doing this.
We shall write $\rho$ for a primitive $p$-th root of unity in $k$.

\subsection{The space $X$.}

We fix once and for all a lattice
 $L\subset V$ which is free as a $\Z[\mu_{m}]$-module.
For any finite place $v$ we shall write $L_v$
 for the closure of $L$ in $V_v$.
We shall also write $L_{\AS}$ for the closure of $L$
 in $V_{\AS}$.
Hence $L_{\AS}=\prod_{v\notin S} L_{v}$.

Let
$$
 V_m
 =
 V
 \cap
 \bigcap_{v|m} L_v,
$$
 and consider the group $X=V_{m}/L$.
There are two ways of thinking about $X$.
First, the diagonal embedding of $k$ in $\AS$ induces an isomorphism
\begin{equation}
    \label{Xism}
    X \cong V_{\AS}/L_{\AS}.
\end{equation}
Secondly, we can regard $X$ as a dense subgroup of the group
 $X_{\infty}=V_{\infty}/L$.
The two ways of thinking about $X$ give the connection
 between the arithmetic and geometric cocycles defined
 in \S3.

\paragraph{The semi-group $\Upsilon$.}
Consider the following semigroup in $\GL_{n}(k)$:
$$
 \Upsilon
 =
 \{\alpha\in \GL_{n}(k) :
 \alpha L \supseteq L\}.
$$

Let $f$ be a fundamental function on $X\setminus\{0\}$,
 and define a fundamental function $f_{\AS}$ on
 $V_{\AS}\setminus\{0\}$ by
$$
 f_{\AS}(v)
 =
 \left\{
 \begin{array}{ll}
     f(v+L_{\AS}) & v\notin L_{\AS},\\
     f_{o}(v) & v\in L_{\AS},
 \end{array}
 \right.
$$
 where $f_{o}$ is any fixed fundamental function.
By abuse of notation we shall write $L$ and $L_{\AS}$ for the
 characteristic functions of these sets.

\begin{lemma}
    \label{Xism3}
    With the above notation we have for $\alpha,\beta\in\Upsilon$:
    $$
     \decAS^{(f_{\AS},L_{\AS})}(\alpha,\beta)
     =
     \langle f\alpha^{-1}-f | \alpha\beta L - \alpha L\rangle_{X}.
    $$
\end{lemma}

\proof
This follows immediately from Lemma \ref{Xism4}.
\sq
\medskip


\subsection{The complex $\gX$.}

Our method of calculating $\decAS(\alpha,\beta)$
 on $\Upsilon$ involves constructing fundamental functions on
 $X$ quite explicitly by embedding $X$ in $X_{\infty}:=V_{\infty}/L$
 and finding a fundamental domain $\scF$ for the
 action of $\mu_{m}$ on $X_{\infty}$.
In this section we will find the fundamental domain $\scF$.

\subsubsection{Parallelotopes and $\Diamond$-products.}

Let $\scT$ be a singular  $r$-cube and $\scU$ a
 singular  $s$-cube in $V_\infty$.
We can define an $(r+s)$-cube $\scT\Diamond \scU$ by:
$$
 (\scT\Diamond\scU)(x_1,...,x_r,y_1,...,y_s)
 =
 \scT (x_1,...,x_r) + \scU (y_1,...,y_s).
$$
Note that for any $v\in V_{\infty}$,
 the cell $[v]\Diamond\scT$ is simply a
 translation of $\scT$ by $v$.
Let $v,a_{1},\ldots,a_{r}\in V_{\infty}$.
By the \emph{parallelotope} $\parat(v,a_{1},\ldots,a_{r})$ in
 $V_{\infty}$
 we shall mean the following cell $I^{r}\to X_{\infty}$:
$$
 \Big(\parat(v,a_{1},\ldots,a_{r})\Big)(\xx)
 =
 v+\sum_{i=1}^r x_{i}a_i.
$$
Hence this is just a $\Diamond$-product of line-segments.
We shall more often deal with
 the projections $\para$ of $\parat$ in $X_{\infty}$.
We do not assume that the vectors $a_{i}$
 are linearly independent or even non-zero.

\subsubsection{Construction of $\gX$.}

We shall require a cell decomposition $\gX$ of $X_{\infty}$
 in which the cells are parallelotopes.
To describe this cell decomposition
 it is sufficient to describe the
 highest dimensional cells.

We begin with a cell decomposition
 of $\Q(\rho)_{\infty}/\Z[\rho]$.
The highest dimensional cells are of the form
$$
 \scP
 =
 \para\left(0,\frac{\rho^{i}}{1-\rho},\frac{\rho^{i+1}}{1-\rho},
 \ldots,\frac{\rho^{i+p-2}}{1-\rho}\right),
 \quad
 i=1,\ldots,p.
$$

\begin{lemma}
    The cells $\scP$ above form the highest dimensional cells of
    a cell decomposition of $\Q(\rho)_{\infty}/\Z[\rho]$.
\end{lemma}

\proof
This is Theorem 1.1 of \cite{kubota2}.
\sq
\medskip

We shall refer to the corresponding cell decomposition
 of $\Q(\rho)_{\infty}/\Z[\rho]$ as $\gX(p)$.

We next introduce a cell decomposition
 of $\Q(\zeta)_{\infty}/\Z[\zeta]$.
We have a decomposition
$$
 \Q(\zeta)_{\infty}/\Z[\zeta]
 =
 \bigoplus_{i=1}^{m/p}
 \zeta^{i} \cdot \Q(\rho)_{\infty}/\Z[\rho].
$$
We may therefore define a cell decomposition
 $\gX(p^a)$ of $\Q(\zeta)_{\infty}/\Z[\zeta]$ by
$$
 \gX(p^a)
 =
 \prod_{i=1}^{m/p}\zeta^{i}\cdot\gX(p).
$$
Thus the cells of $\gX(p^a)$ are of the form
 $\bigDiamond_{i=1}^{m/p} \zeta^i \scP_i$
 with $\scP_i$ a cell of $\gX(p)$.

As we are assuming that $L$ is free over $\Z[\zeta]$,
 there is a basis $\{b_{1},\ldots,b_{n}\}$ for $V$ over $\Q(\zeta)$
 such that $L=\sum_{i=1}^{n} \Z[\zeta]b_{i}$.
Again we have a decomposition
$$
 X_\infty
 =
 \bigoplus_{i=1}^{n}(\Q(\zeta)_{\infty}/\Z[\zeta])\cdot b_i.
$$
We may therefore define
$$
 \gX
 =
 \prod_{i=1}^{n}\gX(p^a)\cdot b_i.
$$

\begin{lemma}
\label{constructX}
    The cell decomposition $\gX$ of $X_{\infty}$
     has the following properties:
    \begin{itemize}
	\item[(i)]
	The group $\mu_{m}$ permutes the cells of $\gX$.
	\item[(ii)]
	Every positive dimensional cell has trivial
         stabilizer in $\mu_m$.
	\item[(iii)]
	Every $r$-cell $\scP$ in $\gX$ is of the form
	$$
	 \scP
	 =
	 \pr \scPt,\quad
	 \scPt
	 =
	 \parat(v_\scP,a_{\scP,1},\ldots,a_{\scP,r}),
	$$
	 with $v_\scP\in \frac{1}{1-\rho}L$ and
	 $a_{\scP,i}\in \frac{1}{1-\rho}L\setminus L$.
        For any $\scP$
	 the set $\{a_{\scP,1},\ldots,a_{\scP,r}\}$
	 is linearly independent over $\R$.
	\item[(iv)]
	We have $|\scPt|\cap L \subseteq\{0\}$
	 and $v_{\scP}=0$ if and only if $0\in|\scPt|$.
    \end{itemize}
\end{lemma}

Much of this result is contained in Theorem 1.3 of \cite{kubota2},
 although it is stated there in a rather different notation.
A sketch of the proof is included for completeness.
\medskip

\proof
(i) and (iii) are clear from the construction.
It is sufficient to verify (iv) for $\gX(p)$
 and this is not difficult.
It remains to show that no positive dimensional cell
 is fixed by a non-trivial subgroup of $\mu_m$.
Let $\scP$ be a positive dimensional cell
 and suppose $\xi\scP=\scP$ for some $\xi\in\mu_{m}\setminus\{1\}$.
We have an expression for $\scP$ as a parallelotope:
$$
 \scP
 =
 \para(v_{\scP},a_{\scP,1},\ldots,a_{\scP,r}).
$$
The set $\{a_{\scP,1},\ldots,a_{\scP,r}\}$ is permuted by $\xi$.
This implies that the elements $\xi^{i}a_{\scP,1}$ are all in this
 set.
However, the sum of these elements is zero.
This contradicts the fact (iii),
 that the $a_{\scP,i}$ are linearly independent.
\sq
\medskip

We shall write $\gX_{\bullet}$ for the corresponding chain complex with
 coefficients in $\Z/m$.
Thus $\gX_{r}$ is the free $\Z/m$-module on the $r$-cells of the
 decomposition.
By part (i) of the lemma, we have an action of $\mu_{m}$ on 
 $\gX_{\bullet}$.

\begin{lemma}
    For $r=1,2,\ldots,d$ the $(\Z/m)[\mu_{m}]$-module $\gX_{r}$
     is free.
\end{lemma}

\proof
A basis consists of a set of representatives
 for the $\mu_{m}$-orbits of $r$-cells in $\gX$.
To show that this is a basis 
 we use the fact that such cells have trivial stabilizer.
\sq
\medskip

\subsubsection{The fundamental function $f$.}

We shall fix an orientation $\ord_{V}$ on $V_{\infty}$.
Using this, we define a corresponding orientation $\ord_{X}$
 on $X_{\infty}$ by the formula
$$
 \ord_{X,x}(\pr(\scT))
 =
 \sum_{v\in V_{\infty} :\\ \pr(v)=x}
 \ord_{V,v}(\scT).
$$
Let $\omega_{X}\in \gX_{d}$ be the generator
 of $H_{d}(\gX)$, for which $\ord_{X,x}(\omega_{X})=1$
 holds for every $x\in X_{\infty}$.
By Lemma \ref{orient}, $(1-[\zeta])\omega_{X}=0$
 holds in $H_{d}(\gX)$.
Since $\gX_{d+1}=0$, we have $(1-[\zeta])\omega_{X}=0$
 in $\gX_{d}$.
Hence by the exact sequence (\ref{star})
 there is an element $\scF\in \gX_{d}$ such that
 $[\mu_{m}]\scF = \omega_{X}$.
We fix such an $\scF$ once and for all.
We may regard $\scF$ as a fundamental domain for
 the action of $\mu_{m}$ on $X_{\infty}$.

Define a function $f:X_{\infty}\setminus |\partial\scF| \to \Z/m$
 by
$$
 f(x)
 =
 \ord_{x}(\scF).
$$

\begin{lemma}
    The function $f$ is fundamental on the set of $x\in X_{\infty}$
     whose $\mu_{m}$-orbit does not meet the boundary of $\scP$.
\end{lemma}

\proof
For such an $x$, we have by Lemma \ref{orient}:
$$
 \sum_{\xi\in\mu_{m}} f(\xi x)
 =
 \sum_{\xi\in\mu_{m}} \ord_{\xi x}\scF
 =
 \sum_{\xi\in\mu_{m}} \ord_{x}\left(\xi^{-1}\cdot\scF\right).
$$
By linearity of $\ord_{x}$ we have:
$$
 \sum_{\xi\in\mu_{m}} f(\xi x)
 =
 \ord_{x}\left([\mu_{m}]\scF\right)
 =
 \ord_{x}\left(\omega_{X}\right)
 =
 1.
$$
\sq
\medskip

We shall worry about how to extend the definition
 of $f$ to $|\partial\scF|$ in \S5.3.
The solution will be to take a limit over fundamental
 domains tending to $\scF$.

\subsection{The complex $\gS$.}

Now that we have a fundamental function $f$ at least generically,
 we shall describe the fundamental domain $\scD$ and the cell
 complex $\gS\subset V_{\infty}\setminus\{0\}$
 used in the definition of $\dec_\infty^{(\scD)}$
 in \S3.5.

We have a cell decomposition $\gX$ of $X_\infty$.
Each $r$-cell in this decomposition is of the form
$$
 \scP
 =
 \para(v_\scP,a_{\scP,1},\ldots,a_{\scP,r}).
$$
Corresponding to each such cell,
 we define an $r-1$-chain
 $\gs(\scP)\in C_{r-1}(V_\infty\setminus\{0\})$ as follows.
If $v_\scP\ne 0$ then we simply define $\gs\scP=0$.
If $v_\scP=0$
 then we define $\gs\scP:\Delta^{r-1}\to V_\infty\setminus\{0\}$
 by $\gs\scP=[a_{\scP,1},\ldots,a_{\scP,r}]$.
Roughly speaking $\gs\scP$ is the set of unit tangent vectors
 to $\scP$ at $0$.

The cells $\gs\scP$ form a cell complex,
 which we shall denote $\gS$.
It follows easily that $\gS$ satisfies the conditions of \S3.5.
We extend $\gs$ by linearity to a map
 $\gs:\gX_\bullet\to \gS_{\bullet-1}$.

\begin{lemma}
\label{sigma}
 The map $\gs:\gX_\bullet\to \gS_{\bullet-1}$
  commutes with the action of $\mu_{m}$
  and anticommutes with $\partial$.
 That is, for any cell $\scP$ of $\gX$, we have
 $$
  \partial \gs \scP
  =
  -\gs\partial\scP.
 $$
(Here the minus sign is acting on the coefficients,
 rather than on $V_{\infty}$).
\end{lemma}

\proof
The first statement is clear.
For the second we need to
 examine some separate cases.
Let $\scP$ be an $r$-cell in $\gX$.
First suppose $v_\scP\ne 0$.
In this case, we know by part (iv) of Lemma \ref{constructX},
 that 0 is not in the base set of $\scP$.
It follows that $0$ is not in the base set of any face
 $\scQ$ of $\scP$.
Hence $v_\scQ\ne0$ for every face $\scQ$ of $\scP$.
We therefore have $\gs\scP=0$ and $\gs\scQ=0$ for every face $\scQ$.
As $\partial\scP$ is a sum of faces, the result follows in this case.

Now suppose $v_\scP=0$.
Let $\scQ$ be a face of $\scP$.
If $\scQ$ is a front face then we have $v_\scQ=0$
 and if $\scQ$ is a back face we have $v_\scQ\ne 0$.
Thus, when calculating $\gs\partial\scP$,
 we need only take into account the front faces of $\scP$.
It follows that we have
$$
 \gs\partial\scP
 =
 \sum_{i=1}^r
 (-1)^{i}
 \gs\left(\bigDiamond_{j\ne i} [0,a_{\scP,j}]\right).
$$
By the definition of $\gs$, this gives
\begin{eqnarray*}
 \gs\partial\scP
 &=&
 \sum_{i=1}^r
 (-1)^{i}
 [a_{\scP,1},\ldots,a_{\scP,i-1},a_{\scP,i+1},\ldots,a_{\scP,r}]\\
 &=&
 -\partial[a_{\scP,1},\ldots,a_{\scP,r}]
 =
 -\partial\gs\scP.
\end{eqnarray*}
\sq
\medskip

We define $\scD=\gs\scF$,
 where $\scF$ is the fundamental domain in $\gX_d$ chosen
 in \S4.2.

\begin{lemma}
    The element $\scD$ satisfies the conditions of
     \S3.5.
\end{lemma}

\proof
By Lemma \ref{sigma} we have:
$$
 \partial[\mu_{m}]\scD
 =
 \partial \gs[\mu_{m}]\scF
 =
 \partial \gs\omega_{X}
 =
 -\gs\partial \omega_{X}
 =
 0.
$$
Hence $[\mu_{m}]\scD$ is a cycle. 
It remains to check
 that the image of $[\mu_{m}]\scD$ under the maps
$$
 H_{d-1}(V_\infty\setminus\{0\})
 \stackrel{\partial^{-1}}{\to}
 H_{d}(V_\infty,V_\infty\setminus\{0\})
 \stackrel{\ord_0}{\to}
 \Z/m,
$$
 is 1.
For any cell $\scP$ of $\gX$ centred at 0 we shall
 use the notation
$$
 \gt\scP
 =
 [0,a_{\scP,1},\ldots,a_{\scP,r}].
$$
We extend $\gt$ by linearity to a map
 $\gX_{\bullet}\to C_{\bullet}(V_{\infty})$.
With this notation we have:
$$
 \partial(\gt\scP)
 =
 \gs\scP - \gt(\partial\scP).
$$
Applying this relation to $[\mu_m]\scF$,
 we obtain by Lemma \ref{sigma}:
$$
 \partial([\mu_m]\gt\scF)
 =
 [\mu_m]\scD.
$$
We are therefore reduced to showing
 that $\ord_{V,0}([\mu_m]\gt\scF)=\ord_{X,0}([\mu_m]\scF)$.
This is a little messy, but it can be proved by induction
 on the dimension $d$ of $V_{\infty}$.
\sq
\medskip


\subsection{Modified Parallelotopes.}

In this section we shall discuss a method for
 constructing more general fundamental functions on $X_{\infty}$.

We have a cell decomposition $\gX$ of $X_{\infty}$ in which the $r$-cells
 are of the form
$$
 \scP
 =
 \para(v_\scP,a_{\scP,1},\ldots,a_{\scP,r}).
$$
Recall that $\gg=\End_{\mu_m}(V_{\infty})$.
We shall write $1$ for the identity matrix in $\gg$.
Suppose $\path$ is a path from $0$ to $1$ in $\gg$
 and $\scP$ is an $r$-cell in the complex $\gX$.
We define an $r$-cell $\path\bowtie \scP:I^r\to X_\infty$ as follows:
$$
 (\path\bowtie\scP)(x_{1},\ldots,x_{r})
 =
 \pr\left(
 v_{\scP}+\sum_{i=1}^{r} \path(x_{i}) \cdot a_{\scP,i}
 \right).
$$
We extend the operator ``$\path\bowtie$'' by linearity to
 a map $(\path\bowtie) : \gX_{r} \to C_{r}(X_{\infty})$.
Following Kubota, \cite{kubota2} we shall refer to
 $\path\bowtie\scP$ as a \emph{modified parallelotope}.

\begin{lemma}
    The maps $(\path\bowtie) : \gX_{r} \to C_{r}(X_{\infty})$
     commute with the boundary maps and with the
     action of $\mu_{m}$.
\end{lemma}

\proof
This is a routine verification.
\sq
\medskip

Suppose $\path_{1}$ and $\path_{2}$ are two paths from 0 to 1
 in $\gg$ and $\scH$ is a homotopy from $\path_{1}$ to $\path_{2}$.
By this,
 we shall mean $\scH:I^2\to\gg$ is a continuous map,
 satisfying for all $t,x\in I$:
$$
 \scH(0,x)=\path_1(x),\quad
 \scH(1,x)=\path_2(x),\quad
 \scH(t,0)=0,\quad
 \scH(t,1)=1.
$$
Let $\scP$ be an $r$-cell in the complex $\gX$:
$$
 \scP
 =
 \para(v_\scP,a_{\scP,1},\ldots,a_{\scP,r}).
$$
We define an $r+1$-cell $(\scH\bowtie\scP):I^{r+1}\to V_{\infty}$
 by
$$
 (\scH\bowtie\scP)(t,x_{1},\ldots,x_{r})
 =
 \pr\left(
 v_{\scP}+ \sum_{i=1}^{r} \scH(t,x_{i}) \cdot a_{\scP,r}
 \right).
$$
We extend the operators ``$\scH\bowtie$'' to linear maps
 $(\scH\bowtie):\gX_{r}\to C_{r+1}(X_{\infty})$.

\begin{lemma}
    The maps $\scH\bowtie$ commute with the action of $\mu_{m}$.
    Let $\scH$ be a homotopy from $\path_{1}$ to $\path_{2}$.
    Then $\scH\bowtie$ is a chain homotopy
     from $\path_{1}\bowtie$ to $\path_{2}\bowtie$.
    In other words, for any $\scP\in \gX_{r}$ we have in
     $C_{r}(X_{\infty})$:
    $$
     \partial(\scH\bowtie\scP)+\scH\bowtie\partial\scP
     =
     \path_{2}\bowtie\scP-\path_{1}\bowtie\scP.
    $$
\end{lemma}

\proof
This is proved by calculating $\partial(\scH\bowtie\scP)$
 directly for a cell $\scP$ of $\gX$.
\sq
\medskip

\begin{lemma}
    For any path $\path$, the following holds in $H_{d}(X_{\infty})$:
    $$
     [\mu_m]\path\bowtie\scF
     =
     [\mu_m]\scF.
    $$
\end{lemma}

(This lemma is implicit in \cite{kubota2}).
\medskip

\proof
This follows immediately from the previous two lemmata.
\sq
\medskip

We define a function
 $f^{\path}:X_{\infty}\setminus |\path\bowtie\partial\scF|\to\Z/m$ by
$$
 f^{\path}(x)
 =
 \ord_{x}(\path\bowtie\scF).
$$
It follows from the previous lemma,
 that $f^{\path}$ is fundamental at all $x$, whose
 $\mu_{m}$-orbit avoids the base set of $\path\bowtie\partial\scF$.
However $f^{\path}$ need not be
 the characteristic function of a fundamental domain,
 since it may take other values apart from $0$ and $1$.

\paragraph{Equivalence of paths.}
Suppose $\path$ and $\path^{\prime}$ are two
 paths from $0$ to $1$ in $\gg$.
We shall say that $\path$ and $\path^{\prime}$
 are \emph{equivalent} if there is an increasing continuous
 bijection $\phi:I\to I$,
 such that for all $x\in I$ we have:
$$
 \path(\phi(x))
 =
 \path^{\prime}(x).
$$

\begin{lemma}
    If $\path$ and $\path^{\prime}$ are equivalent paths
    then the corresponding functions $f^{\path}$ and
    $f^{\path^{\prime}}$ have the same domains of definition and are equal.
\end{lemma}

\proof
To prove that $f^{\path}$ and $f^{\path^{\prime}}$
 have the same domains of definition,
 we need only show that
 $|\partial(\path\bowtie\scF)|=|\partial(\path^{\prime}\bowtie\scF)|$.
This follows by breaking the boundary into faces and
 using the fact that $\phi$ is bijective.

The map $\phi$ is homotopic to the identity map.
In other words there is a map $\psi:I\times I\to I$
 such that $\psi(1,x)=\phi(x)$, $\psi(0,x)=x$,
 $\psi(t,0)=0$ and $\psi(t,1)=1$ for all $x,t\in I$.
Using the function $\psi$ we define a homotopy $\scH$
 from $\path$ to $\path^\prime$ by $\scH(t,x)=\path(\psi(t,x))$.
As $\scH\bowtie$ is a chain homotopy we have
$$
 \path^{\prime}\bowtie\scF-\path\bowtie\scF
 =
 \partial(\scH\bowtie\scF)
 +
 \scH\bowtie\partial\scF.
$$
Note that for any cell $\scP$ of $\gX$ we have
 $|\scH\bowtie\scP|=|\path\bowtie\scP|=|\path^\prime\bowtie\scP|$.
Therefore for $x\notin |\path\bowtie\partial\scF|$ we have:
$$
 \ord_x(\path^{\prime}\bowtie\scF)-\ord_x(\path\bowtie\scF)
 =
 0.
$$
In other words $f^{\path}(x) = f^{\path^{\prime}}(x)$.
\sq
\medskip

In view of the above lemma,
 we may specify piecewise linear paths
 simply as sequences of line segments:
$$
 \path
 =
 [0,a_1]+[a_1,a_2]+\cdots+[a_s,1],
$$
 without worrying about the precise parametrization.

\subsection{Statement of the results in the generic case.}

\paragraph{The $d-1$-chain $\scG$.}
For the moment we shall assume that $d\ge 2$.
Thus we are ruling out the case $k=\Q$, $n=1$.
By the definition of $\scF$, we have
$$
 \partial([\mu_m] \scF)
 =
 0.
$$
Thus $[\mu_m](\partial\scF)=0$.
It follows from the exact sequence (\ref{star}),
 that there is an element $\scG\in \gX_{d-1}$ satisfying
$$
 \partial \scF
 =
 (1-[\zeta])\scG.
$$
We shall fix such a $\scG$.
Note that by Lemma \ref{sigma},
 the $d-2$-chain $\scE$ used in the definition of $\deci$
 may be taken to be $-\gs\scG$.

\paragraph{The semigroup $\Upsilon_{\gf}$.}
For an ideal $\gf$ of the ring of integers in $k$,
 let $G_{\gf}$ be the subgroup of $\GL_{n}(k)$ consisting
 of matrices which are integral at every prime dividing $\gf$
 and are congruent to the identity matrix modulo $\gf$.
We shall fix $\gf=(1-\rho)m^{2}$ if $m$ is odd and
 $\gf=4m^{2}$ if $m$ is even.
Next let $\Upsilon_{\gf}=\Upsilon\cap G_{\gf}$.

\paragraph{The results.}
For $\alpha\in\Upsilon_{\gf}$, we define a path
$$
 \path(\alpha)
 =
 [0,\alpha]+[\alpha,1],
$$
 and a homotopy $\scH^1_\alpha$ from
 $\path(1)$ to $\path(\alpha)$:
$$
 \scH^1_\alpha(t,x)
 =
 (1-t)x+t\path(\alpha)(x).
$$
Finally we define
$$
 \tau(\alpha)
 =
 \zeta^{\textstyle{\{\scH^{1}_{\alpha}\bowtie\scG|\alpha L\}}}.
$$
 (One may interpret this definition as a lifted map
 in the sense of Lemma \ref{easy})
 
We shall prove that $\deci^{(\scD)}$ and $\decAS^{(f,L)}$ are related on
 $\Upsilon_{\gf}$ by the following formula:
$$
 \decAS^{(f,L)}(\alpha,\beta)
 \deci^{(\scD)}(\alpha,\beta)
 =
 \frac{\tau(\alpha)\tau(\beta)}{\tau(\alpha\beta)},
 \quad
 \alpha,\beta\in\Upsilon_{\gf}.
$$
As a consequence, we will deduce that for totally complex $k$,
 the Kubota symbol on $\SL_n(\go,\gf)$ is given by
$$
 \kappa(\alpha)
 =
 \frac{\tau(\alpha)}{hw(\alpha)}
$$
Here $hw$ is the splitting of $\deci$
 of \S3.9, Theorem \ref{complexsplit}.

\subsection{A generic formula for the pairing.}

Given paths $\path^{1}$ and $\path^{2}$ from $0$ to $1$ in $\gg$,
 we have constructed fundamental functions
 $f^{1}$ and $f^{2}$ in $\cF_{X}$.
We now describe a geometric method
 for calculating the pairing $\langle f^1-f^2|M-L\rangle_{X}$,
 where $M\supset L$ is a lattice contained in $V_{m}$.

\begin{proposition}
    \label{productformula}
    Suppose $M\setminus L$ does not intersect the
     base set of $\partial(\scH\bowtie\scG)$.
    Then we have:
    $$
     \langle f^{\path_{2}}-f^{\path_{1}}|M-L\rangle_{X}
     =
     \zeta^{\textstyle{\{\scH\bowtie\scG|M-L\}}}.
    $$
\end{proposition}

The right hand side is a lifted map in the sense
 of Lemma \ref{easy}.
\medskip

\proof
As $\scH\bowtie$ is a chain homotopy from $\path_{1}\bowtie$ to
 $\path_{2}\bowtie$, the following holds in $C_{d}(X_{\infty})$:
$$
 \path_{2}\bowtie\scF - \path_{1}\bowtie\scF
 =
 \partial(\scH\bowtie\scF)+\scH\bowtie\partial\scF.
$$
This implies for $x\in M\setminus L$:
$$
 f^{\path_{2}}(x)-f^{\path_{1}}(x)
 =
 \ord_{x}(\path_{2}\bowtie\scF - \path_{1}\bowtie\scF)
 =
 \ord_{x}(\scH\bowtie\partial\scF).
$$
By the definition of $\scG$ and Lemma \ref{orient}, we have:
$$
 f^{\path_{2}}(x)-f^{\path_{1}}(x)
 =
 \ord_{x}(\scH\bowtie(1-[\zeta])\scG)
 =
 \ord_{x}(\scH\bowtie\scG)
 -
 \ord_{\zeta^{-1}x}(\scH\bowtie\scG).
$$
The proposition now follows from the
 definition (\ref{simplepairing}) in \S3.2
 of the pairing $\langle-|-\rangle_{X}$.
\sq
\medskip

\subsection{The order of $\scH\bowtie\scG$ at $0$.}

The statement of the results in \S4.5 involves
 numbers of the form $\ord_{0}(\scH\bowtie \scG)$.
However, this is not as yet defined,
 since $0$ is always in the
 base set of $\partial(\scH\bowtie \scG)$.
To avoid this problem, for any $r$-cell $\scP$
 in $\gX$ containing $0$
 we cut $\scH\bowtie\scP$ into a singular
 part $\scH\bowtie\scP^{0}$ and a non-singular part $\scH\bowtie\scP^{+}$.

Let $U$ be a small neighbourhood of $0$ in $I^{d-1}$.
We define $\scH\bowtie\scP^{0}$ to be the restriction of
 $\scH\bowtie\scP$ to $I\times U$ and we define
 $\scH\bowtie\scP^{+}$ to be the restriction of
 $\scH\bowtie\scP$ to $I\times (I^{d-1}\setminus U)$.
We define the order of $\scH\bowtie\scP$ at 0 to be the limit as $U$
 gets smaller of the order of $\scH\bowtie\scP^{+}$ at 0.

To make things a little more precise we define for $\epsilon>0$
 an $r$-cell $\scP^\epsilon$ in $X_{\infty}$.
If $v_{\scP}=0$ then we define $\scP^{\epsilon}$
 to be the restriction of $\scP:I^r\to X_\infty$
 to the set
$$
 \left\{(x_1,\ldots,x_r)\in I^r: \sum x_i \le \epsilon\right\}.
$$
If $v_{\scP}=0$ then we define $\scP^{\epsilon}=0$.

\begin{lemma}
\label{Pepsilon}
 For $\epsilon>0$ sufficiently small we have
 $$
  \partial(\scP^\epsilon)-(\partial\scP)^\epsilon
  =
  \pr( \epsilon\cdot\gs\scP ).
 $$
\end{lemma}

\proof
We shall suppose that $\scP$ has 0 as its origin,
 since otherwise both sides of the formula are zero.
Under this assumption, we have:
$$
 \scP
 =
 \bigDiamond_{i=1}^r
 [0,a_{\scP,i}],
 \quad
 \scP^{\epsilon}
 =
 [0, \epsilon a_{\scP,1},\epsilon a_{\scP,2},\ldots,\epsilon a_{\scP,r}].
$$
Therefore
\begin{eqnarray*}
 \partial(\scP^{\epsilon})
 &=&
 [\epsilon a_{\scP,1},\ldots,\epsilon a_{\scP,r}]
 +
 \sum_{i=1}^r
 (-1)^i
 [0, \epsilon a_{\scP,1},\ldots,
  \epsilon a_{\scP,i-1},\epsilon a_{\scP,i+1}\ldots,
  \epsilon a_{\scP,r}]\\
 &=&
 [\epsilon a_{\scP,1},\ldots,\epsilon a_{\scP,r}]
 +
 (\partial\scP)^\epsilon.
\end{eqnarray*}
The result follows.
\sq
\medskip

\subsection{Dependence of $\ord_{0}(\scH\bowtie\scG)$ on $\scH$.}

Let $\path_{1}$ and $\path_{2}$ be two paths from $0$ to
 $1$ in $\gg$,
Suppose we have two homotopies $\scH$ and
 $\scI$ from $\path_{1}$ to $\path_{2}$.
In this section, we investigate the relation between
 $\ord_{0,X}(\scH\bowtie\scG)$ and $\ord_{0,X}(\scI\bowtie\scG)$.

We begin by choosing a homotopy $\scU$ from $\scH$ to $\scI$.
Thus, $\scU:I^{3}\to \gg$ satisfies the following conditions:
$$
 \scU(0,t,x)=\scH(t,x),\quad
 \scU(1,t,x)=\scI(t,x),
$$
$$
 \scU(u,0,x)=\path_{1}(x),\quad
 \scU(u,1,x)=\path_{2}(x),
$$
$$
 \scU(u,t,0)=0,\quad
 \scU(u,t,1)=1.
$$
We shall also suppose that there is some
 $\epsilon>0$ such that for $x<1$ we have
$$
 \scU(u,t,\epsilon x)
 =
 x\scU(u,t,\epsilon).
$$
We shall define, under this assumption, maps
 $(\gs\scU):I^2\to\gg$, $(\gs\scH):I\to\gg$
 and $(\gs\scI):I\to\gg$ by
$$
 (\gs\scU)(u,t)
 =
 \epsilon^{-1}\scU(u,t,\epsilon)
 =
 \frac{\partial}{\partial x} \scU(u,t,x),
$$
$$
 (\gs\scH)(t)
 =
 \epsilon^{-1}\scH(t,\epsilon),\quad
 (\gs\scI)(t)
 =
 \epsilon^{-1}\scI(t,\epsilon).
$$
We define maps $(\scU\bowtie):\gX_{r}\to C_{r+2}(X_{\infty})$ by
$$
 (\scU\bowtie\scP)(u,t,x_{1},\ldots,x_{r})
 =
 \pr\left(
 v_{\scP}
 +
 \sum_{i=1}^{r}\scU(u,t,x_{i})\right).
$$
The next lemma says that $\scU\bowtie$ is in some sense
 a chain homotopy between $\scH\bowtie$ and $\scI\bowtie$
 (although these are themselves chain homotopies, not
 chain maps).

\begin{lemma}
\label{Ubowtie}
    The map ``$\scU\bowtie$'' commutes with the action of $\mu_{m}$.
    For every $\scP\in\gX_{r}$
     the following holds in $C_{r+1}(X_{\infty})$:
    $$
     \scI\bowtie\scP - \scH\bowtie\scP
     =
     \partial(\scU\bowtie\scP)
     -
     \scU\bowtie\partial\scP.
    $$
\end{lemma}

\proof
By definition of the boundary map we have
$$
 \partial(\scU\bowtie\scP)
 =
 (\scI\bowtie\scP-\scH\bowtie\scP)
 -
 \hbox{degenerate cells}
 +
 \scU\bowtie(\partial\scP).
$$
\sq
\medskip

%

The crucial point in relating the two cocycles
 is the following formula.

\begin{proposition}
\label{XV}
Suppose that for every $d-1$-cell $\scP$ in $\gX$,
 $\ord_{0,X}(\scH\bowtie\scP^+)$
 and $\ord_{0,X}(\scI\bowtie\scP^+)$
 are defined and for every $d-2$-cell $\scP$ in $\gX$,
 $\ord_{0,X}(\scU\bowtie\scP^+)$ is defined.
Then so is $\ord_{V,0}(\gs\scU\cdot\gs\scG)$
 and we have:
$$
 \ord_{X,0}(\scI\bowtie\scG^{+})-\ord_{X,0}(\scH\bowtie\scG^{+})
 =
 \ord_{V,0}(\gs\scU\cdot\gs\scG).
$$
\end{proposition}

\proof
Applying lemma \ref{Ubowtie} to $\scG$ we obtain
$$
 \scI\bowtie\scG - \scH\bowtie\scG
 =
 \partial(\scU\bowtie\scG)
 -
 \scU\bowtie\partial\scG
$$
Recall that $\scG$ is chosen to satisfy the relation
 $(1-[\zeta])\scG=\partial\scF$.
This implies $(1-[\zeta])\partial\scG=0$.
By the exact sequence (\ref{star}),
 there is a $d-2$-chain $\scQ$ in $\gX$
 such that
$$
 \partial\scG=[\mu_{m}]\scQ.
$$
We therefore have
$$
 (\scI-\scH)\bowtie\scG
 =
 \partial(\scU\bowtie\scG)
 -
 [\mu_{m}](\scU\bowtie\scQ).
$$
We cannot define the order at $0$ of either side of the above equation.
We therefore break $\scG$ and $\scQ$ into their singular
 and non-singular parts.
This gives:
$$
 (\scI-\scH)\bowtie\scG^{+}
 =
 \partial(\scU\bowtie\scG)
 -
 [\mu_{m}]\scU\bowtie\scQ^{+}
 -
 [\mu_{m}]\scU\bowtie\scQ^{\epsilon}
 -
 (\scI-\scH)\bowtie\scG^{\epsilon}.
$$
Note that $\ord_0(\scU\bowtie\scQ^{+})$ is defined,
 so we have by Lemma \ref{orient} $\ord_0([\mu_m]\scU\bowtie\scQ^{+})=0$.
Hence the following holds in $\Z/m$:
$$
 \ord_{0,X}((\scI-\scH)\bowtie\scG^{+})
 =
 -\ord_{0,X}((\scU\bowtie\partial\scG)^{\epsilon}
 +(\scI-\scH)\bowtie\scG^{\epsilon}).
$$
As every cell on the right hand side is contained
 in a small neighbourhood of $0$,
 we may replace $\ord_{0,X}$ by $\ord_{0,V}$.

Now choose an $r$-cell $\scP$ in $\gX$.
We have for $x_i\le\epsilon$:
$$
 (\scU\bowtie\scP)(u,t,x_1,\ldots,x_r)
 =
 \sum_{i=1}^{r} \scU(u,t,x_i)a_{\scP,i}
 =
 \gs\scU(u,t) \sum_{i=1}^{r} x_i a_{\scP,i}.
$$
Thus
$$
 (\scU\bowtie\scP)^\epsilon=\gs\scU\cdot(\scP^\epsilon),
$$
 and similarly,
$$
 (\scH\bowtie\scP)^\epsilon
 =
 \gs\scH\cdot(\scP^\epsilon),
 \quad
 (\scI\bowtie\scP)^\epsilon
 =
 \gs\scI\cdot(\scP^\epsilon).
$$
We therefore have
$$
 \ord_{0,X}((\scI-\scH)\bowtie\scG^{+})
 =
 -\ord_{0,V}(\gs\scU\cdot(\partial\scG)^{\epsilon}
 +(\gs\scI-\gs\scH)\cdot\scG^{\epsilon}).
$$
On the other hand, modulo degenerate cells, we have :
$$
 \partial\gs\scU=\gs\scI-\gs\scH.
$$
This implies
$$
 \ord_{0,X}((\scI-\scH)\bowtie\scG^{+})
 =
 -\ord_{0,V}(\gs\scU\cdot(\partial\scG)^{\epsilon}
 +(\partial\gs\scU)\cdot(\scG^{\epsilon})).
$$
Adding $\partial(\gs\scU\cdot(\scG^{\epsilon}))$ we obtain:
$$
 \ord_{0,X}((\scI-\scH)\bowtie\scG^{+})
 =
 -\ord_{0,V}(
 \gs\scU\cdot (\partial\scG)^{\epsilon}
 -\gs\scU\cdot \partial(\scG^{\epsilon}) ).
$$
Now by Lemma \ref{Pepsilon} we have:
$$
 \ord_{0,X}((\scI-\scH)\bowtie\scG^{+})
 =
 \ord_{0,V}(\epsilon\gs\scU\cdot\gs\scG).
$$
The right hand side is however independent of $\epsilon$
 so the result follows.
\sq
\medskip

\begin{corollary}
    \label{XVcor}
    Suppose there is an $\epsilon>0$ such that for all $0<x<\epsilon$
     we have $\scH(t,x)=\scI(t,x)$.
    Then $\ord_{0}(\scH\bowtie\scG)=\ord_{0}(\scI\bowtie\scG)$.
    (By this we understand that if both sides are defined then they are equal).
\end{corollary}

\proof
in this case we can choose $\scU$
 so that $\gs\scU$ is degenerate.
\sq
\medskip

\section{Deformation of fundamental functions.}

Given a path $\path$,
 we have constructed a function $f^{\path}$ on
 $X_\infty \setminus |\partial(\path\bowtie\scP) |$,
 which is fundamental on all $\mu_{m}$-orbits
 which do not intersect $|\partial(\path\bowtie\scP)|$.
In this chapter we describe a method for extending $f^\path$
 to all but a finite number of points of $X_\infty$.
Recall that we have fixed an ordered basis $\{b_1,\ldots,b_{r}\}$
 for $\gg$ as a vector space over $\R$.
The general idea is that if a point $x\in X_{\infty}$
 is on the boundary of $\path\bowtie\scP$,
 then we move the path $\path$ a little in the direction of $b_{1}$;
 if $x$ is still on the boundary then we move the path in the
 direction $b_{2}$ and so on.
Thus, we define $f^\path$ as a limit of the
 form $\blim_{\epsilon\to 0^+}$ in the notation
 of \S3.5.

We begin with some general results on the existence
 of such limits, and then prove the specific results
 which we need.

\subsection{Existence of certain limits.}

We shall call a subset $Z\subset\R^{a}$
 a \emph{small} subset of $\R^{a}$, if there is
 a $b\ge 0$ and a Zariski-closed subset $Z^{\prime}\subset \R^{a+b}$
 of codimension $\ge b+1$,
 such that $Z$ is contained in
 the archimedean closure of the projection of $Z^{\prime}$
 in $\R^{a}$.
If $Z_{1},\ldots,Z_{c}$ are small subsets of
 $\R^{a}$ then $Z_{1}\cup\ldots\cup Z_{c}$ is a
 small subset of $\R^{a}$.
Note that if $Z$ is small then it has codimension at least 1
 in $\R^{a}$.

\begin{lemma}
    \label{zariskilimit}
    Let $Z$ be a small subset of $\R^{a}$.
    Then for any locally constant function
     $\psi:\R^{a}\setminus Z\to \Z$
     the following (ordered) limit exists:
    $$
     \lim_{\epsilon_{1}\to 0^{+}}
     \lim_{\epsilon_{2}\to 0^{+}}
     \ldots
     \lim_{\epsilon_{a}\to 0^{+}}
     \psi(\epsilon_{1},\ldots,\epsilon_{a}).
    $$
\end{lemma}

\proof
We shall prove this by induction on $a$.
When $a=1$ the set $Z^{\prime}\subset \R^{b+1}$ has codimension $\ge b+1$,
 and is therefore finite.
It follows that $Z$ is finite,
 and it is clear that the limit exists in this case.
Now suppose $a>1$.
We shall decompose the Zariski-closed set $Z^{\prime}\subset \R^{a+b}$
 into irreducible components:
$$
 Z^{\prime} = Z^{\prime}_{1}\cup\ldots\cup Z^{\prime}_{r}.
$$
We shall write $Z_{i}$ for the archimedean closure of the
 projection of $Z^{\prime}_{i}$ in $\R^{a}$.
Let $\{b_{1},\ldots,b_{a+b}\}$ be the standard basis of $\R^{a+b}$.
Write $H^{\prime}$ for the hyperplane in $\R^{a+b}$ spanned
 by $\{b_{1},\ldots,b_{a-1},\ b_{a+1},\ldots,b_{a+b}\}$
 and let $H=\Span\{b_{1},\ldots,b_{a-1}\}$
 be the projection of $H$ in $\R^{a}$.
For each component $Z^{\prime}_{i}$ of $Z^{\prime}$,
 we define a subset $W^{\prime}_{i}\subset H^{\prime}$ by
$$
 W^{\prime}_{i}
 =
 \left\{
 \begin{array}{ll}
     Z^{\prime}_{i}\cap H & \hbox{if $Z^{\prime}_{i}\not\subset H'$,}\\
     \emptyset   & \hbox{if $Z^{\prime}_{i}\subset H'$.}
 \end{array}
 \right.
$$
Thus $W^{\prime}:=W^{\prime}_{1}\cup\ldots \cup W^{\prime}_{r}$
 is Zariski-closed in $H^{\prime}$
 and has codimension $\ge b+1$ in $H^{\prime}$.
Let $W$ be the archimedean closure of the
 projection of $W^{\prime}$ in $H$.
Thus $W$ is a small subset of $H$.
By the inductive hypothesis,
 it is sufficient to show that the limit
$$
 \Psi(v)
 =
 \lim_{\epsilon_{a}\to 0^{+}}\psi(v+\epsilon_{a} b_{a})
$$
 exists and is locally constant for
 $v\in H\setminus W$.

Let $v\in H\setminus W$.
Choose a compact, connected, archimedean neighbourhood $U$ of $v$
 in $H$ such that $U\cap W=\emptyset$.
We shall prove that the limit $\Psi$ exists on $U$ and
 is constant there.
Let $w\in U$.
For any $i$ we either have $Z^{\prime}_{i}\subset H$
 or $w\notin Z_{i}$.
In either case there is a $\delta(w,i)>0$ sufficiently small so that
 we have $w+\epsilon_{a}b_{a}\notin Z_{i}^{\prime}$
 whenever $0<\epsilon_{a}<\delta(w,i)$.
The $\delta(w,i)$ may be chosen to be continuous in $w$.
As $U$ is compact the $\delta(w,i)$ are bounded
 below by some positive $\delta$.
This means that the subset $\tilde U=U\times (0,\delta)b_{a}$ of $\R^{a}$
 does not intersect $Z$.
The function $\psi$ is therefore defined on $\tilde U$ and since
 $\tilde U$ is connected, $\psi$ is equal to a constant $c$ on $\tilde U$.
It follows that $\Psi(w)=c$ for all $w\in U$.
\sq
\medskip

Note that if the order of the limits is changed in the above Lemma,
 then the value of the limit may change.

\subsection{Deformations of cells.}

Given a $d$-chain $\scT\in C_d(V_\infty)$
 we will describe a method for defining $\ord_v(\scT)$ for
 $v\in |\partial\scT|$.
The general idea is to replace $\scT$ by a map
 $\scT:\R^b \to C_d(V_\infty)$ so that our original $\scT$ is $\scT(0)$.
We then have a function $\psi$, defined on part of $\R^b$ by
$$
 \psi(\epsilon)
 =
 \ord_x\scT(\epsilon).
$$
If the function $\scT$ is sufficiently nice then
 we may use Lemma \ref{zariskilimit} to define
$$
 \ord_x(\scT)
 =
 \lim_{\epsilon_{1}\to 0^{+}}
 \lim_{\epsilon_{2}\to 0^{+}}
 \ldots
 \lim_{\epsilon_{b}\to 0^{+}}
 \psi(\epsilon_{1},\ldots,\epsilon_{b}).
$$
In this section,
 we investigate the conditions, which we must impose on $\scT$,
 in order for Lemma \ref{zariskilimit} to be applicable.

\subsubsection{Deformable $d-1$-cells.}

Let $\scT:\R^{b}\times \R^{d-1}\to V_{\infty}$
 be an algebraic function satisfying the following conditions:
\begin{itemize}
    \item[(A)]
    For every $\xx\in\R^{d-1}$,
     the map $\R^{b}\to V_{\infty}$ defined
     by $\epsilon \mapsto \scT(\epsilon,\xx)$
     is affine.
    We shall write $\rk(\xx)$ for the rank of the
     linear part of this map.
    \item[(B)]
    For any $i=1,\ldots,d$,
     the set of $\xx\in\R^{d-1}$ such that
     $\rk(\xx)= i$ has dimension $\le i-1$.
\end{itemize}
We shall call such a $\scT$ a \emph{deformable $d-1$-cell}.
If $\scT(\epsilon,\xx)$ is constant for a certain $\xx$ then the
 value of the constant will be called a \emph{vertex} of $\scT$.
For any $\epsilon\in\R^{b}$ we shall
 consider the $d-1$-cell $\scT(\epsilon):I^{d-1}\to V_{\infty}$
 defined by:
$$
 \Big(\scT(\epsilon)\Big)(\xx)
 :=
 \scT(\epsilon,\xx).
$$

\begin{lemma}
    Let $\scT$ be a deformable $d-1$-cell.
    Then, for any $v\in V_{\infty}$ which is not a vertex of $\scT$,
     the set $\{ \epsilon\in\R^{b} : v\in |\scT(\epsilon)|\}$
     is a small subset of $\R^{b}$.
\end{lemma}

\proof
The above set is contained in the
 projection to $\R^{b}$ of the following set:
$$
 Z^{\prime}
 =
 \{ (\epsilon,\xx)\in \R^{b}\times \R^{d-1}:
 \scT(\epsilon,\xx)=v\}.
$$
As $Z^{\prime}$ is algebraic, it remains only to show that
 $Z^{\prime}$ has codimension at least $d$ in $\R^{b}\times \R^{d-1}$.
To show this we shall write $Z^{\prime}$ as a finite union of subsets
 and bound the dimension of each of the subsets.
Define
$$
 Z^{\prime}_{i}
 =
 \{(\epsilon,\xx)\in \R^{b}\times \R^{d-1}:
 \scT(\epsilon,\xx)=v, \rk(\xx)=i\}
 \quad
 (i=0,1,2,\ldots,d).
$$
The set $Z^{\prime}$ is the union of the subsets $Z^{\prime}_{i}$.

As we are assuming that $v$ is not a vertex of $\scT$,
 it follows that $Z^{\prime}_{0}$ is empty.
For $i>0$,
 our assumption (B) on $\scT$ implies that the projection
 of $Z^{\prime}_i$ in $\R^{d-1}$ has dimension $\le i-1$.
However each fibre of this projection is an affine subspace
 with codimension $i$ in $\R^{b}$.
Therefore the dimension of $Z^{\prime}_{i}$ is $\le b-1$.
This proves the lemma.
\sq
\medskip

\subsubsection{Deformable $d$-cells}

Now suppose we have a continuous map
 $\scT:\R^b\times I^d \to V_\infty$.
Thus, for each $\epsilon\in \R^b$ we have a $d$-cell
 $\scT(\epsilon)$.
We shall call $\scT$ a \emph{deformable $d$-cell}
 if the faces of $\scT$ are deformable $d-1$-cells.
By a \emph{vertex of $\scT$},
 we shall mean a vertex of a face of $\scT$.

\begin{lemma}
    \label{pipedlim}
    Let $\scT$ be a deformable $d$-cell.
    Then for any $v\in V_{\infty}$ which is not a vertex
     of $\scT$, the following limit exists:
    $$
     \lim_{\epsilon_{1}\to 0^{+}}
     \ldots
     \lim_{\epsilon_{b}\to 0^{+}}
     \ord_{v}(\scT(\epsilon)).
    $$
    If $v\notin |\partial \scT(0)|$ then the limit
     is equal to $\ord_{v}(\scT(0))$.
\end{lemma}

\proof
We first prove the existence of the limit.
Consider the set
$$
 Z
 =
 \{\epsilon\in\R^{b}:v\in|\partial\scT(\epsilon)|\}.
$$
Since $\scT(\epsilon)$ tends uniformly to $\scT(0)$,
 the set $Z$ is archimedeanly closed in $\R^{b}$.
On $\R^{b}\setminus Z$ we have a function
$$
 \psi(\epsilon)= \ord_{v}(\scT(\epsilon)).
$$
To prove the existence of the limit,
 it is sufficient, by Lemma \ref{zariskilimit},
 to show that $Z$ is small and $\psi$ is locally constant
 on $\R^{b}\setminus Z$.

It follows from the previous lemma that $Z$ is small.
We shall show that $\psi$ is locally constant on
 $\R^{b}\setminus Z$.
Choose $\epsilon\notin Z$.
As $Z$ is closed
 there is a path connected neighbourhood $U$ of $\epsilon$
 in $\R^{b}$ which does not intersect $Z$.
We shall show that $\psi$ is constant on $U$.
Let $\epsilon^{\prime}\in U$
 and choose a path $p$ from $\epsilon$ to $\epsilon^{\prime}$ in $U$.
Now consider the $(d+1)$-chain
$$
 \scV(t,\xx)
 =
 \scT(p(t),\xx).
$$
The boundary of $\scV$ is
$$
 \partial \scV
 =
 \scT(\epsilon^{\prime})-\scT(\epsilon)
 -
 \sum_{\scU} \scV_{\scU},
$$
 where $\scV_{\scU}(t,\xx)=\scU(p(t),\xx)$
 and $\scU$ runs over the faces of $\scT$.
As $p(t)\notin Z$ we have $\scU(p(t),\xx)\ne v$.
Therefore $v\notin |\scV_{\scU}|$ so in
 $H_{d}(V_{\infty},V_{\infty}\setminus\{v\})$ we have:
$$
 \scT(\epsilon)=\scT(\epsilon^{\prime}).
$$
This implies $\psi(\epsilon)=\psi(\epsilon^{\prime})$,
 so $\psi$ is locally constant
 and we have proved the existence of the limit.

Now suppose that $v\notin |\partial \scT(0)|$.
This means that $0\notin Z$.
As $Z$ is archimedianly closed,
 there is a neighbourhood
 of $0$ on which $\psi$ is constant.
We therefore have as required:
$$
 \lim_{\epsilon_{1}\to 0^{+}}
 \ldots
 \lim_{\epsilon_{b}\to 0^{+}}
 \psi(\epsilon)
 =
 \psi(0).
$$
\sq
\medskip

\subsubsection{Deformations with respect to $\gg$.}

We shall now specialize the above result to the case which we
 require.
Recall that $\gg=\End_{\mu_m}(V_\infty)$.
Thus we have $\gg=M_s(\Q(\zeta)\otimes_\Q\R)$.
We shall write $S_\infty(\Q(\zeta))$ for the set of
 archimedean places of $\Q(\zeta)$.
There is a decomposition:
$$
 \gg
 =
 \bigoplus_{v\in S_\infty(\Q(\zeta))}
 M_s(\Q(\zeta)_v).
$$
As a $\gg$-module, $V_\infty$ decomposes as a sum of
 simple modules:
$$
 V_\infty
 =
 \bigoplus_{v\in S_\infty(\Q(\zeta))}
 \Q(\zeta)_v^s.
$$
For any subset $T\subseteq S_\infty(\Q(\zeta))$
 we shall write $V_T$ for the sum of the $\Q(\zeta)_v^s$
 for $v\in T$.
Every $\gg$-submodule of $V_\infty$ is one of the submodules $V_T$.

We consider real-algebraic functions $\scT:\gg^a\times I^{d-e}\to V_{\infty}$,
 which satisfy the following conditions:
\begin{itemize}
    \item[(C)]
    The map $\scT$ is of the form:
    $$
     \scT(\epsilon_1,\epsilon_2,\ldots,\epsilon_a,\xx)
     =
     \scT(0,\xx)+\sum_{i=1}^{a}\alpha_i\epsilon_i\phi_i(\xx),
    $$
     where the functions $\phi_i:\R^{d-1}\to V_{\infty}$ are algebraic
     and $\alpha_{1},\ldots,\alpha_{a}\in\GL_n(k_\infty)$.
    \item[(D)]
    For every non-empty subset $T\subseteq S_\infty(\Q(\zeta))$,
     the set
    $$
     \left\{
     \xx\in \R^{d-e}:
     \phi_1(\xx),\ldots,\phi_a(\xx)\in V_{T}
     \right\},
    $$
     has dimension $\le \max\{\dim_{\R}(V_{T})-e,0\}$.
\end{itemize}
Such a function $\scT$ will be called a $\gg$-deformable
 $d-e$-cell.

\begin{lemma}
    Let $\scT$ be a $\gg$-deformable $d-1$-cell.
    Then $\scT$ is a deformable $d-1$-cell when
     regarded as a function $\R^{ar}\times I^{d-1}\to V_\infty$.
    The vertices of $\scT$ are the points $\scT(\xx)$, where $\xx$
     is a solution to $\phi_1(\xx)=\ldots=\phi_i(\xx)=0$.
\end{lemma}

\proof
The statement on the vertices is clear and 
 condition (A) follows immediately from condition (C)
It remains to verify condition (B).
Let $Z_{i}=\{\xx\in\R^{d-1} : \rk(\xx)=i\}$.
We must show that for $i=1,\ldots,d$,
 the set $Z_{i}$ has dimension $\le i-1$.

For any $\xx\in \R^{d-1}$,
 the image of $\scT(\xx)$
 is a translation of a $\gg$-submodule of $V_{\infty}$.
The $\gg$-submodules of $V_{\infty}$ are of the form $V_{T}$
 for subsets $T$ of $S_{\infty}(\Q(\zeta))$.
If $V_{T}$ is the submodule corresponding to $\xx$,
 then we clearly have $\rk(\xx)=\dim(V_{T})$.

Given $T\subseteq S_{\infty}(\Q(\zeta))$,
 let $Z_{T}$ denote the set of $\xx$,
 for which the corresponding submodule is $V_{T}$.
With this notation we have:
$$
 Z_{i}
 =
 \bigcup_{T\ :\ \dim V_{T}=i} Z_{T}.
$$
As this is a finite union,
 it is sufficient to show that for any non-empty $T$,
 the set $Z_{T}$ has dimension $\le \dim(V_{T})-1$.

For a particular $\xx$,
 the corresponding submodule $V_{T}$ is the
 $\gg$-span of the vectors $\phi_{1}(\xx),\ldots,\phi_{a}(\xx)$.
Hence
$$
 Z_{T}
 \subseteq
 \{\xx\in\R^{d-1} :\phi_i(\xx)\in V_{T}\}.
$$
The result now follows from condition (D).
\sq
\medskip

We have fixed an ordered basis
 $\{b_{1},\ldots,b_{r}\}$ for $\gg$
 as a vector space over $\R$.
Let $\epsilon=\epsilon_{1}b_{1}+\ldots+\epsilon_{r}b_{r}\in \gg$.
Recall the abbreviation:
$$
 ``\blim_{\epsilon\to 0^{+}}"
 :=
 \lim_{\epsilon_{1}\to 0^{+}}
 \lim_{\epsilon_{2}\to 0^{+}}
 \ldots
 \lim_{\epsilon_{r}\to 0^{+}}.
$$
Consider a real-algebraic map
 $\scT:\gg^a\times I^{d}\to V_{\infty}$.
If the faces of $\scT$
 satisfy (C) and (D) above then $\scT$ is a
 deformable $d$-cell.
Hence by Lemma \ref{pipedlim}, the limit
$$
 \blim_{\epsilon_{1}\to 0^{+}}
 \blim_{\epsilon_{2}\to 0^{+}}
 \ldots
 \blim_{\epsilon_{a}\to 0^{+}}
 \ord_0\scT(\epsilon_{1},\ldots,\epsilon_{r})
$$
 exists for all $v$ apart from the vertices of $\scT$.

Condition (D) above is rather technical.
To be able to verify it in practice
 we shall use the following lemmata.

\begin{lemma}
\label{technical}
Let $W$ be a $d-1$-dimensional $\Q$-subspace of $V$
 and let $W_\infty$ be the closure of $W$ in $V_\infty$.
Then for any non-empty set $T$ of archimedean places of $\Q(\zeta)$,
 we have
$$
 \dim_\R(W_\infty\cap V_{T})
 =
 \dim_\R(V_{T})-1.
$$
\end{lemma}

\proof
As $W_{\infty}$ is a hyperplane in $V_{\infty}$,
 it is sufficient to show that $V_{T}$ is not contained
 in $W_{\infty}$.
It is sufficient to prove this in the case
 that $T$ consists of a single place $v$.

We have a non-degenerate $\Q$-bilinear form on $V$ given by:
$$
 <v,w>
 =
 \sum_{i=1}^s\tr_{\Q(\zeta)/\Q}(v_i w_i);
$$
 here we are identifying $V$ with $\Q(\zeta)^{s}$.
Extending the form to $V_\infty$, the various subspaces
 $V_v$ are orthogonal.
The subspace $W_\infty$ is the orthogonal complement of
 some $w\in V\setminus\{0\}$.
As the coordinate of $w$ in $V_v$ is non-zero,
 it follows that $w$ is not orthogonal to $V_v$.
Therefore $V_v$ is not a subspace of $W_\infty$ and the result follows.
\sq
\medskip

\begin{lemma}
\label{verytechnical}
Let $\scP$ be a $d-2$-dimensional cell in $\gX$
 and let $W_\scP$ be the $\R$-span of the vectors $a_{\scP,i}$.
Then for any non-empty set $T$ of archimedean places of $\Q(\zeta)$ we have
$$
 \dim_\R(W_\scP\cap V_T)
 =
 \max\{\dim_\R(V_T)-2,0\}.
$$
\end{lemma}

More general statements than the above seem to be false.
\medskip

\proof
As with the previous lemma,
 it is sufficient to prove this in the case $T=\{v\}$.
If $m=2$ then $V_v=V_\infty$ and there is nothing to prove.
We therefore assume $m>2$,
 so $\Q(\zeta_m)$ is totally complex.

The subspace $W_\scP$ is the orthogonal complement of $\{v,w\}$
 some $v,w\in V\setminus\{0\}$.
If we show that the coordinates of $v$ and $w$ in $V_v$ are linearly
 independent over $\R$, then the result follows.

Our strategy for finding the vectors $v$ and $w$ is as follows.
The cell $\scP$ is a $d-2$-dimensional face of some
 $d$-dimensional cell $\scQ$ in $\gX$.
There are two $d-1$-dimensional
 faces of $\scQ$ containing $\scP$,
 each of which is obtained by removing one of the basis elements
 $\{a_{\scQ,i}\}$.
For each $i=1,\ldots,d$
 we shall find a non-zero vector $v(i)$,
 which is orthogonal to the vectors $\{a_{\scQ,j}:j\ne i\}$.
The vectors $v,w$ may be taken to be
 $v(i),v(j)$,
 where $a_{\scQ,i}$ and $a_{\scQ,j}$ are the removed basis vectors.

We recall the construction of the basis $\{a_{\scQ,i}\}$.
We begin with a basis $\{d_1,\ldots,d_s\}$
 for the lattice $L$ over $\Z[\zeta]$.
For each $i=1,\ldots,s$,
 we choose a set of representatives
 $\zeta^{a(i,1)},\ldots,\zeta^{a(i,m/p)}$
 for $\mu_p$-cosets in $\mu_m$.
Then the basis $\{a_{\scQ,i}\}$ is as follows:
$$
 \left\{\textstyle{\frac{\rho^k}{1-\rho}}\zeta^{a(i,j)}d_i
 : i=1,\ldots,s,\; j=1,\ldots,m/p,\; k=1,\ldots,p-1\right\}.
$$
To ease notation we shall use the index set:
$$
 \cI
 =
 \{(i,j,k):i=1,\ldots,s,\; j=1,\ldots,m/p,\; k=1,\ldots,p-1\}.
$$
For $(i,j,k)\in\cI$ we define
 $a_{i,j,k}=\rho^k\zeta^{a(i,j)}a_i$.
Then the basis $\{a_{\scQ,i}\}$ is simply
 $\{\frac{1}{1-\rho}a_{\bi}:\bi\in\cI\}$.

We shall use the following Hermitean form on $V_\infty$:
$$
 \left\langle\sum v_i d_i,\sum w_i d_i\right\rangle
 =
 \frac{p}{m}
 \sum_{i=1}^s
 \tr^{\Q(\zeta)_\infty}_\R(v_i \overline{w_i}),
$$
 where $\overline{w_i}$ denotes the complex conjugate
 of $w_i$ in $\Q(\zeta)_\infty$.
We shall write $A$ for the $\cI\times \cI$ matrix
 $\langle a_\bi,a_\bj\rangle$.
One can show
 that the entries of $A$ are as follows:
$$
 \langle a_{i,j,k},a_{i',j',k'}\rangle
 =
 \left\{
 \begin{array}{ll}
 0 & (i,j)\ne(i',j'), \\
 -1 & (i,j)=(i',j'),\; k\ne k', \\
 p-1 & (i,j,k)=(i',j',k').
 \end{array}
 \right.
$$

Now consider a vector $v=(1-\overline\rho)\sum v_{\bi} a_\bi$
 and let $[v]$ be the column vector of coefficients $v_\bi$.
The vector $v$ is orthogonal to $\frac{1}{1-\rho}a_\bi$ if and only
 if the $\bi$-th row of $A[v]$ is zero.
Fix an $\bi$ and suppose $v$ is orthogonal to
 $\frac{1}{1-\rho}a_\bj$ for all $\bj\ne\bi$.
This means that all but the $\bi$-th row
 of $A[v]$ is zero, so $[v]$ is a multiple of the
 $\bi$-th column of $A^{-1}$.
We shall write $v(\bi)$ for the element of $V_\infty$,
 for which $[v(\bi)]$ is the $\bi$-th column of $A^{-1}$.
To prove the theorem we need to show show that for
 $\bi\ne \bj$ the coordinates of $v(\bi)$ and $v(\bj)$
 in $V_v$ are linearly independent over $\R$.

By finding $A^{-1}$ one obtains:
$$
 v(i,j,k)
 =
 (1-\overline\rho)(\rho^k-1)\zeta^{a(i,j)}d_i.
$$
Let the place $v$ correspond to the
 embedding $\iota:\Q(\zeta)\hookrightarrow\C$.
We must show that for $(i,j,k)\ne(i',j',k')$
 the vectors $\iota(v(i,j,k)),\iota(v(i',j',k'))\in\C^s$
 are linearly independent over $\R$.
If $i\ne i^\prime$ then this is clearly the case as they
 are independent over $\C$.
We therefore assume $i=i'$.
We are reduced to showing that the complex number
$$
 z
 =
 \iota\left(
 \frac{(\rho^k-1)\zeta^{a(i,j)}}{(\rho^{k'}-1)\zeta^{a(i,j')}}
 \right).
$$
 is not real.
We let $z_1=\iota\left(\frac{\rho^k-1}{\rho^{k'}-1}\right)$
 and $z_2=\iota(\zeta)^{a(i,j)-a(i,j')}$.
The argument of $z_1$ is of the form
$$
 \pi (k-k')\frac{r}{p},
 \quad
 r\in\{1,2,\ldots,p-1\}.
$$
Therefore $z_1^p\in\R$.
If $j\ne j'$ then $z_2^p\notin\R$ so we are done.
Finally, if $j=j'$ then $z_2=1$ but $k\ne k'$,
 so $z_1\notin\R$ and again we are done.
\sq
\medskip

Up until now,
 we have examined deformable cells in $V_\infty$.
However,
 we need to define the order of a cell in $X_\infty$
 at points on its boundary.
We call a map $\scT:\R^b\times I^{d}\to X_\infty$ a deformable cell
 in $X_{\infty}$,
 if one, or equivalently all,
 of its lifts to $V_\infty$ are deformable cells.
We define the vertices of such a $\scT$ to be the projections
 in $X_\infty$ of the vertices of a lift of $\scT$.

\begin{lemma}
\label{Xlimit}
	Let $\scT:\R^b\times I^d\to X_\infty$ be a deformable $d$-cell
         in $X_\infty$.
	If $x\in X_\infty$ is not a vertex of $\scT$
	 then the following limit exists:
	$$
	 \lim_{\epsilon_{1}\to 0^{+}}
         \ldots
         \lim_{\epsilon_{b}\to 0^{+}}
         \ord_{x}(\scT(\epsilon)).
	$$
\end{lemma}

\proof
This follows from the previous results using the relation
$$
 \ord_x(\scT(\epsilon))
 =
 \sum_{y\to x}\ord_y(\tilde\scT(\epsilon)),
$$
 where $\tilde\scT$ is a lift of $\scT$ and the sum is over
 the preimages of $x$ in $V_\infty$.
As $|\tilde\scT|$ is compact, the sum is in fact finite
 and hence commutes with the limits.
\sq
\medskip

It will be more convenient to speak of
 \emph{piecewise deformable cells}.
We shall call a map $\R^{b}\times I^d\to X_{\infty}$
 a piecewise deformable cell,
 if there is a subdivision of $I^d$,
 such that the restriction of $\scT$
 to any of the pieces in the subdivision is deformable.
The reason for this is that our functions will be piecewise
 algebraic rather that algebraic.
By a \emph{piecewise deformable chain} we shall simply
 mean a formal sum of piecewise deformable cells.
Thus a piecewise deformable chain will be a map $\R^b\to C_\bullet$.
If $x$ is not a vertex of a piecewise deformable $d$-chain $\scT$
 then we define
$$
 \ord_x(\scT)
 :=
 \lim_{\epsilon_{1}\to 0^{+}}
 \ldots
 \lim_{\epsilon_{b}\to 0^{+}}
 \ord_{x}(\scT(\epsilon)).
$$
It follows from the above results that this limit exists.
If $\scT$ is deformable then we shall call $\scT$ a \emph{deformation}
 of $\scT(0)$.

\subsection{Deforming paths.}

\paragraph{The function $\bar f$.}
We now apply Lemma  \ref{pipedlim}
 to the function $f$.
Define a path $\path(\epsilon)$ for $\epsilon\in\gg$ by
$$
 \path(\epsilon)
 =
 \left[0, \textstyle{\frac{1}{2}}+\epsilon\right]
 +
 \left[\textstyle{\frac{1}{2}}+\epsilon,1\right].
$$

\begin{proposition}
For any $d-1$ or $d-2$-dimensional cell $\scP$ in $\gX$
The map $\epsilon \mapsto \path(\epsilon)\bowtie\scP$
 is piecewise $\gg$-deformable.
Its vertices are in $\frac{1}{1-\rho}L$.
\end{proposition}

\proof
Let $\scP$ be any $d-1$-cell in $\gX$.
We first cut $\path(\epsilon)\bowtie\scP$ into its $2^{d-1}$
 algebraic pieces and then prove that each piece is deformable.
Thus for any subset $A\subseteq\{1,2,\ldots,d-1\}$
 we define
$$
 \scA_A(\epsilon)
 =
 [v_\scP]
 \Diamond
 \bigDiamond_{i\in A}
 [0, (\textstyle{\frac{1}{2}}+\epsilon)a_{\scP,i}]
 \Diamond
 \bigDiamond_{i\notin A}
 [(\textstyle{\frac{1}{2}}+\epsilon) a_{\scP,i},a_{\scP,i}].
$$
We shall show that each $\scA_A$ is deformable
 by verifying conditions (C) and (D) above.
We have
$$
 \scA(\epsilon,\xx)
 =
 v_\scP
 +
 \sum_{i\in A}
 (\textstyle{\frac{1}{2}}+\epsilon) x_i a_{\scP,i}
 +
 \displaystyle\sum_{i\notin A}
 ((\textstyle{\frac{1}{2}}+\epsilon) +(\textstyle{\frac{1}{2}}-\epsilon)x_i)
 a_{\scP,i}.
$$
This implies
$$
 \scA(\epsilon,\xx)
 =
 \scA(0,\xx)
 +\epsilon \left(\sum_{i\in A} x_i a_{\scP,i}
 +\sum_{i\notin A}  (1-x_i) a_{\scP,i}\right).
$$
Therefore $\scA$ verifies condition (C) with
$$
 \phi(\xx)
 =
 \sum_{i\in A} x_i a_{\scP,i}+\sum_{i\notin A}  (1-x_i) a_{\scP,i}.
$$
To verify (D) we let $W_\scP$ be the $\R$-span of the
 vectors $a_{\scP,i}$.
Thus $\phi$ maps $\R^{d-1}$ bijectively to $W_{\scP}$.
We must show that for any non-empty set of archimedean places $T$
 we have $\dim_\R(W_\scP\cap V_T)\le\dim_\R(V_T)-1$.
This follows from Lemma \ref{technical}.
It follows from the formula for $\phi$ that the only vertex of
 $\scA$ is $v_\scP+\sum_{i\notin A} a_{\scP,i}$.
By Lemma \ref{constructX} we know that this is in $\frac{1}{1-\rho}L$.
The case of a $d-2$-cell in $\gX$ is similar except
 that one must use Lemma \ref{verytechnical}
 instead of Lemma \ref{technical}.
\sq
\medskip

We may now define
$$
 \bar f(x)
 :=
 \blim_{\epsilon\to 0^{+}} \ord_x(\path(\epsilon)\bowtie\scF).
$$
This limit exists for all $x$ not in $\frac{1}{1-\rho}L$.

\begin{proposition}
    \label{flimit}
    If $f(x)$ is defined then so is $\bar f(x)$ and they are equal.
    However $\bar f(x)$ is defined for all $x\notin\Vert(\scP)$.
    Furthermore if the $\mu_{m}$-orbit of $x$ does not intersect
     $\Vert(\scP)$ then we have
    $$
     \sum_{\zeta\in\mu_{m}}\bar f(\zeta x)
     =
     1.
    $$
    In particular the restriction of $f$ to $X\setminus \{0\}$
     is a fundamental function.
\end{proposition}

\proof
The first two assertions follow
 from Lemma \ref{Xlimit}.
To prove the formula,
 we use the fact that finite sums
 commute with limits as follows:
$$
 \sum_{\zeta\in\mu_{m}}
 \bar f(\zeta x)
 =
 \sum_{\zeta\in\mu_{m}}
 \blim_{\epsilon\to 0^{+}}
 f^{(\path(\epsilon))}(\zeta x)
 =
 \blim_{\epsilon\to 0^{+}}
 \sum_{\zeta\in\mu_{m}}
 f^{(\path(\epsilon))}(\zeta x)
 =
 \blim_{\epsilon\to 0^{+}} 1
 =
 1.
$$
\sq

\paragraph{The function $\bar f^{\alpha}$.}
Recall that for any $\alpha\in \GL_{n}(k_{\infty})$
 we have a path
 $\path^{\alpha}$ from $0$ to $1$ in $\gg$,
 defined by
$$
 \path^{\alpha}
 =
 [0,\alpha] + [\alpha,1].
$$
More precisely, let
$$
 \path^{\alpha}(x)
 =
 \left\{
 \begin{array}{lll}
     2x \alpha & \hbox{for} & x\le \frac{1}{2},\medskip \\
     \alpha+(2x-1)(1-\alpha) & & x \ge \frac{1}{2}.
 \end{array}
 \right.
$$
By a vertex of $\path^{\alpha}$ we shall mean
 one of the following points of $I$:
$$
 \Vert(\path^{\alpha})
 =
 \left\{
 0,\frac{1}{2},\frac{m^{2}+1}{2m^{2}},\frac{m^{2}+2}{2m^{2}}, \ldots,1
 \right\}.
$$
By a vertex of $\path^{\alpha}\bowtie\scF$
 we shall mean a point $v\in V_{\infty}$ of the form
$$
 v
 =
 (\path^{\alpha}\bowtie\scP)(\xx),
 \quad
 \xx
 \in
 \Vert(\path^{\alpha})^{d-1},
$$
 where $\scP$ is a $d-1$-cell in $\gX$.

We shall write $\Vert(\path^{\alpha}\bowtie\scF)$
 for the set of all vertices  of $\path^{\alpha}\bowtie\scF$.
The path $\path^{\alpha}$ gives rise to a fundamental function
 $f^{\alpha}$ away from the boundary of $\path^{\alpha}\bowtie\scP$.
We shall extend the definition of $f^{\alpha}$ to all
 points of $X_{\infty}$ apart from the vertices of $\path^{\alpha}\bowtie\scF$.

Given $\epsilon,\nu\in\gg$ we define a new path
 $\path^{\alpha}(\epsilon,\nu)$ by
\begin{eqnarray*}
 \path^{\alpha}(\epsilon,\nu,x)
 &=&
 \left\{
 \begin{array}{lll}
     \alpha\path(\epsilon,2x) & \hbox{for}& x\le \frac{1}{2},
     \smallskip\\
     \alpha + (2x-1)(1-\alpha) + \phi(2m^{2}x)\nu
                              & &x\ge \frac{1}{2}.
 \end{array}
 \right.
\end{eqnarray*}
Here $\phi:\R/\Z\to \R$ is the $\Z$-periodic function
 defined on the interval $I$ by
$$
 \phi(x)
 =
 \left\{
 \begin{array}{ll}
     x & x\le \frac{1}{2},
     \smallskip\\
     1-x & x\ge \frac{1}{2}.
 \end{array}
 \right.
$$
The path $\path^{\alpha}(\epsilon,\nu)$
 reduces to $\path^{\alpha}$
 when $\epsilon$ and $\nu$ are both zero.
We extend our definition of $f^{\alpha}$ as follows:
$$
 \bar f^{\alpha}(x)
 =
 \blim_{\epsilon\to 0^{+}}
 \blim_{\nu\to 0^{+}}
 \ord_x(\path^{\alpha}(\epsilon,\nu)\bowtie\scF).
$$

\begin{proposition}
    \label{falimit}
    Let $\scP$ be a $d-1$- or $d-2$-cell in $\gX$.
    The map $(\epsilon,\nu)\mapsto
     \path^{\alpha}(\epsilon,\nu)\bowtie\scP$
     is piecewise $\gg$-deformable.
    It is a deformation of $\path^{\alpha}\bowtie\scP$.
    Its vertices are in $\Vert(\path^{\alpha}\bowtie\scP)$.
    Consequently the limit $\bar f^\alpha(x)$ exists for all
     $x\notin \Vert(\path^{\alpha}\bowtie\scF)$.
    The function $\bar{f}^{\alpha}$ is fundamental on all $\mu_{m}$-orbits
     which do not intersect $\Vert(\path^{\alpha}\bowtie\scF)$.
\end{proposition}

\proof
Let $\scP$ be any $d-1$-cell in $\gX$.
As before we must show that the map
$$
 \scT(\epsilon,\nu,\xx)
 =
 (\path^{\alpha}(\epsilon,\nu)\bowtie\scP)(\xx)
$$
 is piecewise deformable.

We cut $\scT$ into its algebraic pieces.
These pieces are translations
 by elements of $\Vert(\path^{\alpha}\bowtie\scP)$
 of the pieces
\begin{eqnarray*}
 \scA
 &=&
 \left(
 \alpha({\textstyle\frac{1}{2}}+\epsilon)\cdot \bigDiamond_{i\in A} [0, a_{\scP,i}]
 \right)\\
 &&
 \Diamond
 \left(
 \alpha({\textstyle\frac{1}{2}}
 -\epsilon) \cdot \bigDiamond_{i\in B} [0, a_{\scP,i}]
 \right)\\
 &&
 \Diamond
 \left(
 ({\textstyle\frac{1-\alpha}{2m^{2}}+\frac{\nu}{2}})
 \cdot
 \bigDiamond_{i\in C} [0, a_{\scP,i}]
 \right)\\
 &&
 \Diamond
 \left(
 ({\textstyle\frac{1-\alpha}{2m^{2}}-\frac{\nu}{2}})
 \cdot
 \bigDiamond_{i\in D} [0, a_{\scP,i}]
 \right),
\end{eqnarray*}
 where the sets $A,B,C,D$ form a partition of $\{1,2,\ldots,d-1\}$.
The piece $\scA$ satisfies condition (C) above
 with
$$
 \phi_1(\xx)
 =
 \left(\sum_{i\in A}-\sum_{i\in B}\right) x_i a_{\scP,i},
 \quad
 \phi_2(\xx)
 =
 \left(\sum_{i\in C}-\sum_{i\in D}\right) x_i a_{\scP,i}.
$$
To prove condition (D) we apply
 lemma \ref{technical} to the subspace $W_{\scP}$
 spanned by $\{a_{\scP,i}\}$.
The case of a $d-2$-cell $\scP$ is similar
 but one must use Lemma \ref{verytechnical}
 instead of Lemma \ref{technical}.

As the $a_{\scP,i}$ are linearly independent
 it follows that the only vertex of $\scA$ is $0$,
 so the vertices of $\scT$ are the translations,
 which are in $\Vert(\path^{\alpha}\bowtie\scP)$.
\sq
\medskip

\paragraph{The function $\bar f^{\alpha\beta,\alpha}$.}
Finally for $\alpha,\beta\in \GL_{n}(k_{\infty})$
 we define a path $\path^{\alpha\beta,\alpha}$ by
$$
 \path^{\alpha\beta,\alpha}
 =
 [0,\alpha\beta]+[\alpha\beta,\alpha]+[\alpha,1].
$$
More precisely let
$$
 \path^{\alpha\beta,\alpha}(x)
 =
 \left\{
 \begin{array}{ll}
     4x \alpha\beta  & x\le \frac{1}{4}, \\
     \alpha\beta + (4x-1)(\alpha-\alpha\beta))
                     & \frac{1}{4}\le x \le \frac{1}{2},\\
     \alpha + (2x-1)(1-\alpha) & x \ge \frac{1}{2}.
 \end{array}
 \right.
$$
We define
$$
 \Vert(\path(\alpha\beta,\alpha))
 =
 \left\{
 0,\frac{1}{4},\frac{m^{2}+1}{4m^{2}},
 \ldots,
 \frac{1}{2},\frac{m^{2}+1}{2m^{2}},\ldots,1
 \right\},
$$
$$
 \Vert(\scP^{\alpha\beta,\alpha})
 =
 \{ \scP_{i}^{\alpha\beta,\alpha}(\xx):
 \xx \in\Vert(\path(\alpha\beta,\alpha))^{d-1}\}.
$$
We shall extend the definition of $f^{\alpha\beta,\alpha}$ to
 $X_{\infty}\setminus\Vert(\scP^{\alpha\beta,\alpha})$.
To do this we define a deformation of
 $\path^{\alpha\beta,\alpha}$ as follows:
\begin{eqnarray*}
 \path^{\alpha\beta,\alpha}(\epsilon,\nu,\xi,x)
 &=&
 \left\{
 \begin{array}{ll}
     \alpha\path^{\beta}(\epsilon,\nu,2x)
     &
     x\le \frac{1}{2},
     \smallskip\\
     \path^{\alpha}(0,\xi,x)
     &
     x\ge \frac{1}{2}.
 \end{array}
 \right.
\end{eqnarray*}
Again, for any $x\in X_{\infty}$,
 which is not a vertex of $\scP^{\alpha\beta,\alpha}$,
 we may define
$$
 \bar f^{\alpha\beta,\alpha}(x)
 =
 \blim_{\epsilon\to 0^{+}}
 \blim_{\nu\to 0^{+}}
 \blim_{\xi\to 0^{+}}
 \ord_x(\path^{\alpha\beta,\alpha}(\epsilon,\nu,\xi)\bowtie \scF).
$$

\begin{proposition}
    \label{fabalimit}
    Let $\scP$ be a $d-1$- or $d-2$-cell in $\gX$.
    The map
    $$
     (\epsilon,\nu,\xi)
     \mapsto
     \path^{\alpha\beta,\alpha}(\epsilon,\nu,\xi)\bowtie\scP
    $$
     is piecewise
     a $\gg$-deformation of $\path^{\alpha\beta,\alpha}\bowtie\scP$.
    Its vertices are in $\Vert(\path^{\alpha\beta,\alpha}\bowtie\scF)$.
    Consequently the limit $\bar f^{\alpha\beta,\alpha}(x)$ exists for all
     $x\notin \Vert(\path^{\alpha\beta,\alpha}\bowtie\scF)$ and is
     fundamental there.
\end{proposition}

This is proved in a similar way to Proposition \ref{falimit}.

\subsection{Deformation of homotopies.}

Again to take account of points on the boundary of $\scH\bowtie\scG$
 we must construct a deformation of $\scH\bowtie\scG$.

\begin{proposition}
    Suppose that for
     $\epsilon_{1},\ldots,\epsilon_{r},\nu_{1},\ldots,\nu_{s}\in\gg$
     we have paths $\path_{1}(\epsilon_{1},\ldots,\epsilon_{r})$
     and $\path_{2}(\nu_{1},\ldots,\nu_{s})$.
    Suppose further that for any $d-1$- or $d-2$-cell $\scP$
     in $\gX$,
     the maps $\path_1\bowtie \scP$ and $\path_2\bowtie \scP$
     are piecewise $\gg$-deformable.
    Define a homotopy
     $\scH(\epsilon_{1},\ldots,\epsilon_{r},\nu_{1},\ldots,\nu_{s})$
     from $\path_{1}(\epsilon_{1},\ldots,\epsilon_{r})$
     to $\path_{2}(\nu_{1},\ldots,\nu_{s})$ by
    $$
     \scH(t,x)=(1-t)\path_1(x)+t\path_2(x).
    $$
    Then $\partial(\scH\bowtie\scG)$ is piecewise $\gg$-deformable.
\end{proposition}

\proof
Recall that $\scG\in \gX_{d-1}$.
We must therefore show that for any $d-1$-cell $\scP$ in $\gX$,
 the $d-1$-chain $\partial(\scH\bowtie\scQ)$ is
 piecewise $\gg$-deformable.
We have
$$
 \partial(\scH\bowtie\scP)
 =
 (\path_2-\path_1)\bowtie\scP +\scH\bowtie\partial\scP.
$$
As we already know that $\path_i\bowtie\scP$ is
 piecewise $\gg$-deformable,
 it is sufficient to show that for any $d-2$-cell $\scQ$
 the cell $\scH\bowtie\scQ$ is piecewise $\gg$-deformable.
We have
$$
 (\scH\bowtie\scQ)(t,\xx)
 =
 (1-t)(\path_1\bowtie\scQ)(\xx)+t(\path_2\bowtie\scQ)(\xx).
$$
This implies
$$
 \rk_{\scH\bowtie\scQ}(t,\xx)
 \ge
 \min\{ \rk_{\path_1\bowtie\scQ}(\xx),\rk_{\path_2\bowtie\scQ}(\xx)\}.
$$
The result now follows as $\path_i\bowtie\scQ$ is
 piecewise $\gg$-deformable.
\sq
\medskip

\subsection{The limit defining $\deci$.}

Recall that we have defined $\deci(\alpha,\beta)$ in terms
 of the limit
$$
 \blim_{\epsilon_1\to 0^+}
 \blim_{\epsilon_2\to 0^+}
 \blim_{\epsilon_3\to 0^+}
 \;
 \ord_{0,V}
 \Big(
 [1+\epsilon_1,\alpha(1+\epsilon_2),\alpha\beta(1+\epsilon_3)]\cdot\scE
 \Big).
$$
To show that this limit exists we must prove:

\begin{proposition}
The map
$$
 (\epsilon_1,\epsilon_2,\epsilon_3)
 \mapsto
 [1+\epsilon_1,\alpha(1+\epsilon_2),\alpha\beta(1+\epsilon_3)]\cdot\scE
$$
 is a $\gg$-deformation of $[1,\alpha,\alpha\beta]\cdot\scE$
 with no vertices.
\end{proposition}

\proof
We have $\scE=\gs \scG$,
 where $\scG$ is an element of $\gX_{d-1}$.
We therefore fix a $d-1$-cell $\scP$
 in $\gX$ and define
$$
 \scA
 =
 [1+\epsilon_1,\alpha(1+\epsilon_2),\alpha\beta(1+\epsilon_3)]
 \cdot\gs\scP.
$$
We must show that the boundary of $\scA$ is deformable.
As $\gs$ anticommutes with $\partial$ we have:
$$
 \partial\scA
 =
 \partial([1+\epsilon_1,\alpha(1+\epsilon_2),\alpha\beta(1+\epsilon_3)])
 \cdot\gs\scP
 -
 [1+\epsilon_1,\alpha(1+\epsilon_2),\alpha\beta(1+\epsilon_3)]
 \cdot\gs\partial\scP.
$$
We must show that both the summands above are deformable.
We shall show that all the cells in the above expression
 satisfy conditions (C) and (D) above
 and have no vertices.

We begin with the first summand.
This is made up of cells of the form
$$
 [(1+\epsilon_1),\alpha(1+\epsilon_2)]\cdot\gs\scP:
 \Delta^{1}\times\Delta^{d-2}\to V_{\infty}.
$$
Condition (C) is satisfied with
$$
 \phi_i(\xx,\yy)
 =
 x_i \gs\scP(\yy),
 \quad
 (\xx,\yy)\in\Delta^{1}\times\Delta^{d-2}.
$$
As $0\notin|\gs\scP|$ and the $x_i$ are not all zero
 it follows that there are no vertices.
Let $H$ be the hyperplane in $W_\scP$ containing $|\gs\scP|$.
To verify condition (D) we must show that
 $H\cap V_T$ has dimension $\le \dim_\R(V_T)-2$.
This follows from Lemma \ref{technical} since $H$ does not contain $0$.

The second summand contains cells of the form:
$$
 [(1+\epsilon_1),\alpha(1+\epsilon_2),\alpha\beta(1+\epsilon_3)]
 \cdot
 \gs\scQ:
 \Delta^{2}\times\Delta^{d-3}\to V_{\infty},
$$
 with $\scQ$ a $d-2$-cell of $\gX$.
Again condition (C) is satisfied with
$$
 \phi_i(\xx,\yy)
 =
 x_i \gs\scQ(\yy),
 \quad
 (\xx,\yy)\in\Delta^{2}\times\Delta^{d-3}
$$
As the $x_i$ are never all zero and $0\notin|\gs\scQ|$,
 it follows that the functions $\phi_i$
 are never simultaneously all $0$.
This shows that there are no vertices.

The base set of $\gs\scQ$ lies in a hyperplane
 $H$ in $W_{\scQ}$.
To verify condition (D) we must show that for any set $T$
 of archimedean places of $\Q(\zeta)$ the dimension of
 $H\cap V_T$ is $\le \dim_\R(V_T)-3$.
As $H$ does not go through $0$ this reduces to proving that
 $\dim_\R(W_{\scP}\cap V_T)\le \dim_\R(V_T)-2$.
However this follows from Lemma \ref{verytechnical}.
\sq
\medskip

\section{The relation between the arithmetic and geometric cocycles.}

\subsection{Main Results.}

We now state our main results.
Recall that we have a homotopy $\scH^{1}_{\beta}$ from $\path(1)$ to
 $\path(\beta)$ defined by
$$
 \scH^{1}_{\beta}(x,t)
 =
 t\path(\beta)(x) + (1-t)x.
$$
From this we have constructed a homotopy
 $\scH^{\alpha}_{\alpha\beta,\alpha}$ from
 $\path(\alpha)$ to $\path(\alpha\beta,\alpha)$ as follows:
$$
 \scH^{\alpha}_{\alpha\beta,\alpha}(x,t)
 =
 \left\{
 \begin{array}{ll}
     \alpha \scH^{1}_{\beta}(2x,t)
     & x\le \frac{1}{2},\smallskip\\
     \path(\alpha)(x)
     & x\ge \frac{1}{2}.
 \end{array}
 \right.
$$
Finally we have a homotopy
 $\scH^{\alpha\beta,\alpha}_{\alpha\beta}$
 from $\path^{\prime}(\alpha\beta,\alpha)$ to $\path(\alpha\beta)$ satisfying
$$
 \scH^{\alpha\beta,\alpha}_{\alpha\beta}(x,t)
 =
 2x\alpha\beta,\quad x\le \frac{1}{2}.
$$
 and constructed by splitting the triangle
 $[\alpha\beta,\alpha,1]$ into $m^{2}$ smaller triangles.

Recall that our geometric cocycle is given by the formula:
$$
 \dec_\infty^{(\gs\scF)}(\alpha,\beta)
 =
 \zeta^{
 \textstyle{\ord_{0,V}
 \left([1,\alpha,\alpha\beta]\cdot\gs\scG
 \right)}}.
$$
By Proposition \ref{XV} we have
\begin{equation}
 \label{geometric}
 \dec_\infty^{(\gs\scF)}(\alpha,\beta)
 =
 \zeta^{
 \textstyle{\ord_{0,X}\left(
 (\scH^{1}_{\alpha}+\scH^{\alpha}_{\alpha\beta,\alpha}
 +\scH^{\alpha\beta,\alpha}_{\alpha\beta}
 -\scH^{1}_{\alpha\beta})\bowtie\scG
 \right)}}.
\end{equation}
One may check that this remains true even if the quantities
 are defined as limits in the sense of \S5.
Looked at from this point of view, it seems natural to
 divide out a certain coboundary from this cocycle.
For $\alpha\in\Upsilon_\gf$ define
$$
 \tau(\alpha)
 =
 \zeta^{
 \textstyle{\left\{
 \scH^{1}_{\alpha}\bowtie\scG
 |\alpha L
 \right\}}}.
$$
We shall prove the following.

\begin{theorem}
\label{main}
For $\alpha,\beta\in\Upsilon_\gf$ we have
$$
 \deci(\alpha,\beta)\decAS(\alpha,\beta)
 =
 \frac{\tau(\alpha)\tau(\beta)}{\tau(\alpha\beta)}.
$$
\end{theorem}

We shall show that there is a continuous
 cocycle $\dec_m$ on $\SL_n(k_m)$ and an extension of $\tau$ to
 $\SL_n(k)$ such that for $\alpha,\beta\in\SL_n(k)$ we have
\begin{equation}
    \label{globalsplit}
    \deci(\alpha,\beta)\dec_m(\alpha,\beta)\decAS(\alpha,\beta)
    =
    \frac{\tau(\alpha)\tau(\beta)}{\tau(\alpha\beta)}.
\end{equation}

The proof of Theorem \ref{main}
 will be broken down into the following three lemmata.

\begin{lemma}
    \label{A}
    Let $\alpha\in\Upsilon_{\gf}$.
    Then for any $\mu_{m}$-invariant lattice $M\subset V_{m}$
     containing $L$, we have
    $$
     \langle f\alpha^{-1}-f^{\alpha} |\alpha M-\alpha L\rangle_{X}
     =
     1.
    $$
\end{lemma}

\begin{lemma}
    \label{C}
    Let $\alpha\in G_{\gf}$ and $\beta\in\GL_n(k_\infty)$.
    Then for any $\mu_{m}$-invariant lattice $M\subset V_{m}$
     containing both $L$ and $\alpha^{-1}L$,
     the following holds in $\Z/m$:
    $$
     \left\{\scH^1_\beta\bowtie\scG|M\right\}
     =
     \left\{\scH^{\alpha}_{\alpha\beta,\alpha}\bowtie\scG|\alpha M\right\}.
    $$
\end{lemma}

\begin{lemma}
    \label{D}
    Let $\alpha,\beta\in G_{\gf}$.
    For any $\mu_{m}$-invariant lattice $M\subset V_{m}$
     containing $L$, $\alpha L$ and $\alpha\beta L$,
     the following holds in $\Z/m$:
    $$
     \left\{\scH^{\alpha\beta}_{\alpha\beta,\alpha}\bowtie\scG|M\right\}
     =
     0.
    $$
\end{lemma}

The proofs of these lemmata in \S6.4-6.6 are essentially exercises in
 the geometry of numbers.
In each case one is reduced to showing
 that the number of lattice points in a certain
 set is a multiple of $m$.
One achieves this by cutting the
 set into $m$ pieces which are translations of one another
 by elements of the lattice.

\subsection{A formula for the Kubota symbol.}

Assume in this section that $k$ is
 totally complex.
We shall describe the Kubota symbol on the
 group
$$
 \Gamma_{\gf}
 =
 \{\alpha\in \SL_{n}(k): \alpha L=L,\; \alpha\equiv I_{n} \bmod \gf\}
 =
 \Upsilon_{\gf}\cap \SL_{n}(\go).
$$

\begin{corollary}
    \label{kubota}
    For $\alpha\in\Gamma_{\gf}$ the Kubota symbol
     $\kappa_{m}(\alpha)$ is given by the formula:
    $$
     \kappa_{m}(\alpha)
     =
     \frac{\zeta^{\textstyle{\ord_{0,X}\left(\scH^1_\alpha\bowtie\scG\right)}}}
     {h(w(\alpha))}.
    $$
\end{corollary}

\proof
We shall regard $\SL_{n}(k)$ as a dense subgroup of $\SL_{n}(\Af)$,
 where $\Af$ denotes the ring of finite ad\`eles of $k$.
Let $U$ be the closure of $\Gamma_{\gf}$ in $\SL_{n}(\Af)$.
This is a compact open subgroup of $\SL_{n}(\Af)$.
Consider the extension
$$
 1 \to
 \mu_{m} \to
 \SLt_{n}(\Af) \to
 \SL_{n}(\Af) \to
 1,
$$
 corresponding to the cocycle $\decAS\decm$.
The cocycle $\decAS\decm$ is 1 on $U\times U$.
Therefore on $U$ the map $\alpha\mapsto(\alpha,1)$
 is a splitting of the extension.
On $\SL_{n}(k)$ we have by (\ref{globalsplit})
 and Theorem \ref{complexsplit}:
$$
 \decAS\decm(\alpha,\beta) \partial hw(\alpha,\beta)
 =
 \partial\tau(\alpha,\beta).
$$
This implies that on $\SL_{n}(k)$ the map
 $\alpha\mapsto (\alpha, \tau(\alpha)/hw(\alpha))$
 splits the extension.
As the Kubota symbol is the ratio
 of these two splittings,
 the result follows.
\sq
\medskip

\begin{remark}
    Some authors speak of the ``Kubota symbol on $\GL_n$''.
    By this they are in effect choosing an embedding of $\GL_n$
     in $\SL_{n+r}$ and pulling back the Kubota symbol.
    As the above formula is valid for $n$ arbitrarily large
     there is no need to have a separate formula for $\GL_n$.
    It is worth mentioning that the formula of the above corollary
     gives a homomorphism
    $$
     \GL_{n}(\go,\gf) \to \Mu,\quad
     \alpha
     \mapsto
     \frac{\zeta^{\textstyle{\ord_{0,X}\left(\scH^1_\alpha\bowtie\scG\right)}}}
     {w(\alpha)},     
    $$
     which is a rather more canonical extension
     of the Kubota symbol on $\GL_{n}$.
    Here $\GL_{n}(\go,\gf)$ represents the principal congruence
     subgroup modulo $\gf$.
\end{remark}

\subsection{Proof of Theorem \ref{main}.}

We shall deduce the theorem from lemmata \ref{A}, \ref{C} and \ref{D}.

Let $\alpha,\beta\in \Upsilon_\gf$.
We begin with the definition of $\decAS$:
$$
 \decAS(\alpha,\beta)
 =
 \langle f-f\alpha | \beta L - L \rangle
 =
 \langle f\alpha^{-1}-f | \alpha\beta L - \alpha L \rangle.
$$
By Lemma \ref{A} we have:
$$
 \decAS(\alpha,\beta)
 =
 \langle f^\alpha-f | \alpha\beta L - \alpha L \rangle.
$$
By Proposition \ref{productformula} we have:
$$
 \decAS(\alpha,\beta)
 =
 \zeta^{
 \textstyle{\left\{\scH^\alpha_1\bowtie\scG | \alpha\beta L - \alpha
L\right\}}}.
$$
This implies
\begin{eqnarray*}
 \decAS(\alpha,\beta)
 &=&
 \tau(\alpha)
 \zeta^{
 \textstyle{\left\{\scH^\alpha_1\bowtie\scG | \alpha\beta L\right\}}}\\
 &=&
 \frac{\tau(\alpha)}{\tau(\alpha\beta)}
 \zeta^{
 \textstyle{
 \left\{(\scH^\alpha_1+\scH^1_{\alpha\beta})\bowtie\scG |
 \alpha\beta L\right\}}}.
\end{eqnarray*}
On the other hand we have by Lemma \ref{C}:
$$
 \tau(\beta)
 =
 \zeta^{
 \textstyle{\{\scH^1_\beta\bowtie\scG|\beta L\}}}
 =
 \zeta^{
 \textstyle{\{\scH^\alpha_{\alpha\beta,\alpha}\bowtie\scG| \alpha\beta L\}}}.
$$
Taking the last two formulae together we obtain:
$$
 \decAS(\alpha,\beta)
 =
 \frac{\tau(\alpha)\tau(\beta)}{\tau(\alpha\beta)}
 \zeta^{
 \textstyle{\left\{
 (\scH^\alpha_1
 +\scH^1_{\alpha\beta}
 -\scH^\alpha_{\alpha\beta,\alpha})
 \bowtie\scG | \alpha\beta L\right\}}}.
$$
By Lemma \ref{D} we have
$$
 \decAS(\alpha,\beta)
 =
 \frac{\tau(\alpha)\tau(\beta)}{\tau(\alpha\beta)}
 \zeta^{
 \textstyle{\left\{
 (\scH^\alpha_1
 +\scH^1_{\alpha\beta}
 +\scH^{\alpha\beta}_{\alpha\beta,\alpha}
 +\scH^{\alpha\beta,\alpha}_\alpha)
 \bowtie\scG | \alpha\beta L\right\}}}.
$$
This implies by (\ref{geometric}):
$$
 \decAS(\alpha,\beta)
 \deci(\alpha,\beta)
 =
 \frac{\tau(\alpha)\tau(\beta)}{\tau(\alpha\beta)}
 \zeta^{\textstyle{\left\{
 (\scH^\alpha_1
 +\scH^1_{\alpha\beta}
 +\scH^{\alpha\beta}_{\alpha\beta,\alpha}
 +\scH^{\alpha\beta,\alpha}_\alpha)
 \bowtie\scG | \alpha\beta L-L\right\}}}.
$$
On the other hand by Proposition \ref{productformula}
 the final term in the above is equal to:
$$
 \langle
 f^\alpha -f
 +f-f^{\alpha\beta}
 +f^{\alpha\beta}-f^{\alpha\beta,\alpha}
 +f^{\alpha\beta,\alpha}-f^\alpha|
 \alpha\beta L-L\rangle
 =
 1.
$$
\sq

\subsection{Proof of Lemma \ref{A}.}

We have by definition,
\begin{eqnarray*}
 \langle f^{\alpha}-f\alpha^{-1}| \alpha M - \alpha L \rangle_{X}
 &=&
 \prod_{\zeta\in\mu_{m}}
 \zeta^{\textstyle{\sum_{x\in T} f^{\alpha}(x)f(\zeta\alpha^{-1} x)}},
\end{eqnarray*}
 where the sums are over all $x$ in the finite subset $T$ of $X$
 given by:
$$
 T
 =
 ((\alpha M)/L) \setminus ((\alpha L)/L).
$$
The functions $f$ and $f^{\alpha}$ are defined as limits
 of the functions $f^{(\epsilon)}$ and $f^{\alpha,(\epsilon,\nu)}$.
As $T$ is finite we may choose $\epsilon,\nu$ small enough
 so that for every $x\in T$ we have
 $f(\alpha^{-1}x)=f^{(\epsilon)}(\alpha^{-1}x)$
 and
 $f^{\alpha}(x)=f^{\alpha,(\epsilon,\nu)}(x)$.
Fix $\zeta\ne 1$.
It is sufficient to show that in $\Z/m$ we have
$$
 \sum_{x\in T} f^{\alpha,(\epsilon,\nu)}(x)f^{(\epsilon)}(\zeta\alpha^{-1} x)
 =
 0.
$$
By definition we have
$$
 f^{\alpha,(\epsilon,\nu)}(x)
 =
 \ord_x(\path(\alpha,(\epsilon,\nu))\bowtie\scF)
$$
It is therefore sufficient to show that for any $d$-cell $\scP$
 in $\gX$ we have:
$$
 \sum_{x\in T}
 \ord_x(\path(\alpha,(\epsilon,\nu))\bowtie\scP)
 f^{(\epsilon)}(\zeta\alpha^{-1} x)
 =
 0.
$$
Fix such a $\scP$.
We cut $\path(\alpha,(\epsilon,\nu))\bowtie\scP$
 into $2^{d}$ smaller pieces.
Define for each subset $A\subseteq\{1,2,\ldots,d\}$
$$
 \scA_A(x_{1},\ldots,x_{d})
 =
 \sum_{i\notin A} \alpha\path(\epsilon,x_{i})\cdot a_{\scP,i}
 +
 \sum_{i\in A} (\alpha+(1-\alpha)x_{i} + \phi(m^{2}x_{i})\nu) a_{\scP,i}.
$$
Here $\phi$ is as in \S5.3.
Then for $x\in T$ we have
$$
 \ord_x(\path(\alpha,(\epsilon,\nu))\bowtie\scP)
 =
 \sum_{A\subseteq \{1,\ldots,d\}}
 \ord_x(\scA_A).
$$
It is therefore sufficient to prove that for any $A$ we have
 in $\Z/m$:
$$
 \sum_{x\in T}
 \ord_x(\scA_A)
 f^{(\epsilon)}(\zeta\alpha^{-1} x)
 =
 0.
$$

There are two cases which we must consider.
In the case that $A$ is empty, we have
$$
 \alpha^{-1}\scA_A
 =
 \path(\epsilon)\bowtie\scP.
$$
Therefore if $x$ is in the base set of $\scA_{A}$,
 then $\alpha^{-1}x$ (which is well defined as $\alpha\in\Upsilon$)
 is in the base set of $\path(\epsilon)\bowtie\scP$.
If this is the case then as $\zeta\ne 1$, we have
 $f^{(\epsilon)}(\alpha^{-1}\zeta x)=0$.

Next suppose that there is an index $j\in A$.
Without loss of generality assume $1\in A$.
Then we have a decomposition $\scA_A = \scB\Diamond \scC$,
 where $\scB:I\to X_{\infty}$ is the 1-cell given by
$$
 \scB(x)
 =
 ((1-\alpha)x + \phi(m^{2}x)\nu)\cdot a_{\scP,1},
$$
 and $\scC$ is a $d-1$-cell.
As the function $\phi$ is periodic,
 we may write $\scB$ as a sum of $m^{2}$ translations of
$$
 \scB_{1}(x)
 =
 \left(\frac{1-\alpha}{m^{2}}x + \phi(x)\nu \right) \cdot a_{\scP,1},
$$
 where the translations are by multiples of
 $\frac{1-\alpha}{m^{2}} a_{\scP,1}$.
It follows from our congruence condition of $\alpha$
 that these translations are
 in $\alpha L/L$.

Our expression for $\scB$ implies a similar expression for
 $\scA_{A}$:
$$
 \scA_A
 =
 \sum_{l=1}^{m^{2}} \scS(l),
$$
 where each $\scS(l)$ is a translation of $\scS(1)$ by
 an element of $\alpha L/L$.
As both the set $T$ and the function $f^{(\epsilon)}(\zeta\alpha^{-1}x)$
 are invariant under such translations, we have
\begin{eqnarray*}
 \sum_{x\in T}\ord_{x}(\scA_{A})f(\zeta\alpha^{-1}x)
 &=&
 \sum_{l=1}^{m^{2}}
 \sum_{x\in T}\ord_{x}(\scS(l))f(\zeta\alpha^{-1} x)\\
 &=&
 m^{2}
 \sum_{x\in T}\ord_{x}(\scS(1))f(\zeta\alpha^{-1} x)
 \equiv
 0 \bmod m.
\end{eqnarray*}
\sq

\subsection{Proof of Lemma \ref{C}.}

Choose a lattice $M\subset V_{m}$ containing
 $L$ and $\alpha^{-1} L$.
It is sufficient to show that for
 every $d-1$-cell $\scP$ in $\gX$ we have in $\Z/m$:
$$
 \{
 \scH^\alpha_{\alpha\beta,\alpha}\bowtie \scP
 |
 \alpha M
 \}
 =
 \{
 \scH^1_\beta\bowtie \scP
 |
 M
 \}.
$$
To prove this, we cut $\scH^\alpha_{\alpha\beta,\alpha}\bowtie \scP$
 into $2^{d-1}$ pieces.
One of the pieces will be precisely $\alpha\cdot\scH^1_\beta\bowtie \scP$;
 for the other pieces $\scA$ we will show that $\{\scA|M\}=0$.

The various pieces of $\scH^\alpha_{\alpha\beta,\alpha}\bowtie \scP$
 will be indexed by the subsets $A\subseteq\{1,2,\ldots,d-1\}$.
Recall that $\scH^\alpha_{\alpha\beta,\alpha}\bowtie\scP$ is defined by
$$
 (\scH^\alpha_{\alpha\beta,\alpha}\bowtie\scP)
 (t,x_1,\ldots,x_{d-1})
 =
 v_\scP
 +
 \sum_{i=1}^{d-1}
 \scH^\alpha_{\alpha\beta,\alpha}(t,x_i)\cdot a_{\scP,i}.
$$
For any subset $A\subset\{1,2,\ldots,d-1\}$ we define
$$
 \scA_A
 (t,x_1,\ldots,x_{d-1})
 =
 v_\scP
 +
 \sum_{i\in A}
 \scH^\alpha_{\alpha\beta,\alpha}(t,x_i/2)\cdot a_{\scP,i}
 +
 \sum_{i\not\in A}
 \scH^\alpha_{\alpha\beta,\alpha}(t,(x_i+1)/2)\cdot a_{\scP,i}.
$$
As no element of $\alpha M$ lies in any set $|\partial\scA_A|$,
 we have:
$$
 \{\scH^\alpha_{\alpha\beta,\alpha}\bowtie\scP|\alpha M\}
 =
 \sum_A \{\scA_A|\alpha M\}.
$$

We now examine the pieces $\scA_A$.
First consider $\scA_{\{1,2,\ldots,d-1\}}$.
For $t,x\in I$ we have
$$
 \alpha^{-1}\cdot\scH^\alpha_{\alpha\beta,\alpha}(t,x/2)
 =
 \scH^1_\beta (t,x).
$$
This implies
$$
 \alpha^{-1}\cdot\scA_{\{1,2,\ldots,d-1\}}
 =
 \scH^1_\beta\bowtie\scP.
$$
Therefore
$$
 \{\scA_{\{1,2,\ldots,d-1\}}|\alpha M\}
 =
 \{\scH^1_\beta\bowtie\scP|M\}.
$$
The lemma will now follow when we show that for any
 proper subset $A\subset\{1,2,\ldots,d-1\}$
 we have $\{\scA_A|\alpha M\}=0$.

Fix a proper subset $A\subset \{1,2,\ldots,d-1\}$.
Without loss of generality we have $1\notin A$.
Note that for $t,x\in I$ we have
$$
 \scH^{\alpha}_{\alpha\beta,\alpha}(t,(x+1)/2)
 =
 \scM(x),\quad
 \hbox{where $\scM(x)=\path(\alpha,(x+1)/2)$.}
$$
We therefore have
$$
 \scA_A
 =
 (\scM \cdot a_{\scP,1}) \Diamond \scB,
$$
 for some $d-1$-cell $\scB$.
We shall now use this fact to cut $\scA_A$ into $m$
 pieces, which are translations of each other by elements of $M$.
We note that from the definition of $\path(\alpha)$ we have
$$
 \hbox{$\scM\left(x+\frac{1}{m}\right)
 =
 \scM(x) + \frac{1-\alpha}{m},\quad
 0\le x\le 1-\frac{1}{m}$.}
$$
We define for $x\in I$:
$$
 \scN(x)=\scM(\hbox{$\frac{x}{m}$})\cdot a_{\scP,1}.
$$
With this notation we have
$$
 \{\scA_A|\alpha M\}
 =
 \sum_{i=0}^{m-1}
 \left\{
 \left[i\frac{1-\alpha}{m}a_{\scP,1}\right]\Diamond\scN \Diamond \scB 
 \Big|
 \alpha M\right\}.
$$
To prove the lemma it only remains to show that all the
 terms in the above sum are equal.
We have
$$
 \left\{
 \left[i\frac{1-\alpha}{m}a_{\scP,1}\right]\Diamond\scN \Diamond \scB
 \Big|
 \alpha M\right\}
 =
 \left\{ \scN \Diamond \scB
 \Big|
 \alpha M-i\frac{1-\alpha}{m}a_{\scP,1}\right\},
$$
 so we need only show that the translation $\frac{1-\alpha}{m}a_{\scP,1}$
 is an element of the lattice $\alpha M$.
This follows from the congruence condition on $\alpha$, the condition on $M$
 and the fact that $(1-\rho)a_{\scP,1}\in L$.
\sq
\medskip

\subsection{Proof of Lemma \ref{D}.}

We must show that
 $\{\scH^{\alpha\beta,\alpha}_{\alpha\beta}\bowtie\scG|M\}=0$
 in $\Z/m$.
However by Corollary \ref{XVcor} it is sufficient to prove this
 formula with the homotopy $\scH^{\alpha\beta,\alpha}_{\alpha\beta}$
 replaced by another homotopy $\scH$ from $\path(\alpha\beta,\alpha)$
 to $\path(\alpha,\beta)$ as long as we have for $x$ close to $0$:
 $\scH(t,x)=2\alpha\beta x$.
We shall choose such a homotopy $\scH$ for which the calculation
 is easier.

\paragraph{The homotopy $\scH$.}
Before beginning we shall fix a parametrization of the
 path $\path(\alpha\beta,\alpha)$ as follows:
$$
 \path(\alpha\beta,\alpha)(x)
 =
 \left\{
 \begin{array}{ll}
     \alpha\beta\path(\epsilon,2x)
     &
     x\le \frac{1}{2},\\
     \alpha\path^{\beta}(0,\nu,2x-\frac{1}{2})
     &
     \frac{1}{2}\le x\le \frac{3}{4},\\
     \path^{\alpha}(0,\xi,2x-1)
     &
     x\ge \frac{3}{4}.
 \end{array}
 \right.
$$
This has the advantage that it agrees with $\path^{\alpha\beta}(x)$
 for $x\le \frac{1}{2}$.
We fix $\epsilon$, $\xi$ and $\nu$
 sufficiently small for out purposes,
 and we forget about these variables for the rest of the proof.

We introduce a sequence of
 auxiliary paths
$$
 \path(\alpha\beta,\alpha)=\path_{0},
 \ldots,
 \path_{m^{2}}=\path(\alpha\beta).
$$
These are defined as follows.
For any $i=0,\ldots,m^{2}$ we define
$$
 \path_{i}(x)=2x\alpha\beta, \quad x\le \frac{1}{2}.
$$
The triangle in $\gg$ with vertices $\alpha\beta$, $\alpha$ and $1$
 is cut into $m^{2}$ similar triangles, or $4m^{2}$ if $m$ is even.
These similar triangles are numbered as in the following diagram.
\medskip


\begin{center}
    \includegraphics{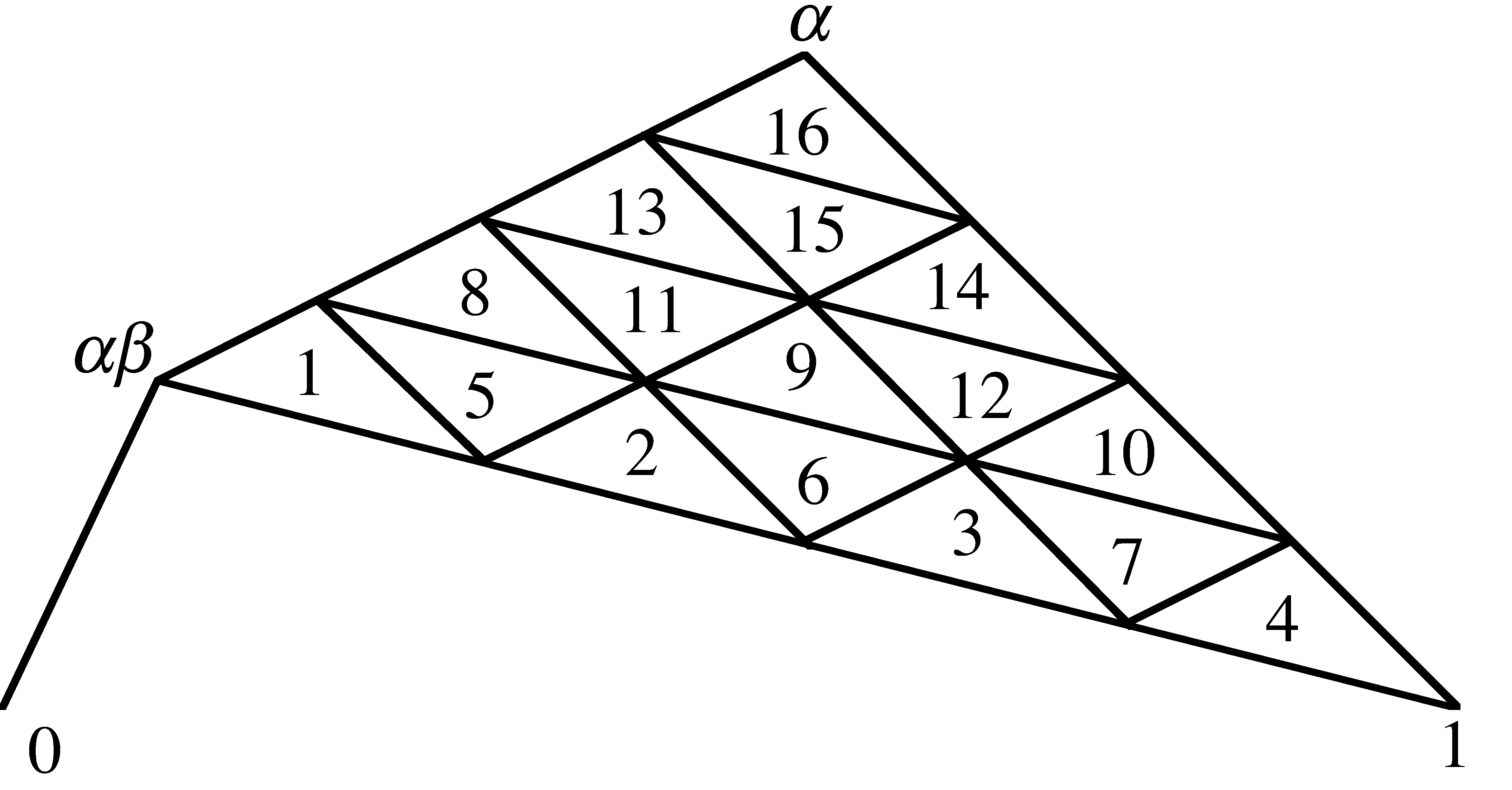}
\end{center}
\medskip

We define $\path_{i}$ to be the path from $0$ to $1$
 in the diagram which passes below the triangles numbered
 $1,\ldots,i$ and above the triangles numbered $i+1,\ldots,m^{2}$.
This means that paths $\path_{i-1}$ and $\path_{i}$ differ only
 by the $i$-th triangle.
We shall parametrize these paths as follows.
For any $x\in I$ which is not mapped into the edge of the
 $i$-th triangle we define $\path_{i-1}(x)=\path_{i}(x)$.
Let $[a_{i},b_{i}]\subset I$ be the subinterval mapped to the $i$-th triangle
 by both $\path_{i}$ and $\path_{i-1}$.
The interval $[a_{i},b_{i}]$ will have length $\frac{1}{2m}$.
One of the two paths will go around two edges of the triangle and
 the other will go around the third edge.
We parametrize the path which goes around only one of the edges
 so that the path is affine there.
The other path will be affine on the two
 subintervals $[a_{i},a_{i}+\frac{1}{4m}]$ and
 $[a_{i}+\frac{1}{4m},b_{i}]$, and will map $a_{i}+\frac{1}{4m}$ to the
 vertex of the triangle.

We define homotopies $\scH_{1},\ldots,\scH_{m^{2}}$ as follows:
$$
 \scH_{i}(x,t)
 =
 (1-t)\path_{i-1}(x) + t\path_{i}(x).
$$
Thus $\scH_{i}$ is a homotopy from $\path_{i-1}$ to $\path_{i}$.
Finally we put all these homotopies together to make the homotopy
 $\scH$:
$$
 \scH(x,t)
 =
 \scH_{[mt+1]}(x,\{mt\}),
$$
 where $[\cdot]$ and $\{\cdot\}$ denote the integer part and fractional part
 respectively.
This is a homotopy from $\path(\alpha\beta,\alpha)$
 to $\path(\alpha\beta)$.

\paragraph{Calculation of $\{\scH\bowtie\scG|M\}$.}
Note that we have
$$
 \{\scH\bowtie\scG|M\}
 =
 \blim_{\epsilon\to 0^{+}}
 \blim_{\nu\to 0^{+}}
 \blim_{\xi\to 0^{+}}
 \blim_{\eta\to 0^{+}}
 \{ \scH(\epsilon,\nu,\xi,\eta) \bowtie\scG|M\}.
$$
We fix $\epsilon,\nu,\xi,\eta$ so that all our functions are
 defined on $M$ and equal to their limits.

It is sufficient to show that for any $d-1$ cell $\scP$ in $\gX$
 we have in $\Z/m$:
$$
 \{\scH\bowtie\scP|M\}
 =
 0.
$$
We now fix such a $\scP$.

Assume for a moment that $m$ is odd.
Recall that to construct the homotopy $\scH$
 we used a sequence of paths
 $\path_{0},\ldots,\path_{m^{2}}$ and homotopies
 $\scH_{1},\ldots,\scH_{m^{2}}$,
 where $\scH_{i}$ is a homotopy from $\path_{i-1}$ to $\path_{i}$.
We therefore have
$$
 \{\scH|M \}
 =
 \sum_{i=1}^{m^{2}}
 \{\scH_i\bowtie\scP|M\}.
$$
Each homotopy $\scH_{i}$ corresponds to one of the $m^{2}$ subtriangles
 $T_{1},\ldots,T_{m^{2}}$ of the triangle with vertices $\alpha\beta,\alpha,1$.
Each of these triangles is either a translation of $T_{1}$
 or a translation of $T_{2}$.
We shall prove that if $T_{i}$ is a translation of $T_{j}$ then we
 have
\begin{equation}
    \label{claim2}
    \{\scH_i\bowtie\scP|M\}
    \equiv
    \{\scH_j\bowtie\scP|M\}
    \bmod m.
\end{equation}
As there are $\frac{m(m+1)}{2}$ triangles of type $T_{1}$ and
 $\frac{m(m-1)}{2}$ of type $T_{2}$, this implies
\begin{eqnarray*}
 \{\scH\bowtie\scP|M \}
 &\equiv&
 \frac{m(m+1)}{2}\{\scH_1\bowtie\scP|M\}
 +
 \frac{m(m-1)}{2}\{\scH_2\bowtie\scP|M\}\\
 &\equiv&
 0
 \qquad\bmod m,
\end{eqnarray*}
 which proves the result.
This is the only place in which we need to assume that $m$ is odd.
In the case that $m$ is even we must cut the large triangle into
 $4m^{2}$ subtriangles instead of just $m^{2}$
 so $m(2m+1)$ of them are of type $T_{1}$ and $m(2m-1)$ are
 of type $T_{2}$.
This is why we need a slightly different congruence condition
 on $\alpha$ and $\beta$ when $m$ is even.

It remains only to prove the congruence (\ref{claim2}).
To do this we shall cut $\scH_i\bowtie\scP$ and $\scH_j\bowtie\scP$
 into pieces.
Some of the pieces of $\scH_i\bowtie\scP$ will be translates by
 elements of $M$
 of pieces of $\scH_j\bowtie\scP$, and so will cancel each other out.
Any piece which does not cancel in this way
 will be the product of a line segment with length in $m\cdot M$
 and a $d-1$-chain.
Thus its contribution to (\ref{claim2}) will vanish modulo $m$.

Without loss of generality we shall assume that
 $T_{i}$ is a translation of $T_{1}$.
We begin by cutting the interval $I$ onto four pieces.
For $x$ in the interval $[0,\frac{1}{2}]$ we have
$$
 \scH_{i}(t,x)
 =
 \path_{i}(x)
 =
 \alpha\beta\path(\epsilon)(2x).
$$
Let $[a_{i},b_{i}]$ be the subinterval of $[\frac{1}{2},1]$ which is
 mapped by $\path_{i}$ and $\path_{i-1}$ to the triangle $T_{i}$.
The other two pieces of $I$ are $[\frac{1}{2},a_{i}]$ and $[b_{i},1]$.
It is possible that one of these two will be a single point.
For $x$ in the interval $[\frac{1}{2},a_{i}]$ we have
 $\scH_{i}(x,t)=\path_{i}(x)$.
In this region, $\path_{i}$ is a sum of line segments
 whose endpoints differ by $\frac{1-\alpha}{m}$,
 $\frac{1-\alpha\beta}{m}$ or $\frac{\alpha-\alpha\beta}{m}$.
We define
$$
 \scU_{i}(x)
 =
 \path_{i}\left((1-x)\frac{1}{2}+xa_{i}\right).
$$
The region of $\path_{i}$ between $b_{i}$ and $1$ is similar
 and we define
$$
 \scU_{i}^{\prime}(x)
 =
 \path_{i}((1-x)b_{i}+x).
$$

Suppose $\{1,\ldots,d-1\}$ is the disjoint union of the
 four sets $A$, $B$, $C$ and $D$.
We shall define $\scA_{i}(A,B,C,D)$ to be the restriction
 of $\scH_{i}\bowtie\scP$ to the subset
$$
 I
 \times
 [0,\textstyle{\frac{1}{2}}]^{A}
 \times
 [\textstyle{\frac{1}{2}},a_{i}]^{B}
 \times
 [a_{i},b_{i}]^{C}
 \times
 [b_{i},1]^{D}
 \subset I^{d}.
$$
In other words we have
\begin{eqnarray*}
 \scA_{i}(A,B,C,D)(t,x_{1},\ldots,x_{d-1})
 &=&
 v_\scP
 +
 \sum_{j\in A}
 \alpha\beta\path(\epsilon)(x_{j})\cdot a_{\scP,j}\\
 &&
 +
 \sum_{j\in B}
 \scU_{i}(x_{j})\cdot a_{\scP,j}\\
 &&
 +
 \sum_{j\in C}
 \scH_{i}((1-x_{j})a_{i}+x_{j}b_{i},t)\cdot a_{\scP,j}\\
 &&
 +
 \sum_{j\in D}
 \scU_{i}^{\prime}(x_{j})\cdot a_{\scP,j}.
\end{eqnarray*}
We have
$$
 \{\scH_{i}\bowtie\scP|M\}
 =
 \sum_{A,B,C,D}
 \{\scA_{i}(A,B,C,D)|M\}.
$$
It is sufficient to prove that for any choice of $A$, $B$, $C$ and $D$
 we have
$$
 \{\scA_{i}(A,B,C,D)|M\}
 =
 \{\scA_{1}(A,B,C,D)|M\}.
$$
To prove this we shall consider three cases.

Case 1. Suppose that $B$ is non-empty
 and let $j\in B$.
We can then decompose $\scA_{i}(A,B,C,D)$ as
$$
 \scA_{i}(A,B,C,D)
 =
 \scV \Diamond \scW,
$$
 where $\scV(x)=\scU_{i}(x)\cdot a_{\scF,j}$
 and $\scW$ is a $d-1$-chain.
The cell $\scV$ is a sum of line segments
 whose length is in $m\cdot M$.
It follows that $\{\scA_{i}(A,B,C,D)|M\}\equiv 0\bmod m$.
Similarly $\{\scA_{1}(A,B,C,D)|M\}\equiv 0\bmod m$.

Case 2.
Suppose that $D$ is non-empty.
We may reason as in case 1 to show that
 both $\{\scA_{i}(A,B,C,D)|M\}$ and $\{\scA_{1}(A,B,C,D)|M\}$
 are congruent to 0 modulo $m$.

Case 3.
Suppose $B$ and $D$ are empty.
We shall prove that $\scA_{i}(A,B,C,D)$ is a translation
 of $\scA_{1}(A,B,C,D)$ by an element of $M$.
We shall assume without
 loss of generality that $C=\{1,\ldots,r\}$
 and $A=\{r+1,\ldots d-1\}$.
We may then decompose $\scA_{i}(A,B,C,D)$ as follows:
$$
 \scA_{i}(A,B,C,D)
 =
 \scB_{i}
 \Diamond
 \scC,
$$
 where $\scB_{i}:I^{r+1}\to X_{\infty}$ is given by
$$
 \scB_{i}(t,x_{1},\ldots,x_{r})
 =
 \sum_{j=1}^{r}
 \scH_{i}(t,(1-x_{j})a_{i}+x_{j}b_{i})\cdot a_{\scP,j}
$$
 and $\scC:I^{d-r-1}\to X_{\infty}$ is given by:
$$
 \scC(\xx)
 =
 \sum_{j=1}^{d-r-1}
 \alpha\beta\path(\epsilon, x_{j})\cdot a_{\scP,j}.
$$
It is therefore sufficient to prove that
 $\scB_{i}$ is a translation of $\scB_{1}$
 by an element of $M$.
Recall that the triangle $T_{i}$ is a translation
 of $T_{1}$ by some vector $v\in\gg$.
Furthermore $v$ is of the form
$$
 v
 =
 r\frac{\alpha\beta-1}{m}
 +
 s\frac{\alpha-1}{m},
 \qquad
 r,s\in\Z.
$$
It follows that the restriction of $\scH_{i}$ to $[a_{i},b_{i}]$
 and the restriction of $\scH_{1}$ to $[a_{1},b_{1}]$ also
 differ by the translation $v$.
In other words we have for $x\in I$,
$$
 \scH_{i}(t,(1-x)a_{i}+x b_{i})
 =
 \scH_{i}(t,(1-x)a_{1}+x b_{1})+v.
$$
This implies
$$
 \scB_{i}(t,x_{1},\ldots,x_{r})
 =
 \scB_{1}(t,x_{1},\ldots,x_{r})
 +
 v\sum_{j=1}^{r}a_{\scP,j}.
$$
It remains to show that the translations $v\cdot a_{\scP,j}$ are in $M$.
This reduces to showing that $\frac{\alpha-1}{m}a_{\scP,j}$
 and $\frac{\alpha\beta-1}{m}a_{\scP,j}$ are in $M$.
However this is true by our congruence conditions on $\alpha$
 and $\beta$, the fact that $(1-\rho)a_{\scP,j}\in L$
 and the assumption that $M$ contains $L$, $\alpha L$ and $\alpha\beta L$.
\sq

\section{Extending the cocycle to $\GL_n(\A)$.}

\subsection{The cocycle on $\SL_{n}(k_{m})$.}

Let $k_{m}$ be the sum of the fields $k_{v}$ for all
 finite places $v$ dividing $m$.
We therefore have $\A=\AS\oplus k_{\infty}\oplus k_{m}$.
So far, we have a cocycle $\decAS$ on $\GL_{n}(\AS)$
 and a cocycle $\deci $ on $\GL_{n}(k_{\infty})$.
We shall now find a continuous cocycle $\dec_{m}$
 on $\SL_{n}(k_{m})$ so that $\decAS\deci \dec_{m}$ is
 metaplectic on $\SL_{n}(\A)$.

\paragraph{Extending $\tau$.}
Recall that $G_{\gf}$ is the subgroup of $\GL_{n}(k)$ consisting of
 matrices which are integral at all places dividing $m$
 and congruent to the identity modulo $\gf$.
This is the subgroup generated by the semigroup $\Upsilon_{\gf}$.
We have a function $\tau:\Upsilon_{\gf}\to \mu_{m}$
 such that for $\alpha,\beta\in\Upsilon_{\gf}$ the following holds:
\begin{equation}
    \label{splitting}
    \decAS(\alpha,\beta)\deci (\alpha,\beta)
    =
    \partial\tau(\alpha,\beta).
\end{equation}
For any matrix $\alpha\in G_{\gf}$ there is a natural number $N$
 such that $\alpha/N\in\Upsilon_{\gf}$.
We extend $\tau$ to a function on $G_{\gf}$ by defining
$$
 \tau(\alpha)
 =
 \frac{\decAS(N,N^{-1}\alpha)\deci (N,N^{-1}\alpha)}
 {\decAS(N,N^{-1})\deci (N,N^{-1})}
 \tau(N^{-1}\alpha),
 \quad
 \alpha\in G_{\gf},
$$
where $N\in\N$ is chosen so that $N^{-1}\alpha\in\Upsilon_{\gf}$.
It follows from the cocycle relation and the fact that
 $\tau(N^{-1})=1$ for $N\in\N\cap\Upsilon^{-1}$,
 that this does not depend on our choice of $N$.
Furthermore it is easy to check that on the whole
 group $G_{\gf}$ the formula (\ref{splitting}) still holds.

We extend $\tau$ to a function on $\GL_{n}(k)$.
To do this we choose a set $\rep$ of representatives
 for cosets $\GL_{n}(k)/G_{\gf}$.
Thus every element of $\GL_{n}(k)$ may be uniquely
 expressed in the form $r\alpha$ with $r\in\rep$ and $\alpha\in G_{\gf}$.
We define
$$
 \tau(r\alpha)
 =
 \decAS(r,\alpha)\deci (r,\alpha)\tau(\alpha).
$$

\paragraph{The cocycle $\dec_{m}$.}
We define the cocycle $\dec_{m}$ on the dense subgroup $\SL_{n}(k)$
 of $\SL_{n}(k_{m})$ by
\begin{equation}
    \label{decm}
    \dec_{m}(\alpha,\beta)
    =
    \frac{\partial\tau(\alpha,\beta)}
    {\decAS(\alpha,\beta)\deci (\alpha,\beta)}.
\end{equation}
We shall prove that $\dec_{m}$ extends to a continuous cocycle
 on $\SL_{n}(k_{m})$.
We then define for $\alpha,\beta\in\SL_{n}(\A)$
$$
 \decA(\alpha,\beta)
 =
 \decAS(\alpha,\beta)
 \deci (\alpha,\beta)
 \dec_{m}(\alpha,\beta).
$$
It follows immediately from the definition (\ref{decm})
 that the restriction of $\decA$ to $\SL_{n}(k)$
 is $\partial \tau$.
Therefore $\decA$ is metaplectic.
It remains to prove the following.

\begin{theorem}
    The cocycle $\dec_{m}$ on $\SL_{n}(k)$ extends by
    continuity to $\SL_{n}(k_{m})$.
\end{theorem}

\proof
It follows immediately from the definitions
 that for $\alpha,\beta\in \SL_{n}(k)$ and $\epsilon\in G_{\gf}$ we have
$$
 \dec_{m}(\alpha,\beta\epsilon)
 =
 \dec_{m}(\alpha,\beta).
$$
It is therefore sufficient to prove that for any $\beta\in\SL_{n}(k)$,
 the function $\alpha\mapsto\dec_{m}(\alpha,\beta)$ is continuous.
We fix $\beta$.
Note that for $\epsilon\in G_{\gf}\cap (\beta G_{\gf}\beta^{-1})$ we have
 from the cocycle relation:
$$
 \dec_{m}(\alpha\epsilon,\beta)
 =
 \dec_{m}(\alpha,\beta)\dec_{m}(\epsilon,\beta).
$$
In particular the map $\psi:G_{\gf}\cap \beta G_{\gf}\beta^{-1}\to\mu_{m}$
 given by $\psi(\epsilon)=\dec_{m}(\epsilon,\beta)$ is a homomorphism.
To prove continuity we need only show that $\ker(\psi)\cap \SL_{n}(k)$
 is open in the induced topology from $\SL_{n}(k_{m})$.
However this fact follows immediately from
 Lemma \ref{congruence} below.
\sq
\medskip

\begin{remark}
    It is worth noting that $\ker(\psi)$ is not open in the topology
     induced by $\GL_{n}(k_{m})$ and $\dec_{m}$ does not extend by
     continuity to $\GL_{n}(k_{m})$.
\end{remark}

\paragraph{A weak form of the congruence subgroup problem.}
Let $S_{m}$ be the set of finite primes in $S$ and define
 $R=\{x\in k: \forall v\in S_{m}\; |x|_{v}\le 1\}$.
This is a subring of $k$ whose primes correspond to the elements
 of $S_{m}$.
The ring $R$ is a dense subring of $\go_{m}:=\oplus_{v\in S_{m}}\go_{v}$.

\begin{lemma}
    \label{congruence}
    Every subgroup $H$ of finite index in $\SL_{n}(R)$
     is a congruence subgroup, ie. $H$ contains all matrices
     congruent to the identity modulo some non-zero ideal of $R$.
    Equivalently $\SL_{n}(\go_{m})$ is the profinite completion
     of $\SL_{n}(R)$.
\end{lemma}

\begin{remark}
    This is a very weak statement, which follows immediately
    from \cite{bassmilnorserre} (or \cite{serre2} for $n=2$).
    However since these papers effectively construct the
    universal metaplectic cover of $\SL_{n}$,
    it is worth noting that the limited result required here
    can be obtained in an elementary way.
\end{remark}

\proof
Let $H$ be a subgroup of $\SL_{n}(R)$ of index $d$.
We shall assume without loss of generality that $m$ divides $d$
 and that $H$ is normal in $\SL_{n}(R)$.
We shall show that any matrix congruent to the identity
 modulo $d^{2}R$ is in $H$.
First note that the $d$-th power of any element of $\SL_{n}(R)$
 is in $H$.
Therefore $H$ contains the elementary matrices
$$
 I_{n}+\lambda d e_{i,j},
 \quad
 \lambda\in R,\;
 i\ne j.
$$
Here $e_{i,j}$ denotes the matrix whose $(i,j)$-entry is 1
 and whose other entries are all zero.

By a $d$-operation we shall mean an operation of the form ``add
 $\lambda d$ times row $i$ to row $j$'' ($i\ne j$, $\lambda \in R$).
If a matrix can be reduced by $d$-operations to the identity matrix
 then that matrix must be in $H$ since $d$-operations have the effect of
 multiplying on the left by $I_{n}+\lambda d e_{i,j}$.

Now let $A=(a_{i,j})$ be any matrix congruent to the identity
 modulo $d^{2}R$.
We shall show that $A$ may be reduced to the identity matrix
 by $d$-operations.
Since $m$ divides $d^{2}$, the entries $a_{i,i}$ on the diagonal of $A$
 are units in $R$ (they are congruent to 1 modulo every prime ideal
 of $R$).
The entries off the diagonal are divisible by $d^{2}$.
Therefore we can reduce $A$ by $d^{2}$-operations to a diagonal matrix.
Furthermore the diagonal matrix which we obtain
 will still be congruent to the identity modulo $d^{2}$.

Now let $A$ be diagonal.
It remains to show how $A$ can be reduced to the identity.
We first describe a method for converting $a_{i,i}$ to $1$.
For $i<n$ we may add $da_{i,i}^{-1}$ times row $i$ to row $i+1$.
This gives us a $d$ in the $(i+1,i)$ entry.
Next subtracting $d\frac{a_{i,i}-1}{d^{2}}$ times row $i+1$ from
 row $i$ we obtain a $1$ in the $(i,i)$ entry.
After this we subtract $d$ times row $i$ from row $i+1$ to obtain
 a zero there.
Finally, subtracting a multiple of row $i+1$ from row $i$
 we obtain a diagonal matrix with a 1 in the $(i,i)$ position.
In this process we have only changed the
 $(i,i)$ and $(i+1,i+1)$-entries.
We may perform this process for $i=1,2,\ldots,n-1$
 consecutively to obtain a diagonal matrix with $a_{i,i}=1$ for
 $i=1,2,\ldots,n-1$.
Since the resulting matrix has determinant 1, it follows that
 $a_{n,n}$ is also 1.
\sq

\paragraph{A product formula for $\decA$.}
For any place $v$ of $k$, let $\decv$ be the restriction of
 $\decA$ to $\SL_{n}(k_{v})$.
It is known (see \cite{hill2,hill3}) that
 for almost all places $v$
 the cocycle $\decv$ is trivial on $\GL_{n}(\go_{v})$.
Therefore the product $\prod_{v}\dec_{v}(\alpha,\beta)$
 converges and we have
 up to a coboundary\footnote{A formula for this coboundary
 is given in \cite{hill2}.}
$$
 \decAS(\alpha,\beta)
 =
 \prod_{v\notin S}
 \decv(\alpha,\beta),
 \quad
 \alpha\beta\in\SL_{n}(\AS).
$$
As the groups $\SL_{n}(k_{v})$ ($v\in S$) are perfect,
 the restriction maps give an isomorphism:
$$
 H^{2}(\SL_{n}(k_{S}),\mu_{m})
 \cong
 \bigoplus_{v\in S}
 H^{2}(\SL_{n}(k_{v}),\mu_{m}).
$$
(The corresponding statement for $\GL_{n}$ would be false).
This implies up to a coboundary on $\SL_{n}(\A)$:
$$
 \decA
 =
 \prod_{v}
 \decv.
$$

\subsection{Application : The power reciprocity law.}

We shall now deduce the power reciprocity law
 from our results.

%
%
%
%

Consider the cocycle $\decA$ on $\SL_{3}(\A)$.
For $\alpha\in \A^{\times}$ we define
$$
 \varphi(\alpha)
 =
 \left(
 \begin{array}{ccc}
     \alpha & 0 & 0\\
     0 & \alpha^{-1} & 0\\
     0 & 0 & 1
 \end{array}
 \right),
 \quad
 \varphi^{\prime}(\alpha)
 =
 \left(
 \begin{array}{ccc}
     \alpha & 0 & 0\\
     0 & 1 & 0\\
     0 & 0 & \alpha^{-1}
 \end{array}
 \right),
$$
One knows (see \cite{hill2} Lemma 23)
 that for $\alpha,\beta\in\A^{\times}$ we have
$$
 [\varphi(\alpha),\varphi^{\prime}(\beta)]_{\decAS}
 =
 \prod_{v\notin S}
 (\alpha,\beta)_{v,m}.
$$
On the other hand if $\alpha,\beta\in k^{\times}$, then since
 $\decA$ is metaplectic we have
$$
 [\varphi(\alpha),\varphi^{\prime}(\beta)]_{\decA}
 =
 1.
$$
This implies for $\alpha,\beta\in k^{\times}$,
$$
 \prod_{v\not\in S}(\alpha,\beta)_{v,m}
 \prod_{v\in S}[\varphi(\alpha),\varphi^{\prime}(\beta)]_{\decv}
 =
 1.
$$
Fix a place $v\in S$ and consider the function
 $\psi:k_{v}^{\times}\times k_{v}^{\times}\to \mu_{m}$ defined by
$$
 \psi(\alpha,\beta)
 =
 [\varphi(\alpha),\varphi^{\prime}(\beta)]_{\decv},
$$
 where $\decv$ is the restriction of $\decA$ to $\SL_{3}(k_{v})$.
To prove the reciprocity law
 it remains to show that $\psi$ is the $m$-th power
 Hilbert symbol.

For real $v$ this is a consequence of Corollary 1 (\S3.8).
For compex $v$ both $\psi$ and the Hilbert symbol
 are trivial for topological reasons.
Assume from now on that $v$ is a non-archimedean
 place dividing $m$.
The function $\psi$ is bimultiplicative and continuous
 by the properties of commutators.
Furthermore we have $\psi(\alpha,1-\alpha)=1$ for all
 $\alpha\ne 1$ since this formula holds on the dense subset
 $k^{\times}\setminus\{0\}$.
Therefore $\psi$ is a continuous Steinberg symbol.
As the Hilbert symbol is the universal continuous
 Steinberg symbol on $k_{v}$ we have
 $\psi(\alpha,\beta)=(\alpha,\beta)_{m,v}^{a}$
 for some fixed $a\in\Z/m\Z$.
Substituting $\alpha=\zeta\in\mu_{m}$, $\beta\in\go_{v}^{\times}$ we have
 (see for example \cite{serre3} XIV Proposition 6):
$$
 (\zeta,\beta)_{m,v}
 =
 \zeta^{\textstyle{\frac{1-N^{k_{v}}_{\Q_{p}}(\beta)}{m}}}
$$
Taking $\beta\in\go$ close to 1 in the topology of $\go_{w}$
 for all $w\in S\setminus\{v\}$ we have
$$
 \psi(\zeta,\beta)
 =
 [\varphi(\zeta),\varphi'(\beta)]_{\decAS}^{-1}.
$$
Proposition 2 of \cite{hill2}
 implies
$$
 [\varphi(\zeta),\varphi'(\beta)]_{\decAS}^{-1}
 =
 \zeta^{\textstyle{\frac{1-N^{k}_{\Q}(\beta)}{m}}}
 =
 (\zeta,\beta)_{m,v}.
$$
With a suitable choice of $\zeta,\beta$ this shows that $a=1$.
\sq
\medskip

As was mentioned in the introduction,
 it would be more satisfactory to have a local definition
 of $\decv$ for $v\in S_{m}$
 and a local proof that $\psi$ is the Hilbert symbol.

\subsection{Extending the cocycle to $\GL_{n}$.}

We now have a metaplectic cocycle $\decA$ on $\SL_{n}(\A)$,
 whose restriction to $\SL_{n}(\AS\oplus k_{\infty})$
 extends naturally to $\GL_{n}(\AS\oplus k_{\infty})$.
One might ask whether $\decAS\deci$
 extends to a metaplectic cocycle on $\GL_{n}(\A)$.
The answer to this depends on precisely how one poses the question.
If one asks whether there is a cocycle $\dec_{m}$ on $\GL_{n}(k_{m})$,
 which when multiplied together with $\decAS$ and $\deci$
 gives a metaplectic cocycle then the answer is ``no''.
However it is true that $\decAS\deci$ is the restriction
 of a metaplectic cocycle $\decA$ on $\GL_{n}(\A)$.

\paragraph{The change of base field property.}
By embedding the group $\GL_{n}$ in $\SL_{n+1}$,
 one can obtain a perfectly good
 metaplectic extension of $\GL_{n}$ by restriction.
In fact, this is the metaplectic extension
 which has been most studied.
However, such extensions are badly behaved under change of
 base field (see \cite{paddytori})
 compared with $\decAS\deci$.
For this reason, I shall extend $\decAS\deci$ to $\GL_{n}(\A)$
 to obtain an extension which is well behaved under change of base
 field.

More precisely, if $l$ is a finite extension of $k$
 of degree $d$ then by choosing
 a basis for $l$ over $k$ we may regard $\GL_{n}(\A_{l})$ as a
 subgroup of $\GL_{nd}(\A_{k})$.
With this identification $\GL_{n}(l)$ is a subgroup of $\GL_{nd}(k)$.
Thus metaplectic extensions of $\GL_{nd}/k$ restrict to
 metaplectic extensions of $\GL_{n}/l$.
We shall write $R$ the restriction map from $\GL_{nd}/k$ to
 $\GL_{n}/l$.
The classes $\decAS$ and $\deci$
 have the following ``change of base field property'':
$$
 R(\decAS^{(k)})=\decAS^{(l)},\quad
 R(\deci ^{(k)})=\deci ^{(l)}.
$$
This is clear since the base field never arises
 in the definitions of $\decAS$ or $\deci $.
It is known that the class on $\GL_{n}$ obtained by
 restricting from $\SL_{n+1}$ does not have the
 change of base field property
 (see for example \cite{paddytori}).
We shall extend $\decAS\deci$ to $\GL_{n}(\A)$
 in such a way that it does have this property.
To achieve this it is clearly sufficient to treat
 the case $k=\Q(\mu_{m})$.
From now on we shall restrict ourselves to this case.

\paragraph{The metaplectic kernel of $\GL_{n}$.}
The group $\GL_{n}$ is the semi-direct product of
 $\SL_{n}$ and $\GL_{1}$.
The normal subgroups $\SL_{n}(\A)$
 and $\SL_{n}(k)$ are perfect,
 i.e. they are equal to their own commutator subgroups.
Therefore the Hochschild--Serre spectral sequence
 shows that the restriction maps give isomorphisms:
$$
 H^{2}(\GL_{n}(\A),\mu_{m})
 \cong
 H^{2}(\SL_{n}(\A),\mu_{m})
 \oplus
 H^{2}(\GL_{1}(\A),\mu_{m}),
$$
$$
 H^{2}(\GL_{n}(k),\mu_{m})
 \cong
 H^{2}(\SL_{n}(k),\mu_{m})
 \oplus
 H^{2}(\GL_{1}(k),\mu_{m}).
$$
For an algebraic group $G$ we shall write $\cM(G,\mu_{m})$
 for the kernel of the restriction map
 $H^{2}(G(\A),\mu_{m})\to H^{2}(G(k),\mu_{m})$.
The above isomorphisms imply that we have
$$
 \cM(\GL_{n},\mu_{m})
 \cong
 \cM(\SL_{n},\mu_{m})
 \oplus
 \cM(\GL_{1},\mu_{m}).
$$
We have already constructed a metaplectic extension
 of $\SL_{n}$, so to show that our cocycle extends to a metaplectic
 cocycle of $\GL_{n}$ we need only show that its restriction
 to $\GL_{1}(\AS\oplus k_{\infty})$ extends to a
 metaplectic cocycle on $\GL_{1}(\A)$.
However by Theorem \ref{stability}
 and \cite{hill2}, the restriction of $\decAS\deci$ to $\GL_{1}$ is simply
 the cocycle $\decAS\deci$ constructed in the case $n=1$.

We describe the group $H^{2}(\A^{\times},\mu_{m})$.
For $\sigma\in H^{2}(\A^{\times},\mu_{m})$
 the commutator of $\sigma$ is a continuous
 bimultiplicative, skew symmetric function
 $\A^{\times}\times\A^{\times}\to\mu_{m}$.
We therefore have a map
$$
 H^{2}(\A^{\times},\mu_{m})
 \to
 \Hom(\wedge^{2}\A^{\times},\mu_{m}).
$$
This map is surjective.
We write $H^{2}_{sym}(\A^{\times},\mu_{m})$ for its kernel.
There is an isomorphism given by the restriction maps (see \cite{klose}):
$$
 H^{2}_{sym}(\A^{\times},\mu_{m})
 \cong
 \bigoplus_{v} H^{2}(\mu_{m}(k_{v}),\mu_{m}).
$$
Each of the groups $H^{2}(\mu_{m}(k_{v}),\mu_{m})$ is
 canonically isomorphic to $\Z/m$.
We write $H^{2}_{asym}(\A^{\times},\mu_{m})$
 for the kernel of the restriction map
 $H^{2}(\A^{\times},\mu_{m})
 \to\oplus_{v} H^{2}(\mu_{m}(k_{v}),\mu_{m})$.
Thus the commutator gives an isomorphism:
$$
 H^{2}_{asym}(\A^{\times},\mu_{m})
 \cong
 \Hom(\wedge^{2}\A^{\times},\mu_{m}),
$$
 and we have a decomposition
$$
 H^{2}(\A^{\times},\mu_{m})
 =
 H^{2}_{sym}(\A^{\times},\mu_{m})
 \oplus
 H^{2}_{asym}(\A^{\times},\mu_{m}).
$$
There is a similar decomposition of $H^{2}(k^{\times},\mu_{m})$:
\begin{eqnarray*}
    H^{2}(k^{\times},\mu_{m})
    &=&
    H^{2}_{sym}(k^{\times},\mu_{m})
    \oplus
    H^{2}_{asym}(k^{\times},\mu_{m}),\\
    H^{2}_{sym}(k^{\times},\mu_{m})
    &\cong&
    H^{2}(\mu_{m}(k),\mu_{m})
    \cong
    \Z/m\Z,\\
    H^{2}_{asym}(k^{\times},\mu_{m})
    &\cong&
    \Hom(\wedge^{2}k^{\times},\mu_{m}).
\end{eqnarray*}
Consider the restriction map
 $H^{2}(\A^{\times},\mu_{m})\to H^{2}(k^{\times},\mu_{m})$.
Clearly the restriction to $k^{\times}$ of a symmetric cocycle on
 $\A^{\times}$ is symmetric\footnote{
 However the restriction of an asymmetric cocycle is not necessarily
 asymmetric.}.
The resulting map
 $H^{2}_{sym}(\A^{\times},\mu_{m})\to H^{2}(k^{\times},\mu_{m})$
 corresponds to the map $\oplus_{v}\Z/m\Z \to \Z/m\Z$ given by
$$
 (a_{v})
 \mapsto
 \sum_{v} a_{v}.
$$
We now examine the commutator of $\decAS\deci$.
Under our conditions on $k$ the group
 $\Hom(\wedge^{2}k_{\infty}^{\times},\mu_{m})$ is trivial,
 so $\deci $ has trivial commutator.
Therefore the commutator of $\deci \decAS$ is the
 same as the commutator of $\decAS$.
This has been calculated in \cite{hill2} (Theorems 4 and 5)
 and is given by
$$
 [\alpha,\beta]_{\AS}
 =
 \left\{
 \begin{array}{ll}
     (-1)^{\displaystyle{\frac{(|\alpha|_{\AS}-1)(|\beta|_{\AS}-1)}{m^{2}}}}
     \displaystyle{\prod_{v\notin S} (\alpha,\beta)_{v,m}}
     & \hbox{if $m$ is even,}\\
     \displaystyle{\prod_{v\notin S}(\alpha,\beta)_{v,m}}
     & \hbox{if $m$ is odd.}
 \end{array}
 \right.
$$
Define a map $\chi_{k}:\A_{k}^{\times}/k^{\times}\to\Z_{2}^{\times}$ by
$$
 \chi_{k}(\alpha)
 =
 \chi_{\Q}(N^{k}_{\Q}\alpha),
 \qquad
 \chi_{\Q}(\alpha)
 =
 \sign(\alpha_{\infty})\alpha_{2} \prod_{\hbox{$p$ finite}}|\alpha|_{p}.
$$
Define also a bilinear map
 $\psi_{k}:(\A^{\times}/k^{\times})\times(\A^{\times}/k^{\times})\to\mu_{m}$
 by
$$
 \psi_{k}(\alpha,\beta)
 =
 (-1)^{\frac{(\chi_{k}(\alpha)-1)(\chi_{k}(\beta)-1)}{m^{2}}}.
$$
The point of this is that for $\alpha,\beta\in\AS^{\times}$
 we have
$$
 \psi_{k}(\alpha,\beta)
 =
 (-1)^{\frac{(|\alpha|_{\AS}-1)(|\beta|_{\AS}-1)}{m^{2}}}.
$$

\begin{theorem}
    There is a unique class $\DecA\in \cM(\GL_{1}/\Q(\mu_{m}),\mu_{m})$
     with the following properties:
    \begin{itemize}
	\item
	the restriction of $\DecA$ to $\AS^{\times}$
	is $\decAS$;
	\item
	the restriction of $\DecA$ to $k_{\infty}^{\times}$
	is $\deci $;
	\item
	The commutator of $\DecA$ is
	$$
	 [\alpha,\beta]_{\DecA}
	 =
	 \psi_{\Q(\mu_{m})}(\alpha,\beta)
	 \prod_{v} (\alpha,\beta)_{m,v}.
	$$
    \end{itemize}
\end{theorem}

\proof
It follows from the above discussion that there is a unique
 class with the given commutator and restrictions
 to $\AS^{\times}$ and $k_{\infty}^{\times}$
 and any given symmetric restriction $\mu_{m}(k_{m})$.
As $m$ is a power of a prime, $k_{m}$ is a field,
 so there is a unique choice of restriction to $\mu_{m}(k_{m})$
 for which the restriction to $k$ is asymmetric.
\sq

%

\end{document}